\documentclass[11pt]{article}
\usepackage[T1]{fontenc}
\usepackage{latexsym,amssymb,amsmath,amsfonts,amsthm}
\usepackage{graphics}
\usepackage[demo]{graphicx}
\usepackage[section]{placeins}\usepackage{graphicx}
\usepackage{mathrsfs}
\usepackage{subfigure}
\usepackage{float}
\usepackage{color}

\usepackage{tikz}
\usetikzlibrary{calc}
\usetikzlibrary{intersections,through}

\usepackage{mathtools}

\usepackage{epstopdf}
\newcommand{\R}{{\mathbb R}}

\newcommand{\N}{{\mathbb N}}
\newcommand{\C}{{\mathbb C}}

\newcommand{\be}{\begin{eqnarray}}
\newcommand{\ben}{\begin{eqnarray*}}
\newcommand{\en}{\end{eqnarray}}
\newcommand{\enn}{\end{eqnarray*}}

\newcommand{\real}{{\rm Re\,}}
\newcommand{\ima}{{\rm Im\,}}

\newtheorem{theorem}{Theorem}[section]

\newtheorem{lemma}[theorem]{Lemma}

\newtheorem{definition}[theorem]{Definition}
\newtheorem{remark}[theorem]{Remark}

\definecolor{rot}{rgb}{0,0,0}
\definecolor{hw}{rgb}{0,0,0}

\begin{document}
\renewcommand{\theequation}{\arabic{section}.\arabic{equation}}

\title{Imaging a moving point source from multi-frequency data measured at one and sparse observation points (part II): near-field case in 3D}

	\author{Guanqiu Ma\footnotemark[2], Hongxia Guo\footnotemark[1]\; \footnotemark[2], Guanghui Hu\footnotemark[2]}
	
	\date{}
	\maketitle
	
	\renewcommand{\thefootnote}{\fnsymbol{footnote}}
	\footnotetext[1] {Corresponding author}
	\footnotetext[2]{School of Mathematical Sciences and LPMC, Nankai University, Tianjin, 300071, China( gqma@nankai.edu.cn, hxguo@nankai.edu.cn, ghhu@nankai.edu.cn).}
	\renewcommand{\thefootnote}{\arabic{footnote}}

	\begin{abstract}
In this paper, we introduce a frequency-domain approach to extract information on the trajectory of a moving point source. The method hinges on the analysis of multi-frequency near-field data recorded at one and sparse observation points in three dimensions. The radiating period of the moving point source is supposed to be supported on the real axis and a priori known. In contrast to inverse stationary source problems, one needs to classify observable and non-observable measurement positions. The analogue of these concepts in the far-field regime were firstly proposed in the authors' previous paper (SIAM J. Imag. Sci., 16 (2023): 1535-1571). In this paper we shall derive the observable and non-observable measurement positions for straight and circular motions in $\R^3$. In the near-field case, we verify that the smallest annular region centered at an observable position that contains the trajectory can be imaged for an admissible class of orbit functions.
Using the data from sparse observable positions, it is possible to reconstruct the $\Theta$-convex domain of the trajectory. Intensive 3D numerical tests with synthetic data are performed to show effectiveness and feasibility of this new algorithm.

\vspace{.2in} {\bf Keywords}: {\bf inverse moving source problem, Helmholtz equation, multi-frequency data, factorization method, uniqueness, near-field data}
	\end{abstract}

	\section{Introduction}

	\subsection{Time-dependent model and its inverse Fourier transform}

Assume the entire space $\mathbb{R}^3$ is filled by a homogeneous and isotropic medium. We designate the sound speed of the background medium as the constant $c > 0$. We consider the acoustic radiating problem incited by a moving point source. This source traces a trajectory defined by the $C^1$-smooth function $a(t): [t_{\min}, t_{\max}] \rightarrow \mathbb{R}^3$, with $0 < t_{\min} < t_{\max}$.
The source function $S(x,t)$ is supposed to radiate wave signal at the initial time point $t_{\min}$ and stop radiating at the ending time point $t_{\max}$, i.e., it is supported in the interval $[t_{\min}, t_{\max}] $ with respect to the time variable $t>0$.
More precisely, the source function is supposed to take the form
	\begin{equation}
		S(x,t) = \delta(x-a(t)) \ell(t) \chi(t),
	\end{equation}
	where $\delta$ denotes the Dirac delta function, $\ell(t)\in C(\R)$ is a real-valued function fulfilling the positivity constraint
	\begin{equation} \label{cd:l}
		|\ell(t)| \geq \ell_0 >0, \quad t\in [t_{\min},t_{\max}],
	\end{equation}
	and $\chi(t)$ is the characteristic function over the interval $[t_{\min}, t_{\max}]$, defined by
		\begin{equation*}
		\chi(t) := \left\{
		\begin{aligned}
			& 1, \quad t \in [t_{\min}, t_{\max}],\\
			& 0, \quad t \notin [t_{\min}, t_{\max}].
		\end{aligned}
		\right.
	\end{equation*}
	
Denote the trajectory by $\Gamma \coloneqq \left\{x: x=a(t),\, t\in [t_{\min}, t_{\max}]\right\} \subset \R^3$.
One can easily find Supp $S(\cdot, t) \subset \Gamma$ for all $t\in [t_{\min}, t_{\max}]$ in the distributional sense.
The propagation of the radiated wave fields $U(x,t)$ is governed by the initial value problem
	\begin{equation}
		\left\{
		\begin{aligned}
			&c^{-2}\frac{\partial^2 U}{\partial t^2} = \Delta U + S(x,t), \quad &&(x, t) \in \mathbb{R}^3 \times \mathbb{R}^+, \mathbb{R}^+ \coloneqq \{t\in \R: t>0\},\\
			&U(x,0)=\partial_t U(x,0) = 0, &&x\in \mathbb{R}^3.
		\end{aligned}
		\right.
	\end{equation}
The solution $U$ can be expressed through the convolution of the fundamental solution $G$ of the wave equation with the source term, i.e.,
	\begin{equation}\label{time-solu}
		U(x,t) = G(x;t) * S(x,t) \coloneqq \int_{\mathbb{R}^+} \int_{\mathbb{R}^3}
		G(x-y; t-\tau)S(y,\tau)\,dyd\tau
	\end{equation}
	where
	\begin{equation*}
		G(x;t) = \frac{\delta(t-c^{-1}|x|)}{4\pi |x|}.
	\end{equation*}
In this paper the one-dimensional Fourier and inverse Fourier transforms are defined by
	\ben
		(\mathcal{F}u)(\omega):=\frac{1}{\sqrt{2\pi}}\int_{\R}u(t)e^{-i\omega t}\,dt,\quad
		(\mathcal{F}^{-1}v)(t):=\frac{1}{\sqrt{2\pi}}\int_{\R}v(k)e^{i\omega t}\,d\omega,
	\enn
	respectively.
The inverse Fourier transform of $S$ is thus given by
 	\begin{equation}\label{sourcef}
 		f(x,\omega) \coloneqq (\mathcal{F}^{-1} S(x, \cdot))(\omega)=\frac{1}{\sqrt{2\pi}}\int_{\mathbb{R}} \delta(x-a(t))\ell(t)\chi(t) e^{i\omega t}\, dt = \frac{1}{\sqrt{2\pi}}\int_{t_{\min}}^{t_{\max}} \delta(x-a(t)) \ell(t) e^{i\omega t}\, dt.
 	\end{equation}
 	From the expression \eqref{time-solu}, one deduces the inverse Fourier transform of the wave field $U$,
 	\begin{equation}\label{wFU}
 		\begin{aligned}
 		u(x,\omega)=(\mathcal{F}^{-1}U)(x,\omega) &= \int_{\mathbb{R}^3} (\mathcal{F}^{-1}G)(x-y;\omega) (\mathcal{F}^{-1}S)(y,\omega)\,dy\\
 		&= \frac{1}{\sqrt{2\pi}}\int_{\mathbb{R}^3} \Phi(x-y;\omega /c) f(y,\omega)\,dy.
 		\end{aligned}
 	\end{equation}
 	Here, $\Phi(x;k)$ is the fundamental solution to the Helmholtz equation $(\Delta + k^2)w = 0$, given by
\begin{equation*}
		\Phi(x;k) = \frac{e^{ik|x|}}{4 \pi |x|},\quad x\in \mathbb{R}^3,\, |x| \neq 0.
	\end{equation*}
Taking the inverse Fourier transform on the wave equation yields the inhomogeneous Helmholtz equation
 	\begin{equation}\label{eq1}
 		\Delta u(x,\omega) + \frac{\omega^2}{c^2}  u(x,\omega) = -f(x,\omega), \qquad x\in \R^{3}, \;\omega>0.
 	\end{equation}
From \eqref{wFU} we observe that $u$ satisfies 	
 the Sommerfeld radiation condition
 	\be\label{SRC}
	\lim\limits_{r \to \infty} r (\partial_r u - i\frac{\omega}{c}u) = 0,\quad r = |x|,\en
which holds uniformly in all directions $x/|x|$.

\subsection{Formulation in the frequency domain and literature review}
	Denote by $[\omega_{\min}, \omega_{\max}]$ an interval of frequencies on the positive real axis.
From the time-domain settings we see
	$$ f(x,\omega) = 0, \,\mbox{for all } x\notin \Gamma, \omega\in [\omega_{\min},\omega_{\max}],$$
implying supp $f(\cdot, \omega) = \Gamma$ for all $\omega \in [\omega_{\min}, \omega_{\max}]$.
For every $\omega > 0$, the unique solution $u$ to \eqref{eq1}-\eqref{SRC} is given by \eqref{wFU}, i.e.,
	\begin{equation}\label{expression-w}
		u(x, \omega) = \frac{1}{\sqrt{2\pi}}\int_{\mathbb{R}^3} \Phi(x-y;\omega /c) f(y, \omega) dy=\frac{1}{8\pi^2} \int_{t_{\min}}^{t_{\max}} \frac{e^{i\omega(t+c^{-1}|x-a(t)|)}}{|x-a(t)|}\ell(t)\,dt,  \quad x \notin \Gamma.
	\end{equation}
Noting that the time-dependent source $S$ is real valued, we have $f(x, -\omega)= \overline{f(x,\omega)}$ for all $\omega>0$ and thus $u(x, -\omega)=\overline{u(x, \omega)}$.

\begin{figure}[!ht]
		\centering
		\scalebox{0.6}{
		\begin{tikzpicture}

		\draw [very thick,dotted] (5,0) arc [ start angle = 0, end angle = 360, radius = 5];
		\draw [very thick,smooth] (-1.5,-1.5) .. controls (1,-2) .. (2,0);
		\draw [very thick,dotted] (-5,0) .. controls (-0.5,-1.5) .. (5,0);
		\draw [very thick,dotted] (-5,0) .. controls (0.5,1.5) .. (5,0);
		\draw (0,5) node [above] {$|x^{(j)}|=R,\,u(x^{(j)},\omega)$};

		\draw (1,-2) node [below] {$\Gamma$};

		\draw (-1.5,-1.5) node [left] {$a(t_{\min})$};
		\draw (2,0) node [right] {$a(t_{\max})$};
		
		\fill (5,0) circle (3pt);
		\fill (0,5) circle (3pt);
		\fill (-5,0) circle (3pt);
		\fill (0,-5) circle (3pt);	
		\fill (-0.5,-1.1) circle (3pt);	
		\fill (0.5,1.1) circle (3pt);		
						
		\end{tikzpicture}
		}
		\caption{Imaging the trajectory $\Gamma$ from knowledge of multi-frequency near-field data measured at a finite number of observation points
		$|x^{(j)}|=R$, $j=1,2,\ldots,M$.}
		\label{ip}
	\end{figure}
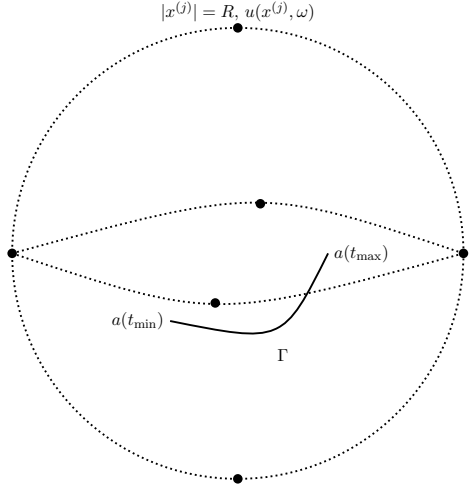

	In this paper we are interested in the following inverse problem (see Fig. \ref{ip}):
	\begin{description}
	\item[(IP):] Recovery the trajectory $\Gamma$ using the multi-frequency near-field data
	$$\{u(x^{(j)},\omega): \omega\in [\omega_{\min},\omega_{\max}],\,j=1,2,\cdots,M\}, \, x^{(j)}\in S_R:=\{x:|x|=R\},$$
where  $R>\sup\limits_{t\in [t_{\min},t_{\max}]} |a(t)|$.
\end{description}
A specific question of interest for (IP) is framed as follows:
\begin{description}
\item What kind information on $\Gamma$ can be extracted
from the the multi-frequency near-field data
	$\{u({x},\omega): \omega\in [\omega_{\min},\omega_{\max}]\}$ at a single observation point ${x}\in S_R$ ?
\end{description}	
	 The above questions are of great significance in various industrial, medical and military applications. This is mainly due to the fact that, in practical scenarios, the number of available measurement positions is inherently quite limited and the multi-frequency data can always be acquired  by inverse Fourier transforming the time-dependent signals.

To the best of the authors' knowledge, mathematical studies on direct and inverse scattering theory for moving targets are relatively scarce when compared to the extensive literature dedicated to scattering by stationary objects (see the monograph \cite{Isa1989}). Cooper \& Strauss \cite{Co79, CoStra} and Stefanov \cite{S1991} have made significant contributions to the rigorous mathematical theory of direct and inverse scattering from moving obstacles. 
Recently, there has been a growing research interest in detecting the motion of a moving point source governed by inhomogeneous wave equations. 
Several inversion algorithms have been proposed to recover the trajectory, profile and magnitude of a moving point source, for example, the algebraic method \cite{NIO2012, T2020, Ohe2011}, the time-reversal method \cite{GF2015}, the method of fundamental solutions \cite{CGMS2020}, matched-filter and correlation-based imaging schemes \cite{FGPT2017}, the iterative thresholding scheme \cite{Liu2021} and the method of Bayesian inference \cite{LGS2021, WKT2022}. In addition, the references \cite{LGY21, HKLZ2019, HKZ2020, HLY20, T2022} provide uniqueness and stability results on the identification of moving sources.

The focus of this paper is to establish a factorization method for imaging the trajectory $\Gamma$ using multi-frequency near-field data measured at sparse positions. The Factorization method, initially proposed by Kirsch in 1998 \cite{KG08}, has found successful applications in various inverse scattering problems involving multi-static data at a fixed energy (or equivalently, the Dirichlet-to-Neumann map). Its multi-frequency version was investigated in \cite{GS2017} and \cite{GGH2022} for inverse stationary source problems. Using the multi-frequency data at a single observation direction, one can reconstruct the smallest strip encompassing the support of the source and perpendicular to the observation direction. Moreover, the data from sparse near-field observations can be used to recover the $\Theta$-convex polygon (i.e., a convex polygonal region with normals aligned to observation directions) of the support. In our previous paper \cite{GHM2023} we have studied the same kind of inverse moving source problems by using the multi-frequency far-field data. The aim of this paper is to carry over the analysis and numerics of \cite{GHM2023} to the near-field case.

Similar to the discussions in \cite{GHM2023}, we will show  that imaging the smallest annular region of the motion from a single receiver is impossible for general orbit functions. This can be achieved only when the observation point is {\em observable} as defined in Definition \ref{obd} and when the orbit function possesses some monotonicity properties (see Theorem \ref{TH4.3} (ii)). In the absence of these conditions, one can solely obtain a slimmer annulus $A_{\Gamma}^{(x)}$ (to be defined  in \eqref{K}), with a width less than the aforementioned smallest annulus. For non-observable points, the test functions cannot lie in the range of the data-to-pattern operator as indicated in Lemma \ref{lem3.4}. Consequently, extracting any information on the motion of a moving source fails in principle, although the numerical reconstructions still display partial information. Utilizing sparse observable points, we propose an indicator function for imaging the analogue of the $\Theta$-convex domain related to the trajectory. Some uniqueness results are summarized in Theorem \ref{TH4.3} as a byproduct of the factorization scheme established in Theorems \ref{Th:factorization} and \ref{TH4.2}.

The remainder of the paper is structured as follows. Section \ref{sec:4} focuses on the factorization of the multi-frequency near-field operator $\mathcal{N}^{({x})}$ where $|x|=R$ is a fixed observation point. The factorization is based upon the data-to-pattern operator $\mathcal{L}^{(x)}$ and a middle operator of multiplication form, following the approach presented in \cite{GHM2023} and \cite{GS2017}. A range identity is presented to establish a connection between the ranges of $\mathcal{N}^{({x})}$ and $\mathcal{L}^{(x)}$.
Section \ref{RangeLx} is dedicated to the selection of appropriate test functions that effectively characterize the annulus $A_{\Gamma}^{(x)}$ through analysis of the range of the data-to-pattern operator $\mathcal{L}^{(x)}$.
In Section \ref{IdF}, we define indicator functions by using the near-field data measured at one or sparse observable positions.
Finally, Section \ref{num} presents numerical tests performed in three dimensions, validating the concepts discussed in the preceding sections.

	\section{Factorization of near-field operator}\label{sec:4}

The objective of this section is to develop a multi-frequency factorization method, employed to recover the trajectory $\Gamma = \text{Supp} f(\cdot, \omega)$ from the near-field data measured at the point ${x}\in S_R$. For this purpose, we will adopt the approach outlined in \cite{GGH2022} to derive a factorization of the near-field operator $\mathcal{N}^{({x})}$.
Motivated by \cite{GS2017}, we introduce two key parameters: the central frequency $\kappa$ and half of the bandwidth of the near-field data denoted as $K$.
	\begin{equation*}
		\kappa \coloneqq \frac{\omega_{\min} + \omega_{\max}}{2}, \quad K \coloneqq \frac{\omega_{\max} - \omega_{\min}}{2}.
	\end{equation*}
These notation enable us to define the linear near-field operator $\mathcal{N}^{({x})}:L^2(0,K) \to L^2(0,K)$ by
	\begin{equation}\label{FarO}
		(\mathcal{N}^{({x})}\phi)(\tau) \coloneqq \int_0^{K} u(x, \kappa + \tau - s)\,\phi(s)\,ds,\qquad \tau\in(0,K).
	\end{equation}
Recall from \eqref{expression-w} that $u$ is analytic in $\omega \in \mathbb{R}$. Hence the near-field operator $\mathcal{N}^{({x})}: L^2(0, K) \to L^2(0, K)$ is bounded. Further, it follows from \eqref{expression-w} that
 	\begin{equation} \label{def:F}
		(\mathcal{N}^{({x})}\phi)(\tau) = \int_0^{K} \int_{t_{\min}}^{t_{\max}} \frac{e^{i(\kappa + \tau - s)(t+c^{-1}|x-a(t)|)}}{8\pi^2 |x-a(t)|} \ell(t)\,dt \,\phi(s)\,ds.		
	\end{equation}
	Below we shall prove a factorization of the above near-field operator.

	\begin{theorem}
		We have $\mathcal{N}^{({x})} =  \mathcal{L}  \mathcal{T}  \mathcal{L}^*$ where $ \mathcal{L} = \mathcal{L}^{(x)} : L^2(t_{\min},t_{\max}) \to L^2(0, K)$ is defined by
		\begin{equation}\label{tildeL}
			( \mathcal{L}\psi)(\tau) \coloneqq \int_{t_{\min}}^{t_{\max}} e^{i\tau(t+c^{-1}|x-a(t)|)}\psi(t)\, dt, \quad \tau \in (0, K)
		\end{equation}
		for all $\psi \in L^2(t_{\min},t_{\max})$. Here the middle operator $ \mathcal{T} :  L^2(t_{\min},t_{\max}) \to  L^2(t_{\min},t_{\max})$ is a multiplication operator defined by
		\begin{equation}\label{oT}
			( \mathcal{T}\varphi)(t) \coloneqq \frac{e^{i\kappa (t+c^{-1}|x-a(t)|)}}{8\pi^2 |x-a(t)|} \ell(t)\varphi(t).
		\end{equation}
	\end{theorem}

	\begin{proof}
		We first show that the adjoint operator $ \mathcal{L}^* : L^2(0, K) \to L^2(t_{\min},t_{\max})$ of $ \mathcal{L}$ can be expressed by
		\begin{equation}\label{aj-tildeL}
			( \mathcal{L}^* \phi)(t) \coloneqq \int_{0}^{K} e^{-i s (t+c^{-1}|x-a(t)|)}\phi(s)\, ds, \quad \phi \in L^2(0, K).
		\end{equation}
		Indeed, for $\psi\in L^2(t_{\min},t_{\max})$ and $\phi \in L^2(0,K)$, it holds that
		\begin{eqnarray*}
			\langle  \mathcal{L}\psi, \phi \rangle_{L^2(0, K)} &=& \int_0^{K} \left(  \int_{t_{\min}}^{t_{\max}} e^{i\tau (t+c^{-1}|x-a(t)|)} \psi(t)\,dt\right) \,\overline{\phi(\tau)}\,d\tau \\
			&=& \int_{t_{\min}}^{t_{\max}} \psi(t) \left(\int_0^{K} \overline{e^{-i\tau(t+c^{-1}|x-a(t)|)}\phi(\tau)}d\tau\right)\,dt \\
			&=&\langle \psi,  \mathcal{L}^* \phi \rangle_{L^2(t_{\min},t_{\max})}.
		\end{eqnarray*}
		which implies \eqref{aj-tildeL}. By the definition of $\mathcal{T}$, we have
		\begin{equation*}
			( \mathcal{T}  \mathcal{L}^* \phi)(t) = \frac{e^{i\kappa (t+c^{-1}|x-a(t)|)}}{8\pi^2 |x-a(t)|} \ell(t) \int_{0}^{K} e^{-is(t+c^{-1}|x-a(t)|)} \phi(s) \,ds,\quad \phi \in L^2(0, K).
		\end{equation*}
		Hence, using \eqref{sourcef} and \eqref{def:F},
		\begin{eqnarray*}
			( \mathcal{L}  \mathcal{T}  \mathcal{L}^* \phi)(\tau) &=& \int_{t_{\min}}^{t_{\max}} e^{i\tau(t+c^{-1}|x-a(t)|)} \left(\frac{e^{i\kappa (t+c^{-1}|x-a(t)|)}}{8\pi^2 |x-a(t)|} \ell(t)  \int_{0}^{K} e^{-is(t+c^{-1}|x-a(t)|)} \phi(s) \,ds\right)dt\\
			&=& \int_0^{K}\int_{t_{\min}}^{t_{\max}} \frac{e^{i(\kappa + \tau - s)(t+c^{-1}|x-a(t)|)}}{8\pi^2 |x-a(t)|} \ell(t)\,dt\,\phi(s)\,ds\\
			&=& (\mathcal{N}^{({x})} \phi)(\tau).
		\end{eqnarray*}
 		This proves the factorization $\mathcal{N}^{({x})} =  \mathcal{L}  \mathcal{T}  \mathcal{L}^*$.
	\end{proof}

\begin{remark}
In the subsequent sections of this paper, we shall designate the operator $\mathcal{L}$ as the data-to-pattern operator associated with the orbit function $a(t)$. It is evident that the near-field data given by Equation \eqref{expression-w} can be represented as $u(x,\omega)=(\mathcal{L}\, \frac{\ell(t)}{8\pi^2 |x-a(t)|})(\omega)$.
\end{remark}

	Denote by $\text{Range} ( \mathcal{L})$ the range of the data-to-pattern operator $ \mathcal{L}=\mathcal{L}^{(x)}$ (see \eqref{tildeL}) acting on $L^2(t_{\min},t_{\max})$.

	\begin{lemma}
		The operator $ \mathcal{L} : L^2(t_{\min},t_{\max}) \to L^2(0, K)$ is compact with dense range.
	\end{lemma}

	\begin{proof}
For any $\psi \in L^2(t_{\min},t_{\max})$, it holds that $ \mathcal{L}\psi \in H^1(0, K)$, which is compactly embedded into $L^2(0, K)$. This proves the compactness of $ \mathcal{L}$. By \eqref{aj-tildeL},
$(\mathcal{L}^* \phi)(t)$ coincides with the Fourier transform of $\phi$ at $t+c^{-1}|x-a(t)|$.
If the set $\{t+c^{-1}|x-a(t)|: t\in [t_{\min}, t_{\max}]\}$ forms an interval of $\R$, the relation $(\mathcal{L}^* \phi)(t)=0$ implies $\phi=0$ in $L^2(0, K)$ with the properties of Fourier transform. When $t+c^{-1}|x-a(t)|$ equals a constant $C$ as $t\in [t_{\min}, t_{\max}]$, $(\mathcal{L}^* \phi)(t)=0$ also implies $\phi=0$ in $L^2(0, K)$. Hence, $ \mathcal{L}^*$ is injective. The denseness of $\text{Range} ( \mathcal{L})$ in $L^2(0, K)$ follows from the injectivity of $\mathcal{L}^*$.
	\end{proof}

Within the framework of factorization method, it is essential to connect the ranges of $\mathcal{N}^{({x})}$ and $\mathcal{L}$. 
	We first recall that, for a bounded operator $F: Y\rightarrow Y$ in a Hilbert space $Y$ the real and imaginary parts of $F$ are defined respectively by
	\begin{equation*}
		\real F=\frac{F+F^*}{2},\quad \ima F=\frac{F-F^*}{2i},
	\end{equation*}
	which are both self-adjoint operators. Furthermore, by spectral representation we define the self-adjoint and positive operator $|\real F|$ as
	\begin{equation*}
		|\real F|=\int_{\R} |\lambda|\, d E_\lambda,\qquad \mbox{if}\quad \real F=\int_{\R} \lambda\, d E_\lambda.
	\end{equation*}
Here $E_\lambda$ represents the projection measure. The self-adjoint and positive operator $|\ima F|$ can be defined analogously.
 Introduce a new operator
	\begin{equation*}
		F_{\#}:=|\real F| +|\ima F|.
	\end{equation*}
	Since $F_{\#}$ is selfadjoint and positive, its square root $F_{\#}^{1/2}$ is defined as
	\begin{equation*}
		F_{\#}^{1/2}:=\int_{\R^+} \sqrt{\lambda}\, d E_\lambda,\qquad \mbox{if}\quad  F_{\#}=\int_{\R^+} \lambda\, d E_\lambda.
	\end{equation*}
	In this paper we need the following result from functional analysis.
	\begin{theorem}(\cite{GGH2022}) \label{range}
		Let $X$ and $Y$ be Hilbert spaces  and let $F: Y\rightarrow Y$, $L: X\rightarrow Y$, $T: X\rightarrow X$ be linear bounded operators such that $F=LTL^*$. We make the following assumptions
		\begin{itemize}
			\item[(i)] $L$ is compact with dense range and thus $L^*$ is compact and one-to-one.
			\item[(ii)] $\real T$ and $\ima T$ are both one-to-one,  and the operator $T_{\#}=|\real T| +|\ima T|: X\rightarrow X$ is coercive, i.e., there exists $c>0$ with
				\begin{equation*}
					\big\langle T_{\#}\, \varphi, \varphi\big\rangle\geq c\,||\varphi||^2\quad\mbox{for all}\quad \varphi\in X.
				\end{equation*}
		\end{itemize}
		Then the operator $F_{\#}$ is positive and  the ranges of $F_{\#}^{1/2}:Y\rightarrow Y$ and  $L: X\rightarrow Y$ coincide.
	\end{theorem}

	To apply Theorem \ref{range} to our inverse problem, we set
	\begin{equation*}
		F=\mathcal{N}^{({x})}, \quad L=\mathcal{L}, \quad T=\mathcal{T}, \quad X=L^2(t_{\min},t_{\max}),\quad Y=L^2(0, K),
	\end{equation*}
	where $\mathcal{T}$ is the multiplication operator of \eqref{oT}. It is easy to see
	\begin{eqnarray*}
		\left[(\real \mathcal{T})\, \varphi\right](t) &=& \frac{\cos [\kappa (t+c^{-1}|x-a(t)|)]}{8\pi^2 |x-a(t)|} \ell(t) \varphi(t),\\
		\left[(\ima \mathcal{T})\, \varphi\right](t) &=& \frac{\sin [\kappa (t+c^{-1}|x-a(t)|)]}{8\pi^2 |x-a(t)|} \ell(t) \varphi(t)
	\end{eqnarray*}
 	are both one-to-one operators from $L^2(t_{\min},t_{\max})$ onto $L^2(t_{\min},t_{\max})$. The coercivity assumption of $\mathcal{N}^{({x})}$ yields the coercivity of  $\mathcal{T}_{\#}$. As a consequence of Theorem \ref{range}, we obtain
	\begin{equation}\label{RI}
		\mbox{Range}\, [(\mathcal{N}^{({x})})_{\#}^{1/2}]=\mbox{Range}\,(\mathcal{L}^{(x)})\quad\mbox{ for any }\, x\in S_R.
	\end{equation}

	Let $\varphi\in L^2(0, K)$ be a test function. We want to characterize the range of $\mathcal{L}^{(x)}$ through the choice of $\varphi$. Denote by $(\lambda_n^{(x)}, \psi_n^{(x)})$ an eigensystem of the positive and self-adjoint operator $ (\mathcal{N}^{({x})})_{\#}$, which is uniquely determined by the multi-frequency near-field data $\{u(x, \omega ) : \omega \in (\omega_{\min}, \omega_{\max})\}$. Applying Picard's theorem and Theorem \ref{range}, we obtain
	\begin{equation}
		\varphi \in \text{Range}( \mathcal{L}^{(x)}) \quad \text{if and only if } \quad \sum\limits_{n=1}^{\infty} \frac{|\langle \varphi, \psi_n^{(x)}\rangle|^2}{|\lambda_n^{(x)}|} < +\infty.
	\end{equation}
	To establish the factorization method, we now need to choose a proper class of test functions which usually rely on a sample variable in $\mathbb{R}^3$. 

	\section{Range of $ \mathcal{L}^{(x)}$ and test functions}\label{RangeLx}

To characterize the range of $ \mathcal{L}^{(x)}$, we need to investigate monotonicity of the function $ h(t):=t+c^{-1}|x-a(t)|\in C^1[t_{\min}, t_{\max}]$. To achieve this goal, we introduce the concept of division points for a continuous function defined over a closed interval.

\begin{definition}(\cite{GHM2023})\label{DIDP}
		Let $f\in C[t_{\min}, t_{\max}]$. The point $t\in (t_{\min}, t_{\max})$ is called a division point if \\
		(1) $f(t)=0$;\\
		(2) There exist an $\epsilon_0>0$ such that either $|f(t+\epsilon)|>0$ or $|f(t-\epsilon)|>0$ for all $0<\epsilon<\epsilon_0$.
	\end{definition}
Obviously, the division points constitute a subset of the zero set of a continuous function. However, a division point cannot be an interior point of the zero set.	
Given that $a(t)\in C^1[t_{\min}, t_{\max}]$, there exist a finite number of division points for the function $h'$. Let's denote these points as $t_1<t_2<\cdots<t_{n-1}$.
This process divides the interval $[t_{\min},t_{\max}]$ into $n$ sub-intervals, namely $[t_{j-1},t_j]$ for $j=1,2,\cdots,n$, where we set $t_{\min}=t_0$ and $t_{\max}=t_n$.
Let $a_j$ and $h_j$ be the restrictions of $a$ and $h$ to $[t_{j-1},t_j]$, respectively.
Let $a_j$ and $h_j$ represent the restrictions of functions $a$ and $h$, to the sub-interval $[t_{j-1},t_j]$ respectively. We then define
 $$\xi^{(x)}_{j,\min} := \inf\limits_{t\in [t_{j-1},t_j]} \{h_j(t)\},\quad \xi^{(x)}_{j,\max} := \sup\limits_{t\in [t_{j-1},t_j]} \{h_j(t)\},\quad j=1,2,\cdots n.$$
In each sub-interval $(t_{j-1}, t_{j})$, one of
 following cases must hold:
	\begin{itemize}
		\item $h'_j(t)>0$ for all $t\in (t_{j-1},t_j)$. There holds
		$$\xi^{(x)}_{j,\min} = t_{j-1}+c^{-1}|x-a_j(t_{j-1})|, \quad \xi^{(x)}_{j,\max} = t_j+c^{-1}|x-a_j(t_j)|;$$
		\item $h'_j(t)<0$ for all $t\in (t_{j-1},t_j)$. We have
		$$\xi^{(x)}_{j,\min} = t_j+c^{-1}|x-a_j(t_j)|, \, \quad\xi^{(x)}_{j,\max} = t_{j-1}+c^{-1}|x-a_j(t_{j-1})|;$$
		\item $h'_j(t)=0$ for all $t\in (t_{j-1},t_j)$. Consequently,
		$$\xi^{(x)}_{j,\min} = \xi^{(x)}_{j,\max} = t+c^{-1}|x-a_j(t)|,\quad t\in[t_{j-1},t_j].$$
	\end{itemize}
Define
	\be\label{def:xi}\xi^{(x)}_{\min} := \min\limits_j \xi^{(x)}_{j,\min}=
	\inf\limits_{t\in [t_{\min},t_{\max}]} \{h(t)\}
	,\quad \xi^{(x)}_{\max} := \max\limits_j \xi^{(x)}_{j,\max}
	=\sup\limits_{t\in [t_{\min},t_{\max}]} \{h(t)\},
	\en
which denote the minimum and maximum of $h$ over $[t_{\min}, t_{\max}]$, respectively.
	 If $|h'_j(t)|>0$, the monotonicity of the function $\xi=h_j(t)$ for $t\in[t_j, t_{j-1}]$ implies the inverse function $t=h_j^{-1}(\xi)\in C^1[\xi^{(x)}_{j,\min}, \xi^{x)}_{j,\max}]$.
Set $$J=\{j\in \N : 1\leq j\leq n, h_j^{\prime}(t)\equiv 0, t\in (t_{j-1},t_j)\}.$$ and assume $h_j(t)\equiv c_j\in \R$ for $j\in J$. Note that it is possible that $J=\emptyset$.

With these notations we can rephrase the operator $\mathcal{L}^{(x)}$ defined by \eqref{tildeL} as
	\begin{equation}\label{pic-L}
		\begin{aligned}
			(\mathcal{L}^{(x)} \psi)(\tau) &= \sum\limits_{j=1}^{n} \int_{t_{j-1}}^{t_j} e^{i\tau h_j(t)} \psi (t) \, dt \\
			&= \sum\limits_{j\notin J} \int_{t_{j-1}}^{t_j} e^{i\tau h_j(t)} \psi (t) \, dt + \sum\limits_{j\in J} e^{i\tau c_j} \int_{t_{j-1}}^{t_j}  \psi (t) \, dt.
		\end{aligned}
	\end{equation}
For $j\in J$, using $e^{i\tau c} = \sqrt{2\pi} \mathcal{F}^{-1}\delta(t-c)$ we can rewrite each term in the second sum as
	\begin{equation}\label{delta-L}
		e^{i\tau c_j} \int_{t_{j-1}}^{t_j}  \psi (t) \, dt = \sqrt{2\pi} \mathcal{F}^{-1}\delta(t-c_j) \int_{t_{j-1}}^{t_j}  \psi (t) \, dt.
	\end{equation}
For $j\notin J$ and $h'_j(t)>0$, the integral in the first summation on the right hand of \eqref{pic-L} takes the form
	\ben
		\int_{t_{j-1}}^{t_j} e^{i\tau h_j(t)}\psi(t)\,dt &=& \int_{\xi^{(x)}_{j,\min}}^{\xi^{(x)}_{j,\max}} e^{i\tau \xi} \psi(h_j^{-1}(\xi))\,(h_j^{-1}(\xi))^{\prime}\,d\xi\\
	&=& \int_{\xi^{(x)}_{j,\min}}^{\xi^{(x)}_{j,\max}} e^{i\tau \xi} \psi(h_j^{-1}(\xi))|(h_j^{-1}(\xi))^{\prime}|\,d\xi.	
	\enn
Note that $[h_j^{-1}(\xi)]^{\prime}>0$, due to the relation $h_j^{\prime}(t) [h_j^{-1}(\xi)]^{\prime} =1$. Analogously, if $h'_j(t)<0$ for some $j\notin J$, we have $[h_j^{-1}(\xi)]^{\prime}<0$ and thus	
\begin{equation*}
		\begin{aligned}
			\int_{t_{j-1}}^{t_j} e^{i\tau h_j(t)}\psi(t)\,dt &= -\int_{\xi^{(x)}_{j,\min}}^{\xi^{(x)}_{j,\max}} e^{i\tau \xi} \psi(h_j^{-1}(\xi))(h_j^{-1}(\xi))^{\prime}\,d\xi \\
			&= \int_{\xi^{(x)}_{j,\min}}^{\xi^{(x)}_{j,\max}} e^{i\tau \xi} \psi(h_j^{-1}(\xi))|(h_j^{-1}(\xi))^{\prime}|\,d\xi.
		\end{aligned}
	\end{equation*}
Now, extending $h_j^{-1}$ by zero from $(\xi^{(x)}_{j,\min},\xi^{(x)}_{j,\max})$ to $\mathbb{R}$ and extending $\psi\in L^2(t_{\min},t_{\max})$ by zero to $L^2(\mathbb{R})$, we can write each term for $j\notin J$ as

	\begin{equation}\label{con-L}
		\int_{t_{j-1}}^{t_j} e^{i\tau h_j(t)}\psi(t)\,dt = \int_{\mathbb{R}} e^{i\tau \xi} \psi(h_j^{-1}(\xi))|(h_j^{-1}(\xi))^{\prime}|\,d\xi.
	\end{equation}

	Combining \eqref{pic-L}, \eqref{delta-L} and \eqref{con-L}, we get
	\begin{equation}\label{Lxi}
			(\mathcal{L}^{(x)}\psi)(\tau) = \int_{\mathbb{R}} e^{i\tau\xi} g(\xi)\,d\xi,	
	\end{equation}
	with $$g(\xi)=\sum\limits_{j\notin J} \psi(h_j^{-1}(\xi))\, |(h_j^{-1}(\xi))'| + \sum\limits_{j\in J} \delta(\xi-c_j)\int_{t_{j-1}}^{t_j} \psi(t)\,dt.$$
	Note that $g$ is a generalized distribution if $J\neq \emptyset$ and that $g$ coincides with the Fourier transform of $\mathcal{L}^{(x)} \psi$ up to some constant.
	Since the inverse function $h_j^{-1}: [\xi^{(x)}_{j,\min},\xi^{(x)}_{j,\max}] \to [t_{j-1},t_j]$ is a bijection, we have supp $h_j^{-1}(\xi) = [\xi^{(x)}_{j,\min},\xi^{(x)}_{j,\max}]$. Additionally, we have $g(c_j)\neq 0$ as $j \in J$, and $c_j$ belongs to the interval $[\xi_{\min}^{(x)},\xi_{\max}^{(x)}]$. Hence, we show the support of the function $g$ as follows:
	$$\mbox{supp}(g(\xi)) \subset \left\{\bigcup\limits_{j\notin J}\mbox{supp} (h_j^{-1})\right\} \bigcup \left\{ c_j,\,j \in J\right\} = [\xi_{\min}^{(x)},\xi_{\max}^{(x)}].$$

Summing up the above arguments we arrive at
	\begin{lemma}\label{lem:supp}
		Let $\Gamma = \{y: y=a(t), t\in [t_{\min},t_{\max}]\} \subset \mathbb{R}^3$ be a $C^1$-smooth curve with $t_{\max}>t_{\min}$.
		Then
		\begin{equation}
			(\mathcal{F}\mathcal{L}^{(x)}\psi)(\xi) = \sqrt{2\pi} \left(\sum\limits_{j\notin J} \psi(h_j^{-1}(\xi))\, |(h_j^{-1}(\xi))'| + \sum\limits_{j\in J} \delta(\xi-c_j)\int_{t_{j-1}}^{t_j} \psi(t)\,dt\right).
\end{equation}	
	
Moreover,
		\begin{equation*}
			{\rm supp} (\mathcal{F}\mathcal{L}^{(x)}\psi) \subset [\xi^{(x)}_{\min},\xi^{(x)}_{\max}].
		\end{equation*}
		\end{lemma}
		Below we provide a sufficient condition to ensure trivial intersections of the ranges of two data-to-pattern operators corresponding to different trajectories. 		
\begin{lemma}\label{lem:interrange}
 		Let $\Gamma_a = \{y: y=a(t), t\in [t_{\min},t_{\max}]\} \subset \mathbb{R}^3$ and $\Gamma_b = \{y: y=b(t), t\in [t_{\min},t_{\max}]\} \subset \mathbb{R}^3$ be $C^1$-smooth curves such that
 		\begin{eqnarray}\nonumber
 			&&\left[\inf\limits_{t\in [t_{\min},t_{\max}]}(t+c^{-1}|x-a(t)|),\sup\limits_{t\in [t_{\min},t_{\max}]}(t+c^{-1}|x-a(t)|)\right] \\ \label{condition-T}
			&\bigcap& \left[\inf\limits_{t\in [t_{\min},t_{\max}]}(t+c^{-1}|x-b(t)|),\sup\limits_{t\in [t_{\min},t_{\max}]}(t+c^{-1}|x-b(t)|)\right] = \emptyset.			
		\end{eqnarray}
Let $ \mathcal{L}^{(x)}_a$ and $ \mathcal{L}^{(x)}_b$ be the data-to-pattern operators associated with  $\Gamma_a$ and $\Gamma_b$, respectively.
		Then $\mbox{Range}( \mathcal{L}^{(x)}_a) \cap \mbox{Range}( \mathcal{L}^{(x)}_b) = \{0\}$.
	\end{lemma}

	\begin{proof}
		Let $f_a,f_b \in L^2(t_{\min},t_{\max})$ be such that $ (\mathcal{L}^{(x)}_a f_a)(\tau)=   (\mathcal{L}^{(x)}_b f_b)(\tau) \coloneqq  Q(\tau, x)$. We need to prove $Q(\cdot, x) \equiv 0$. By the definition of $ \mathcal{L}$ (see \eqref{tildeL}), the function
		\begin{equation*}
			 \tau \to Q(\tau, x) =  \int_{t_{\min}}^{t_{\max}} e^{i\tau (t+c^{-1}|x-a(t)|)} f_a(t)\,dt
			=  \int_{t_{\min}}^{t_{\max}} e^{i\tau (t+c^{-1}|x-b(t)|)} f_b(t)\,dt
		\end{equation*}
		belongs to $L^2(0, K)$.
		Since $Q(\tau, x)$ is analytic in $\tau\in \R$, the previous relation is well defined for any $\tau \in \mathbb{R}$.
		By Definition \ref{DIDP}, we suppose that $\{t_j\}_{j=1}^{n-1}$ and $\{\tilde{t}_j\}_{j=1}^{m-1}$ are division points of the functions $h_{a}(t) = t+c^{-1}|x-a(t)|$ and $h_{b}(t) = t+c^{-1}|x-b(t)|$, respectively.
Analogously we define $h_{j,a}(t):= t+c^{-1}|x-a_j(t)|$, $h_{j,b}(t) := t+c^{-1}|x-b_j(t)|$, and $J_a:=\{j\in \N : 1\leq j\leq n, h_{j,a}^{\prime}(t)\equiv 0, t\in (t_{j-1},t_j)\}$, $J_b:=\{j\in \N : 1\leq j\leq m, h_{j,b}^{\prime}(t)\equiv 0, t\in (\tilde{t}_{j-1},\tilde{t}_j)\}$. Denote $h_{j,a}(t)\equiv c_{j,a}$ for $j\in J_a$ and $h_{j,b}(t)\equiv c_{j,b}$ for $j\in J_b$.

Using the formula \eqref{Lxi}, the function $ Q(\cdot, x)$ can be rewritten as the inverse Fourier transforms:
		\begin{equation}\label{G}
			 {Q}(\tau, x) = \int_{\mathbb{R}} e^{i\tau  \xi}  {g}_a( \xi, x)\,d \xi =\int_{\mathbb{R}} e^{i\tau  \xi}  {g}_b( \xi, x)\,d \xi,
		\end{equation}
		with
		\begin{equation*}
			g_a(\xi,x)=\sum\limits_{j\notin J_a} f_a(h_{j,a}^{-1}(\xi))\, |(h_{j,a}^{-1}(\xi))'| + \sum\limits_{j\in J_a} \delta(\xi-c_{j,a})\int_{t_{j-1}}^{t_j} f_a(t)\,dt,
		\end{equation*}
		\begin{equation*}
			g_b(\xi,x)=\sum\limits_{j\notin J_b} f_b(h_{j,b}^{-1}(\xi))\, |(h_{j,b}^{-1}(\xi))'| + \sum\limits_{j\in J_b} \delta(\xi-c_{j,b})\int_{\tilde{t}_{j-1}}^{\tilde{t}_j} f_b(t)\,dt.
		\end{equation*}
This implies $g_a(\xi,x)=g_b(\xi,x)$ for all $\xi\in R$.	On the other hand,	the support sets of $g_a$ and $g_b$ satisfy
		\begin{equation*}
			\begin{aligned}
				&{\rm supp}\,  {g}_a(\cdot, x) \subset \left[\inf\limits_{t\in [t_{\min},t_{\max}]}(t+c^{-1}|x-a(t)|),\sup\limits_{t\in [t_{\min},t_{\max}]}(t+c^{-1}|x-a(t)|)\right] ,\\
				&{\rm supp}\,  {g}_b(\cdot, x) \subset \left[\inf\limits_{t\in [t_{\min},t_{\max}]}(t+c^{-1}|x-b(t)|),\sup\limits_{t\in [t_{\min},t_{\max}]}(t+c^{-1}|x-b(t)|)\right].
			\end{aligned}
		\end{equation*}
Hence,	by the condition \eqref{condition-T} we obtain $ {g}_a( \xi, x) =  {g}_b( \xi, x) \equiv 0$ for all $ \xi \in \mathbb{R}$ . In view of \eqref{G}, we get $ Q(\cdot,x) \equiv 0$. 	
	\end{proof}

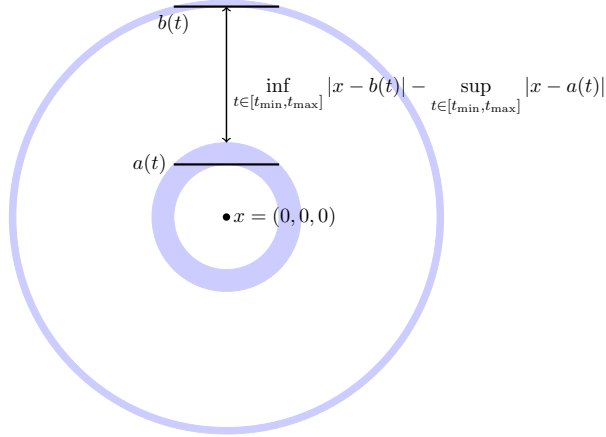
\begin{figure}[!ht]
		\centering
		\scalebox{0.7}{
		\begin{tikzpicture}
		\node (x) at (0,0) {};
		\filldraw[blue!20,even odd rule](x)circle(4.123)(x)circle(4);
		\filldraw[blue!20,even odd rule](x)circle(1.414)(x)circle(1);

		\draw [ very thick] (-1,4) -- (1,4);
		\draw [ very thick] (-1,1) -- (1,1);

		\draw (-1,1) node [left] {$a(t)$};
		\draw (-1,4) node [below] {$b(t)$};
		\draw (x) node [right]{$x = (0,0,0)$};
		\fill (x)	 circle (2pt);


		\draw [<->,thick] (0,1.414) -- (0,4);
		\draw (0,2.3) node [right] {$\inf\limits_{t\in [t_{\min},t_{\max}]}  |x-b(t)|-\sup\limits_{t\in [t_{\min},t_{\max}]} |x-a(t)|$};

		\end{tikzpicture}
		}
		\caption{Illustration of two trajectories $a(t) $ and $b(t)$ in the $Ox_1x_3$-plane such that $\mbox{Range}( \mathcal{L}^{(x)}_a) \cap \mbox{Range}( \mathcal{L}^{(x)}_b) = \{0\}$. Here we set $x=(0,0,0)$, $c=1$,  $a(t)=(t-2,0,1)$ and $b(t) = (t-2,0,4)$ with $t\in [1,3]$ (the black segments). Since $\inf\limits_{t\in [t_{\min},t_{\max}]}  |x-b(t)| =4$ and $\sup\limits_{t\in [t_{\min},t_{\max}]} |x-a(t)|=\sqrt{2}$, the condition \eqref{suff-con} is fulfilled, which implies the condition \eqref{condition-T}. }
		\label{two-tra}
	\end{figure}

\begin{remark}
	A sufficient condition to ensure \eqref{condition-T} is
	\begin{equation}\label{suff-con}
		\inf\limits_{t\in [t_{\min},t_{\max}]}  |x-b(t)| > \sup\limits_{t\in [t_{\min},t_{\max}]} |x-a(t)| +c(t_{\max}-t_{\min}).
	\end{equation}
	  In Fig. \ref{two-tra} we show an example of two orbit functions which satisfy the condition \eqref{suff-con}.
	
\end{remark}

	For any $y\in \R^3$, define the parameter-dependent test functions $\phi^{(x)}_{y}\in L^2(0, K)$  by
	\begin{equation}\label{testfunc}
		\phi^{(x)}_{y}(\omega)=\frac{1}{|t_{\max}-t_{\min}|}\int_{t_{\min}}^{t_{\max}} e^{i\omega  (t+c^{-1}|x-y|)} dt, \quad \omega\in(0, K).
	\end{equation}
Here we emphasize that the test function $\phi^{(x)}_{y}$ depends on both the observation point $x\in S_R$ and the sampling point $y\in \R^3$.
The Fourier transform of the aforementioned test function is given as follows.
	\begin{lemma}\label{INVF-PHI}
		We have
		\begin{equation}
			[\mathcal{F} \phi^{(x)}_{y}](\tau)=
			\left\{\begin{array}{lll}
				\sqrt{2\pi}/|t_{\max}-t_{\min}| \quad &&\mbox{if}\quad \tau\in \left[ t_{\min}+c^{-1}|x-y|,\; t_{\max}+c^{-1}|x-y| \right], \\
				0 \quad&&\mbox{if otherwise}.
			\end{array}\right.
		\end{equation}
	\end{lemma}
	
	\begin{proof}
		Letting $\tau = t+c^{-1}|x-y|$, we can rewrite the function $\phi^{(x)}_{y}$ as
		\begin{equation*}
			\phi^{(x)}_{y}(\omega)=\int_{\R} e^{i\omega \tau}g_y(\tau, x)\,d\tau,
			\end{equation*}
			where\begin{equation*}
			g_y(\tau,x):=\left\{
			\begin{aligned}
			&\frac{1}{|t_{\max}-t_{\min}|} &&{\rm if}\, \tau \in \left[t_{\min}+c^{-1}|x-y|,\; t_{\max}+c^{-1}|x-y|\right],\\
			&0 &&\mbox{ if  otherwise}.
			\end{aligned}\right.
		\end{equation*}
		Therefore, $[\mathcal{F}\phi^{(x)}_{y}](\tau) = \sqrt{2\pi}g_y(\tau,x)$.
	\end{proof}
In the following we present a necessary condition imposed on the observation point $x$ and radiating period $T:=t_{\max}-t_{\min}$ to guarantee that the test function $\phi^{(x)}_y$ lies in the range of the data-to-pattern operator.	
	\begin{lemma}\label{lem3.5}
		If $\phi^{(x)}_y \in {\rm Range} (\mathcal{L}^{(x)})$ for some $y\in \R^3$, we have $\xi^{(x)}_{\max} - \xi^{(x)}_{\min} \geq T$. Here $\xi^{(x)}_{\max}$ and $\xi^{(x)}_{\min}$ are defined by \eqref{def:xi}.
			\end{lemma}

	\begin{proof}
		If $\phi^{(x)}_y \in {\rm Range} (\mathcal{L}^{(x)})$, there exists a function $\psi\in L^2(t_{\min},t_{\max})$ such that $\phi^{(x)}_y = \mathcal{L}^{(x)} \psi$ in $L^2(0,K)$. Since both $\phi^{(x)}_y$ and $\mathcal{L}^{(x)} \psi$ are analytic functions over $\R$, it holds that $\phi^{(x)}_y(\omega) = (\mathcal{L}^{(x)} \psi) (\omega)$ for all $\omega\in \R$. Then their support sets must be identical, i.e., supp$(\mathcal{F}\phi^{(x)}_y)=$ supp$(\mathcal{F}\mathcal{L}^{(x)} \phi) \subset [\xi^{(x)}_{\min},\xi^{(x)}_{\max}]$,
where we have used Lemma \ref{lem:supp}. Hence, the length of supp$(\mathcal{F}\phi^{(x)}_y)$, which can be seen from Lemma \ref{INVF-PHI}, must be less than or equal to that of $[\xi^{(x)}_{\min},\xi^{(x)}_{\max}]$, i.e.,		
\ben \xi^{(x)}_{\max} - \xi^{(x)}_{\min} \geq t_{\max} - t_{\min}=T.\enn		
	\end{proof}	
From the above lemma we conclude that $\phi^{(x)}_y \notin {\rm Range} (\mathcal{L}^{(x)})$ for all $y\in \R^3$, if $\xi^{(x)}_{\max} - \xi^{(x)}_{\min}< T$. Inspired by this fact we introduce the concept of observable points.
 \begin{definition}\label{obd}
		Let $\xi^{(x)}_{\min}$ and $\xi^{(x)}_{\max}$ be the maximum and minimum of the function $h(t)=t+c^{-1}|x-a(t)|\in C^1[t_{\min}, t_{\max}]$ (see \eqref{def:xi}), respectively.
		The point $x\in \R^3$ is called an observable point if $\xi^{(x)}_{\max} - \xi^{(x)}_{\min} \geq T$. The point $x$ is called non-observable if $\xi^{(x)}_{\max} - \xi^{(x)}_{\min} < T$.
	\end{definition}

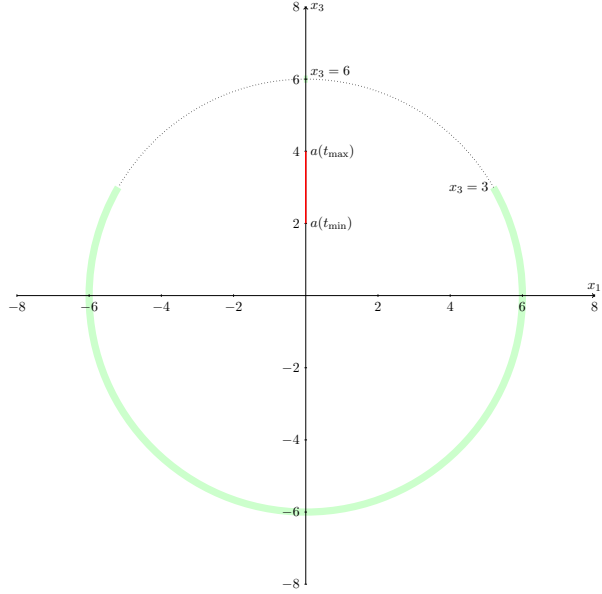
\begin{figure}

	\centering
		\scalebox{0.48}{
		\begin{tikzpicture}
		\fill[green!20] (0:5.9) --(0:6.1) arc (0:30:6.1) --(30:5.9) arc (30:0:5.9);
		\fill[green!20] (150:5.9) --(150:6.1) arc (150:180:6.1) --(180:5.9) arc (180:150:5.9);


		\fill[green!20] (180:5.9) --(180:6.1) arc (180:360:6.1) --(360:5.9) arc (360:180:5.9);
		\fill[green!20] (89.5:5.9) --(89.5:6.1) arc (89.5:90.5:6.1) --(90.5:5.9) arc (90.5:89.5:5.9);

		\draw[->] (-8,0) -- (8,0) node[above] {$x_1$} coordinate(x axis);
		\draw[->] (0,-8) -- (0,8) node[right] {$x_3$} coordinate(y axis);
		\foreach \x/\xtext in {-8,-6,-4,-2, 2, 4, 6, 8}
		\draw[xshift=\x cm] (0pt,1pt) -- (0pt,-1pt) node[below] {$\xtext$};
		\foreach \y/\ytext in {-8,-6,-4,-2, 2, 4, 6, 8}
		\draw[yshift=\y cm] (1pt,0pt) -- (-1pt,0pt) node[left] {$\ytext$};

		\draw [dotted] (5.2,3) arc [ start angle = 30, end angle = 150, radius = 6];
		\draw [very thick,red] (0, 2) -- (0,4);

		\draw (0, 2) node [right] {$a(t_{\min})$};

		\draw (0,4) node [right] {$a(t_{\max})$};
		\draw (30:6) node [left] {$x_3=3$};
		\draw (90:6.2) node [right] {$x_3=6$};

		\end{tikzpicture}
		}

		\caption{Illustration of observable (green arc) and non-observable (dotted arc) points for the trajectory $a(t) = (0,0,2t)$ for $t\in [1,2]$ in the $Ox_1x_3$-plane.	}\label{S2}
\end{figure}

\begin{figure}[ht]
		\centering
		\scalebox{1.2}{
		\begin{tikzpicture}
		\draw[line width=2.5cm,color=red!20] (1.25,-2.5) -- (1.25,2.5);
		\draw[line width=2.5cm,color=green!20] (-1.25,-2.5) -- (-1.25,2.5);
		\draw[->] (-2,0) -- (2,0) node[above] {$x_1$} coordinate(x axis);
		\draw[->] (0,-2) -- (0,2) node[right] {$x_2$} coordinate(y axis);
		\foreach \x/\xtext in {-2,-1, 1, 2}
		\draw[xshift=\x cm] (0pt,1pt) -- (0pt,-1pt) node[below] {$\xtext$};
		\foreach \y/\ytext in {-2,-1, 1, 2}
		\draw[yshift=\y cm] (1pt,0pt) -- (-1pt,0pt) node[left] {$\ytext$};

		\draw [very thick] (-0.5,0) arc [ start angle =180, end angle =360 ,radius =0.5];

		\draw (-0.7, 0) node [above] {$a(t_{\min})$};

		\draw (0.7, 0) node [above] {$a(t_{\max})$};

		\draw (0.5, -0.5) node [below] {$\Gamma$};

		\end{tikzpicture}
		}
		\caption{Illustration of the observable points (green area excluding the trajectory $\Gamma$) and non-observable points (red area excluding the trajectory $\Gamma$) for the trajectory $a(t) = (0.5\cos t,0.5\sin t,0)$ with $t\in [\pi,2\pi]$ in the $Ox_1x_2$-plane.}\label{A1}
		
	\end{figure}
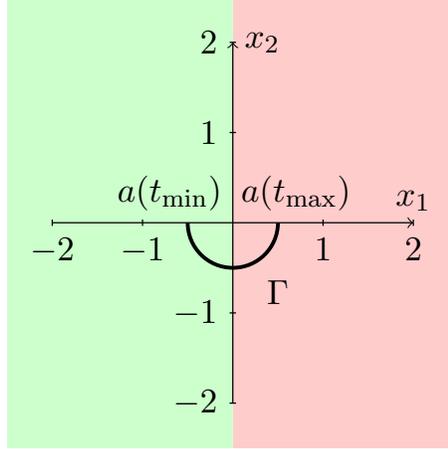

We remark that the set of observable points is uniquely determined by the orbit function $a(t)$ in conjunction with the starting and terminal time points $t_{\min}$ and $t_{\max}$. In the case of non-observable points $x$, our approach does not yield any information about the orbit function, a fact that will be elucidated in the second assertion of Theorem \ref{Th:factorization}.
For an observable point $x$ satisfying $(x-a(t))\cdot a^{\prime}(t) \leq 0$ for all $t \in [t_{\min},t_{\max}]$, we will show that it is possible to reconstruct the smallest annulus encompassing the trajectory and centered at $x$.
However, in cases that $(x-a(t))\cdot a^{\prime}(t) \leq 0$ for all $t \in [t_{\min},t_{\max}]$ is not fulfilled, one can only except to image a slimmer annulus centered at the observable point $x$.
In the subsequent sections, we proceed with the observable points/positions for orbit functions defined by a straight line (see Fig. \ref{S2}) and a semi-circle (see Fig. \ref{A1}) in three dimensions. In both examples, we assume $c=1$.

\vspace{.1in}	
\textbf{Example 1: A straight line segment in $\R^3$.}

Consider an acoustic point source which is moving along a straight line.

\begin{lemma}\label{line-ob}
	Define the orbit function $a(t) := (0,0,2t)\in \R^3$ for $t\in[1,2]$. Then the point $x=(x_1,x_2,x_3)\in \R^3$, $|x|=6$ is observable if $x_3\in [-6,\frac{3-\sqrt{33}}{2}]\bigcup [3,6]$.
\end{lemma}

		\begin{proof}
		From the expression of the orbit function $a(t)$, we have
		\begin{equation*}
			\begin{aligned}
				h(t) = t+|x-a(t)| &= t + \sqrt{x_1^2+x_2^2+(x_3-2t)^2}= t + \sqrt{(2t-x_3)^2+36-x_3^2},\\
				h^{\prime}(t) = 1+|x-a(t)|' &= 1+ \frac{2(2t-x_3)}{\sqrt{(2t-x_3)^2+36-x_3^2} }.
			\end{aligned}
		\end{equation*}
		We notice that $h'(t)\geq 0$ as $t \geq t_0$ and  $h'(t)< 0$ as $t < t_0$, where $t_0 := \frac{x_3}{2}-\sqrt{3-\frac{x_3^2}{12}}$. Hence, there are three cases for the relationship between $t_0$ and $[1,2]$.

	Case (i): If $t_0 \leq 1$, then $x_3 -2 \leq \sqrt{12-\frac{x_3^2}{3}}$, which means $x_3\in [-6,\frac{3+\sqrt{33}}{2}]$. In this case, $h(t)$ is monotonically increasing in $[1,2]$. So, if $x$ is observable, we have
	\begin{equation*}
		h(2)-h(1) \geq 1,
	\end{equation*}	
	that is,
	\begin{equation*}
		(4-x_3)^2 \geq (2-x_3)^2.
	\end{equation*}		
	Thus, $x\in S_6$ is an observable point if $x_3\in [-6,3]$.

Case (ii):	If $t_0 \in [1,2]$, then $x_3 -4 \leq \sqrt{12-\frac{x_3^2}{3}} \leq x_3 -2$, which means $x_3\in [\frac{3+\sqrt{33}}{2},3+\sqrt{6}]$. In this case, $h(t)$ is monotonically decreasing in $[1,t_0]$ and monotonically increasing in $[t_0,2]$. We notice that
	\begin{equation*}
		\max\{h(1),h(2)\} -h(t_0) < 1,
	\end{equation*}	
	for all $x_3\in [\frac{3+\sqrt{33}}{2},3+\sqrt{6}]$.

Case (iii):	If $t_0 \geq 2$, then $x_3 -4 \geq \sqrt{12-\frac{x_3^2}{3}}$, which means $x_3\in [3+\sqrt{6},6]$. In this case, $h(t)$ is monotonically decreasing in $[1,2]$. So, if $x$ is observable, we have
	\begin{equation*}
		h(1)-h(2) \geq 1,
	\end{equation*}	
	that is,
	\begin{equation*}
		\sqrt{(2-x_3)^2+36-x_3^2} - \sqrt{(4-x_3)^2+36-x_3^2} \geq 2.
	\end{equation*}		
	Thus, $x\in S_6$ is an observable point if $x_3=6$.

To sum up, we deduce that an observable point $x\in \R^3, |x|=6$ should fulfill the relation
		$$x_3 \in [-6,3]\bigcup\{6\}.$$
\end{proof}

\vspace{.1in}	
\textbf{Example 2: A semi-circle in $\R^3$.}

Suppose that an acoustic point source moves along a semi-circle centered at $z=(z_1,z_2,z_3)\in \R^3$.
\begin{lemma} \label{lem:circle}
Let the orbit function be $a(t) = (0.5\cos t+z_1,0.5\sin t+z_2,z_3)\in \R^3$ for $t\in[\pi,2\pi]$. Then $x=(x_1,x_2,x_3)\notin \Gamma$ is observable if $x_1 \leq z_1$.
\end{lemma}

\begin{proof}
	From the expression of the orbit function $a(t)$, we have
	\begin{equation*}
		h(t)=t+|x-a(t)| = t+\sqrt{(x_1-z_1-0.5\cos t)^2+(x_2-z_2-0.5\sin t)^2+(x_3-z_3)^2}.
	\end{equation*}
	It is obvious that $|a'(t)|<1$. Then we get $h'(t)>0$ for all $t\in [\pi,2\pi]$, that is, the function $h(t)$ is monotonically increasing in $[\pi,2\pi]$. Hence,
	\begin{eqnarray*}
		\xi_{\min}^{(x)} &&= \pi +\sqrt{(x_1-z_1+0.5)^2+(x_2-z_2)^2+(x_3-z_3)^2},\\
		\xi_{\max}^{(x)} &&= 2\pi +\sqrt{(x_1-z_1-0.5)^2+(x_2-z_2)^2+(x_3-z_3)^2}.
	\end{eqnarray*}
	If $x$ is observable, we have
	\begin{equation*}
	\begin{aligned}
		\xi_{\max}^{(x)}-\xi_{\min}^{(x)} &= \pi +\sqrt{(x_1-z_1-0.5)^2+(x_2-z_2)^2+(x_3-z_3)^2} \\
		& -\sqrt{(x_1-z_1+0.5)^2+(x_2-z_2)^2+(x_3-z_3)^2} \\
		& \geq T =\pi,
	\end{aligned}
	\end{equation*}
	that is,
	\begin{equation*}
		(x_1-z_1-0.5)^2 \geq (x_1-z_1+0.5)^2.
	\end{equation*}
	One can find $x_1 \leq z_1$ through simple calculations.
\end{proof}

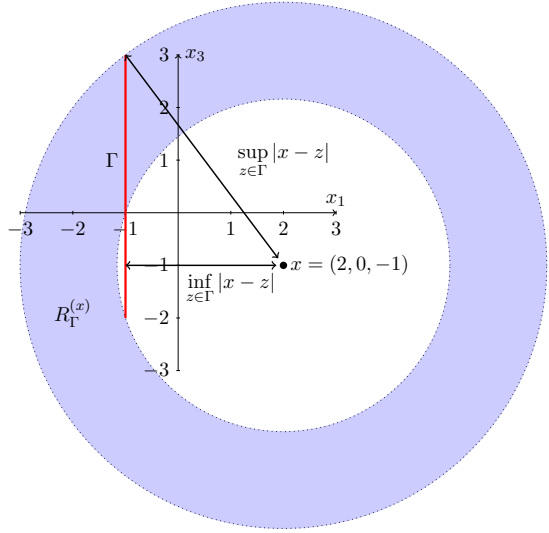
\begin{figure}[!ht]
		\centering
		\scalebox{0.7}{
		\begin{tikzpicture}
		\node (x) at (2,-1) {};
		\filldraw[blue!20,even odd rule](x)circle(5)(x)circle(3.1623);
		\draw[->] (-3,0) -- (3,0) node[above] {$x_1$} coordinate(x axis);
		\draw[->] (0,-3) -- (0,3) node[right] {$x_3$} coordinate(y axis);
		\foreach \x/\xtext in {-3,-2,-1, 1, 2, 3}
		\draw[xshift=\x cm] (0pt,1pt) -- (0pt,-1pt) node[below] {$\xtext$};
		\foreach \y/\ytext in {-3,-2,-1, 1, 2, 3}
		\draw[yshift=\y cm] (1pt,0pt) -- (-1pt,0pt) node[left] {$\ytext$};

		\draw [ very thick,red] (-1,-2) -- (-1,3);

		\draw (-1,1) node [left] {$\Gamma$};
		\draw (-2,-1.5) node [below] {$R_\Gamma^{(x)}$};
		\draw (x) node [right]{$x = (2,0,-1)$};
		\fill (x)	 circle (2pt);

		\draw [dotted] let \p1 = ($ (x) - (-1,-2) $),  \n2 = {veclen(\x1,\y1)}
		in (x) circle (\n2);
		\draw [dotted] let \p1 = ($ (x) - (-1,3) $),  \n2 = {veclen(\x1,\y1)}
		in (x) circle (\n2);

		\draw [<->,thick] (x) -- (-1,3);
		\draw [<->,thick] (x) -- (-1,-1);
		\draw (1,1) node [right] {$\sup\limits_{z\in \Gamma}|x-z|$};
		\draw (1,-1) node [below] {$\inf\limits_{z\in \Gamma}|x-z|$};

		\end{tikzpicture}
		}
		\caption{Illustration of the annulus $A_{\Gamma}^{(x)}$ (blue area) with $x = (2,0,1)$ in the $Ox_1x_3$-plane. The wave propagates with a speed of unity. Here the curve $a(t)=(-1, 0,4-t),\, t\in [1,6]$ denotes the orbit (the red segment) of a point source moving from below to above. There holds $|x-a(t_{\max})| = \sqrt{10}, \,|x-a(t_{\min})| =5,\,\inf\limits_{z\in \Gamma}|x-z|=3,\,\sup\limits_{z\in \Gamma}|x-z|=5$. In this case the annulus $A_{\Gamma}^{(x)}$ is a subset of $\{y\in \R^3: \inf\limits_{z\in \Gamma} |x-z| \leq |x-y| \leq \sup\limits_{z\in \Gamma} |x-z|\}$. }
		\label{kgamma}
	\end{figure}

Given the trajectory $\Gamma = \{y: y=a(t),\, t\in [t_{\min}, t_{\max}]\}$, the set $$\Lambda_\Gamma:= \{y\in \R^3: \inf\limits_{z\in \Gamma} |x-z| \leq |x-y| \leq \sup\limits_{z\in \Gamma} |x-z|\}$$ represents the smallest annulus encompassing $\Gamma$ and centered at the point $x$.
Intuitively, it is reasonable to anticipate recovering this annulus utilizing multi-frequency data collected at a single observation point.
	If $x$ is an observable point, we define the annulus (see Fig. \ref{kgamma})
	\begin{equation}\label{K}
		A_{\Gamma}^{(x)}:=\{y\in \mathbb{R}^3: c(\xi^{(x)}_{\min} - t_{\min} )   \leq |x-y| \leq c(\xi^{(x)}_{\max} - t_{\max})\} \subset \mathbb{R}^3.
	\end{equation}

	\begin{remark}\label{remark3.9}
		If the point source remains stationary at $z\in \mathbb{R}^3$, that is, $\Gamma=\{z\}$, then every point $x\in \mathbb{R}^3\backslash\{z\}$ is observable. In such a scenario, the annulus $A_{\Gamma}^{(x)}$ is reduced to the sphere $\left\{ y\in \mathbb{R}^3: |x-y|=|x-z|\right\}$, since $\xi^{(x)}_{\min} = t_{\min}+c^{-1}|x-z|$ and $\xi^{(x)}_{\max} = t_{\max}+c^{-1}|x-z|$. This signifies that the set of non-observable points is caused by the motion of the point source.
	\end{remark}
	
If $h^{\prime}(t)>0$ for $t \in (t_{\min},t_{\max})$, we have
	\begin{equation*}
		A_{\Gamma}^{(x)}=\{y\in \mathbb{R}^3: |x-a(t_{\min})| \leq |x-y| \leq |x-a(t_{\max})| \}
	\end{equation*}
	which is a subset of $\Lambda_\Gamma$. Moreover, $A_{\Gamma}^{(x)}$ coincides with $\Lambda_\Gamma$ when $(x-a(t))\cdot a^{\prime}(t) \leq 0$ for all $t \in [t_{\min},t_{\max}]$, because $$|x-a(t)|'=-\frac{x-a(t)}{|x-a(t)|}\cdot a^{\prime}(t) \geq 0\quad\mbox {for all}\quad t \in [t_{\min},t_{\max}],$$
implying that
\ben
&\xi^{(x)}_{\min} = t_{\min}+c^{-1}|x-a(t_{\min})|, \quad
\xi^{(x)}_{\max} = t_{\max}+c^{-1}|x-t_{\max}|,\\
&|x-a(t_{\min})| = \inf\limits_{z\in \Gamma} |x-z|, \quad |x-a(t_{\max})| = \sup\limits_{z\in \Gamma} |x-z|.
\enn

	If $h^{\prime}(t)<0$ for $t \in (t_{\min},t_{\max})$, there holds
	\begin{equation*}
		A_{\Gamma}^{(x)}=\{y\in \mathbb{R}^3: |x-a(t_{\max})|+cT \leq |x-y| \leq |x-a(t_{\min})|-cT \},
	\end{equation*}
which is also a subset of
$\Lambda_\Gamma$; see Lemma \ref{lem3.9} below.

	\begin{lemma}\label{lem3.9}
Let $x\in S_R$ be an observable point. We have  $$\inf\limits_{z\in \Gamma} |x-z| \leq |x-y| \leq \sup\limits_{z\in \Gamma} |x-z| \quad\mbox{for all}\quad y\in A_{\Gamma}^{(x)}.$$
\end{lemma}

	\begin{proof}
		Suppose that
		$$\xi^{(x)}_{\min} = t_1+c^{-1}|x-a(t_1)| , \quad \xi^{(x)}_{\max} = t_2+c^{-1}|x-a(t_2)|,\quad\mbox{for some} \quad t_1,t_2\in [t_{\min}, t_{\max}].$$
		Therefore,		\ben
			&&c(\xi^{(x)}_{\min} - t_{\min}) =  ct_1+|x-a(t_1)| -ct_{\min}   \geq |x-a(t_1)| \geq \inf\limits_{z\in \Gamma} |x-z|,\\
			&&c(\xi^{(x)}_{\max} - t_{\max}) =  ct_2+|x-a(t_2)| - ct_{\max}  \leq |x-a(t_2)| \leq \sup\limits_{z\in \Gamma} |x-z|.
		\enn
		This implies that for $ y\in A_{\Gamma}^{(x)}$,
		\ben
		|x-y|\leq c(\xi^{(x)}_{\max} - t_{\max}) \leq \sup\limits_{z\in \Gamma} |x-z|,\quad
	|x-y|\geq c(\xi^{(x)}_{\min} - t_{\min}) \geq \inf\limits_{z\in \Gamma} |x-z|.	
				\enn
	\end{proof}

If $x\in S_R$ is observable, we shall prove that the test function $\phi^{(x)}_y$ defined by \eqref{testfunc} lies in the range of $\mathcal{L}^{(x)}$ if and only if $y \in A_{\Gamma}^{(x)}$. This together with \eqref{RI}
establishes a computational criterion for imaging $A_{\Gamma}^{(x)}$ from the multi-frequency near-field data $u(x,\omega)$ with $\omega \in [\omega_{\min}, \omega_{\max}]$.
We also need to discuss non-observable points.

	\begin{lemma}\label{lem3.4}
		(i) If $x$ is non-observable, we have $\phi^{(x)}_y \notin {\rm Range} (\mathcal{L}^{(x)})$ for all $y \in \mathbb{R}^3$.\\
		(ii) If $x$ is an observable point, we have $\phi^{(x)}_y \in {\rm Range} (\mathcal{L}^{(x)})$ if and only if $y \in A_{\Gamma}^{(x)}$.
	\end{lemma}

	\begin{proof}
	(i) The first assertion follows directly from Lemma \ref{lem3.5} and the Definition \ref{obd} for non-observable points.
		
(ii) If $x$ is an observable point, we have $\xi^{(x)}_{\max} - \xi^{(x)}_{\min} \geq T$. If $\phi^{(x)}_y \in {\rm Range} (\mathcal{L}^{(x)})$, one can find a function $\phi$ satisfying $\phi^{(x)}_y = \mathcal{L}^{(x)} \phi$. Then their support sets  must fulfill the relation supp$(\mathcal{F}\phi^{(x)}_y)=$ supp$(\mathcal{F}\mathcal{L}^{(x)} \phi) \subset [\xi^{(x)}_{\min}, \xi^{(x)}_{\max}]$ by Lemma \ref{lem:interrange}. Using Lemma \ref{INVF-PHI} yields
		$$\left[t_{\min}+c^{-1}|x-y|, t_{\max}+c^{-1}|x-y|\right] \subset [\xi^{(x)}_{\min}, \xi^{(x)}_{\max}].$$
		Hence, $t_{\min}+c^{-1}|x-y| \geq \xi^{(x)}_{\min}$ and $t_{\max}+c^{-1}|x-y| \leq \xi^{(x)}_{\max}$, leading to		
		$$c(\xi^{(x)}_{\min} - t_{\min}) \leq |x-y| \leq c(\xi^{(x)}_{\max} - t_{\max}).$$
This proves $y \in A_{\Gamma}^{(x)}$.

			To prove the reverse direction, we set  $\eta_{\min} = t_{\min}+c^{-1}|x-y|$ and $\eta_{\max} = t_{\max}+c^{-1}|x-y|$. If $y \in A_{\Gamma}^{(x)}$, we have
		\begin{equation*}
			[\eta_{\min},\eta_{\max}] = \left[t_{\min}+c^{-1}|x-y|, t_{\max}+c^{-1}|x-y|\right] \subset [\xi^{(x)}_{\min},\xi^{(x)}_{\max}].
		\end{equation*}
Setting $$\psi_y(t) := \frac{|h'(t)|}{T}\chi_y(t) \in L^2(t_{\min},t_{\max}),$$
where
$$\chi_y(t) = \left\{
\begin{aligned}	
	&1, \, &t\in Y,\\
	&0, \, &t\notin Y,
\end{aligned}
\right.$$
and
$$Y := \left\{t\in [t_{\min},t_{\max}]: h(t)\in [\eta_{\min},\eta_{\max}], h(t) \neq h(\tilde t) \mbox{ for all } \tilde t \in [t_{\min},t)\right\}.$$
One can find that $h(t):Y \to [\eta_{\min},\eta_{\max}]$ is a bijection. Then,
\begin{equation*}
	\begin{aligned}
		(\mathcal{L}^{(x)}\psi_y)(\omega) &= \int_{t_{\min}}^{t_{\max}} e^{i\omega (t+c^{-1}|x-a(t)|)}\psi_y(t)\,dt \\
		&= \int_{Y} e^{i\omega h(t)} \frac{|h'(t)|}{T} \,dt\\
		&= \int_{\eta_{\min}}^{\eta_{\max}} e^{i\omega\xi} \frac{1}{T} \,d\xi = \phi^{(x)}_y(\omega).
	\end{aligned}
\end{equation*}
Therefore, $\phi^{(x)}_y(\omega) \in \mbox{ Range}(\mathcal{L}^{(x)})$.
		
	\end{proof}	
	\begin{remark}
		The proof of Lemma \ref{lem3.4} (ii) corrects a mistake made in Lemma 3.11 of \cite{GHM2023}.
	\end{remark}

	\section{Indicator functions and uniqueness}\label{IdF}
If $x$ is an observable point, by Lemma \ref{lem3.4} and \eqref{RI}, the test functions $\phi^{(x)}_y$ can be effectively employed to define the characteristics function of $A_{\Gamma}^{(x)}$. Introduce the indicator function
	\begin{equation}\label{indicator}
		W^{(x)}(y) \coloneqq \left[\sum_{n=1}^\infty\frac{|\langle \phi^{(x)}_{y}, \psi_n^{(x)} \rangle|_{L^2(0, K)}^2}{ |\lambda_n^{(x)}|}\right]^{-1}, \qquad y\in \mathbb{R}^3.
	\end{equation}
 Combining Theorem \ref{range}, Lemma \ref{lem3.4} and Picard theorem, we obtain
	
	\begin{theorem}\label{Th:factorization}(Single observable point)
	
		If $x$ is an observable point, it holds that
		\ben
		W^{(x)}(y)=\left\{\begin{array}{lll}
		0 &&\quad\mbox{if}\quad y\notin A_{\Gamma}^{(x)},\\
		\mbox{finite positive number} &&\quad\mbox{if}\quad y\in A_{\Gamma}^{(x)}.		
		\end{array}\right.
		\enn
		If $x$ is non-observable, we have $W^{(x)}(y)=0$ for all $y\in \mathbb{R}^3$.

	\end{theorem}
	
	Therefore, for observable points, the values of $W^{(x)}$ within the annulus $A_{\Gamma}^{(x)}$ are expected to be larger compared to those in other regions. However, if $x$ is non-observable, the values of $W^{(x)}$ will identically vanish in $\mathbb{R}^3$.

	\begin{remark}
		The trajectory $\Gamma$ can not be uniquely determined by one observable point $x$. For example, let $\Gamma_1=\{z_1\}$ and $\Gamma_2=\{z_2\}$ be given by two stationary points such that
		\begin{equation*}
			z_1 \neq z_2, \, |x-z_1|=|x-z_2|.
		\end{equation*}
		Then, by Remark \ref{remark3.9}, we have $A_{\Gamma_1}^{(x)}=A_{\Gamma_2}^{(x)}=\{y\in \R^3: |x-y|=|x-z_j|,j=1,2\}$.
	\end{remark}

In the case of sparse observable points $\{x^{(j)}\in S_R: j=1,2,\cdots, M\}$, we shall make use of the following indicator function:
	\begin{equation}\label{W}
		W(y)= \left[\sum_{j=1}^M \frac{1}{W^{(x^{(j)})}(y)}\right]^{-1}=		
		\left[\sum_{j=1}^M\sum_{n=1}^\infty\frac{|\langle \phi^{(x^{(j)})}_{y}, \psi_n^{(x^{(j)})} \rangle|_{L^2(0, K)}^2}{ |\lambda_n^{(x^{(j)})}|}\right]^{-1}, \qquad y\in \mathbb{R}^3.
	\end{equation}
	Define the domain $D_\Gamma$ associated with the observable points $\{x^{(j)}: j=1,2,\cdots, M\}$ as
	\be\label{Theta}
		D_\Gamma \coloneqq \bigcap\limits_{j=1,2,\cdots, M} A_{\Gamma}^{(x^{(j)})}.
	\en
We can reconstruct $D_\Gamma$ from 	
the multi-frequency near-field data measured at sparse observable points.
\begin{theorem}\label{TH4.2}(Finite observable points)
		It holds that $0<W(y) <+\infty$ if $y \in D_\Gamma$ and
		$W(y)=0$ if $y \notin D_\Gamma$.
\end{theorem}

	\begin{proof}
		If $y\in D_\Gamma$, it means that $y\in A_{\Gamma}^{(x^{(j)})}$ for $j=1,2,\cdots,M$. By Theorem \ref{Th:factorization},
		\begin{equation}
			\sum\limits_{n=1}^{\infty} \frac{|\langle \phi^{(x^{(j)})}_{y}, \psi_n^{(x^{(j)})} \rangle|_{L^2(0, K)}^2}{ |\lambda_n^{(x^{(j)})}|} < +\infty \quad\mbox{for all}\quad j = 1,2,\cdots,M.
		\end{equation}
		Then the finite sum over the index $j$ must fulfill the relation $0<W(y)<+\infty$.

		If $y\notin D_\Gamma$, we may suppose without loss of generality that $y\notin A_{\Gamma}^{(x^{(1)})}$.
 By Theorem \ref{Th:factorization},
 		\begin{equation*}
		[W^{(x^{(1)})}(y)]^{-1}=	\sum\limits_{n=1}^{\infty} \frac{|\langle \phi^{(x^{(1)})}_{y}, \psi_n^{(x^{(1)})} \rangle|_{L^2(0, K)}^2}{ |\lambda_n^{(x^{(1)})}|} = \infty.
		\end{equation*}
Together with the definition of $W$, this gives
\ben
W(y)<\left[\sum\limits_{n=1}^{\infty} \frac{|\langle \phi^{(x^{(1)})}_{y}, \psi_n^{(x^{(1)})} \rangle|_{L^2(0, K)}^2}{ |\lambda_n^{(x^{(1)})}|}   \right]^{-1}=0.
\enn

	\end{proof}
Consequently, we arrive at the following uniqueness results, which seem unknown in the literature.
\begin{theorem}\label{TH4.3}(Uniqueness) Denote by $\Gamma=\{a(t): t\in[t_{\min}, t_{\max}]\}$ the trajectory of a moving point source where $a\in C^1[t_{\min}, t_{\max}]$.

(i) The domain $D_\Gamma$ associated with all observable points $x\in S_R$ (see \eqref{Theta}) can be uniquely determined by the multi-frequency data $\{u(x, \omega): x\in S_R, \omega\in(\omega_{\min}, \omega_{\max})\}$.

(ii) Let $x\in S_R$ be an arbitrarily fixed observable point. Then the annulus $A_{\Gamma}^{(x)}$ (see \eqref{K}) can be uniquely determined by the multi-frequency data $\{u(x, \omega):  \omega\in(\omega_{\min}, \omega_{\max})\}$. In particular, the annulus $\Lambda_\Gamma$ can be uniquely recovered if $(x-a(t))\cdot a^{\prime}(t) \leq 0$ for all $t\in [t_{\min}, t_{\max}]$.
\end{theorem}	

\begin{remark} Physically,
the condition $(x-a(t))\cdot a^{\prime}(t) \leq 0$ in the second assertion of Theorem \ref{TH4.3} means that the function $h(t)=t+c^{-1}|x-a(t)|$ is monotonically increasing and the function $|x-a(t)|$ is monotonically non-decreasing in $[t_{\min}, t_{\max}]$. 
\end{remark}

The second assertion of Theorem \ref{TH4.3} provides insight into the nature of information that can be extracted from multi-frequency data obtained at a single observable point. However, in the absence of any prior information on the orbit function, it remains unknown to classify observable and non-observable points.

\section{Numerical implementation}\label{num}

In this section, we will perform a series of numerical experiments to validate our algorithm in 3D. In practical scenarios, time-domain data is often inverse Fourier transformed to yield multi-frequency data. However, to streamline the numerical procedures for simulation purposes, we will exclusively conduct computational tests within the frequency domain, only.
Our primary objective is to extract information on the trajectory of a moving point source. This aim is accomplished through the utilization of multi-frequency near-field data recorded at either a single observation point or sparsely distributed observation points.

 Assuming a wave-number-dependent source term $f(x, k)$, as defined in (\ref{sourcef}), we can synthesize the near-field pattern using equation (\ref{expression-w}) by
\begin{equation}
\begin{aligned}
u(x, \omega)=
\int_{t_{\min}}^{t_{\max}} \frac{e^{i\omega(t+c^{-1}|x-a(t)|)}}{8\pi^2|x-a(t)|} \ell(t) \, dt,\quad x\in S_R \,,\; \omega\in(\omega_{\min}, \omega_{\max}),
\end{aligned}
\end{equation}
where $S_R=\{x\in \R^3:\, |x|=R\}$. The strength function $\ell(t)$ of the signal is defined as $\ell(t)=(t+1)^2$, which obvious fulfills the positivity constraint (\ref{cd:l}).  Below we will describe the process of inversion algorithm.
The frequency interval $(0,K)$ can be discretized by defining $$\omega_n=(n-0.5)\Delta \omega, \quad \Delta \omega:=\frac{K}{N}, \quad n=1,2,\cdots,N.$$
    To approximate the integral in (\ref{FarO}), we adopt {$2N-1$ samples $u(x, \kappa+\omega_n), n=1,2,\cdots,N$ and $u(x, \kappa-\omega_n), n=1,2,\cdots,N-1$}, of the near field and apply the midpoint rule. Therefore, we have
    \begin{equation}
    (\mathcal{N}^{(x)}\phi)(\tau_n) \approx \sum_{m=1}^{N} u(x, \kappa+\tau_n-s_m)\phi(s_m)\Delta \omega,
    \end{equation}
    where $\tau_n:=n\Delta \omega$ and $s_m:=(m-0.5)\Delta \omega$, $n,m=1,2,\cdots,N$.
    The Toeplitz matrix provides a discrete approximation of the near field operator $\mathcal{N}^{(x)}$:
\be \label{matF}
\mathcal N^{(x)}:= \Delta k \begin{pmatrix}
u({x},\kappa+\omega_1) & u({x},\kappa-\omega_1) & \cdots & u({x},\kappa-\omega_{N-2})  & u({x},\kappa-\omega_{N-1})  \\
u({x},\kappa+\omega_2) & u({x},\kappa+\omega_1) & \cdots & u({x},\kappa-\omega_{N-3}) &u({x},\kappa-\omega_{N-2})   \\
\vdots & \vdots  &  &\vdots &\vdots \\
u({x},\kappa+\omega_{N-1}) & u({x},\kappa+\omega_{N-2}) &  \cdots & u({x},\kappa+\omega_1) & u({x},\kappa-\omega_1)\\
u({x},\kappa+\omega_N) & u({x},\kappa+\omega_{N-1}) &  \cdots & u({x},\kappa+\omega_2) & u({x},\kappa+\omega_1)\\
\end{pmatrix}
\en
where $\mathcal{N}^{(x)}$ is a $N\times N$ complex matrix.
For any point $y\in \R^3$ we define the test function vector $\phi_y^{(x)} \in \C^N$ from (\ref{testfunc}) by
\be \label{testn}
\phi_y^{(x)}:= \left(\frac{1}{t_{max}-t_{min}}\int_{t_{min}}^{t_{max}} e^{i\tau_1(t+c^{-1}|x-y|)}dt, \,\cdots, \, \frac{1}{t_{max}-t_{min}}\int_{t_{min}}^{t_{max}} e^{i\tau_N(t+c^{-1}|x-y|)}dt\right).
\en
\noindent Denoting by  $\left\{ ( {\tilde \lambda^{( x)}_n}, \psi^{(x)}_n): n=1,2,\cdots,N \right\}$ an eigen-system of the matrix $\mathcal N^{(x)}$ (\ref{matF}), then one deduces that  an eigen-system of the matrix $(\mathcal N^{(x)})_\#:= |\real \mathcal N^{(x)})|+|\ima(\mathcal N^{(x)})|$ is $\left\{ ( \lambda^{(x)}_n, \psi^{(x)}_n): n=1,2,\cdots,N \right\}$, where $ \lambda^{(x)}_n:=|\real (\tilde \lambda^{(x)}_n)| +|\ima (\tilde \lambda^{( x)}_n)|$. We truncate the indicator function $W^{(x)}$ (\ref{indicator}) by
\be 
W^{({x})}(y):=\left[\sum_{n=1}^N\frac{\left| \phi^{({x})}_{y} \cdot \overline{\psi_n^{({x})} }\right|^2}{ |\lambda_n^{({x})}|}\right]^{-1}, \quad y\in \R^d,
\en
where $\cdot$ denotes the inner product in $\R^3$ and $N$ is consistent with the dimension of the  Toeplitz matrix \eqref{matF}.

The visualization of the annulus $A_\Gamma^{(x)}$ is attainable through the plot of $W^{(x)}(y)$. This visualization carries crucial information about the source trajectory, particularly when  using the data from an observable. Such visualizations can be used to describe the trajectories of moving point sources, no matter  they are characterized by straight lines or arcs. In the following figures, the original trajectory will be highlighted by the red solid line.
Unless otherwise specified, we assume $k_{\min}=0$ for the sake of simplicity.
The bandwidth can be extended from $(0,\omega_{max})$ to $(-\omega_{\max}, \omega_{\max})$ by $u(x, -\omega)=\overline{u(x, \omega)}$. Then, one deduces from these new measurement data with $\omega_{min}=-\omega_{max}$ that $\kappa=0$ and $K=\omega_{max}$.  The frequency band is represented by the interval $(0, 6)$ with $\omega_{\max}=6$, $N=10$ and $\Delta \omega=3/5$.

\subsection{One observation point}

In this subsection, the search domain is selected as a cube of the form $[-2,2]\times[-2,2]\times[-2,2]$ and the observation points are chosen from the set $\{x\in\R^3: |x|=2\}$. The observation points are then set on a sphere with a radius of $2$, such that $x=(2 \sin \varphi \cos \theta ,  2\sin \varphi \sin \theta, 2 \cos \varphi )$  for $\theta\in[0, 2\pi]$   and $\varphi\in[0,\pi]$. Figures below illustrate the slices  at $y_1=0$ or $y_2=0$.

\textbf{Example 1: A straight line segment in $\R^3$} \vspace{.1in}

\begin{figure}[H]
	\centering
    \subfigure[$\theta=0$, $\varphi=6\pi/9$]{
		\includegraphics[scale=0.22]{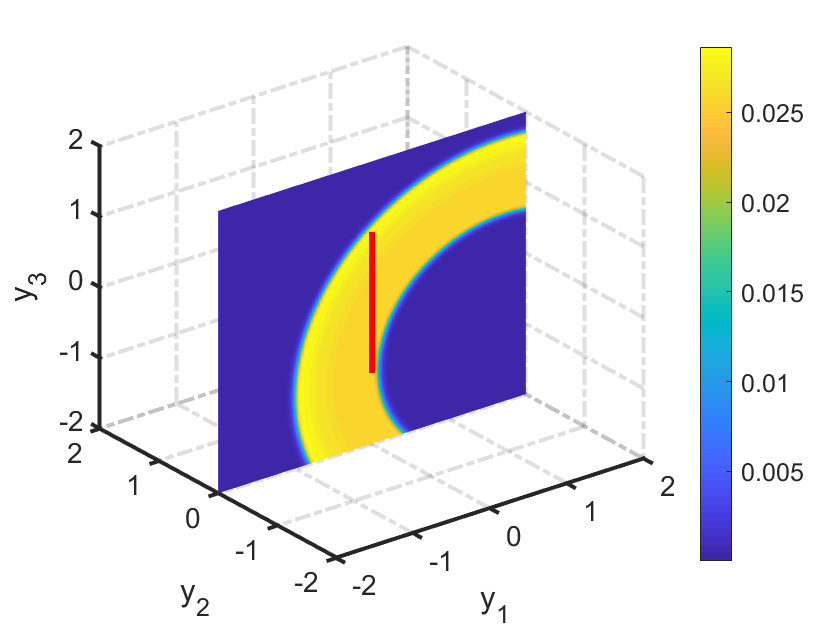}
		
	}
    \subfigure[$\theta=\pi/4$, $\varphi=7\pi/9$]{
		\includegraphics[scale=0.22]{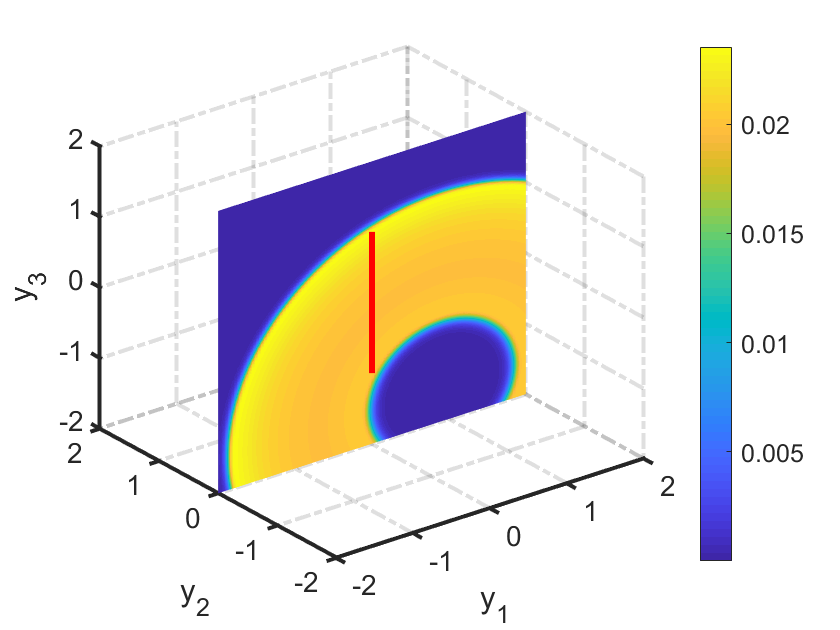}
		
	}
    \subfigure[$\theta=2\pi/4$, $\varphi=8\pi/9$]{
		\includegraphics[scale=0.22]{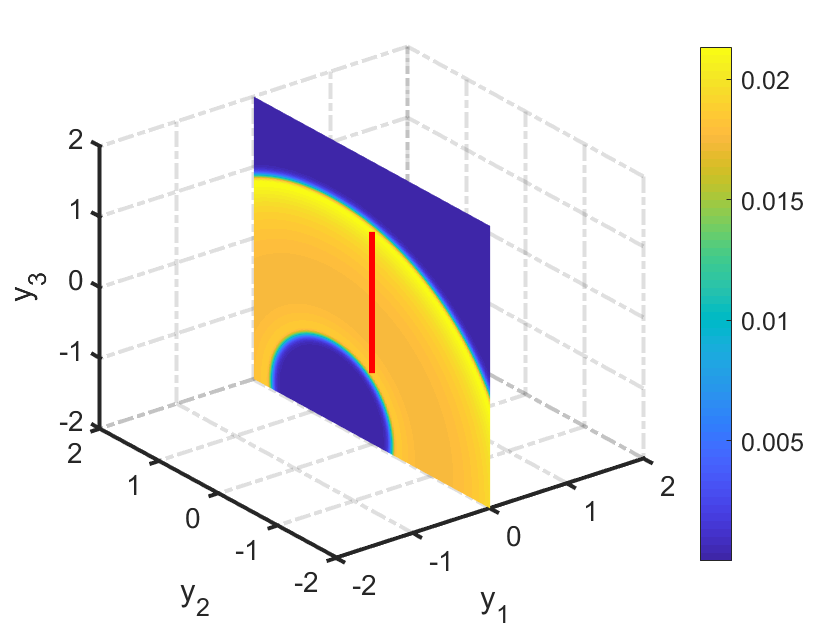}
		
	}
    \subfigure[$\theta=3\pi/4$, $\varphi=6\pi/9$]{
		\includegraphics[scale=0.22]{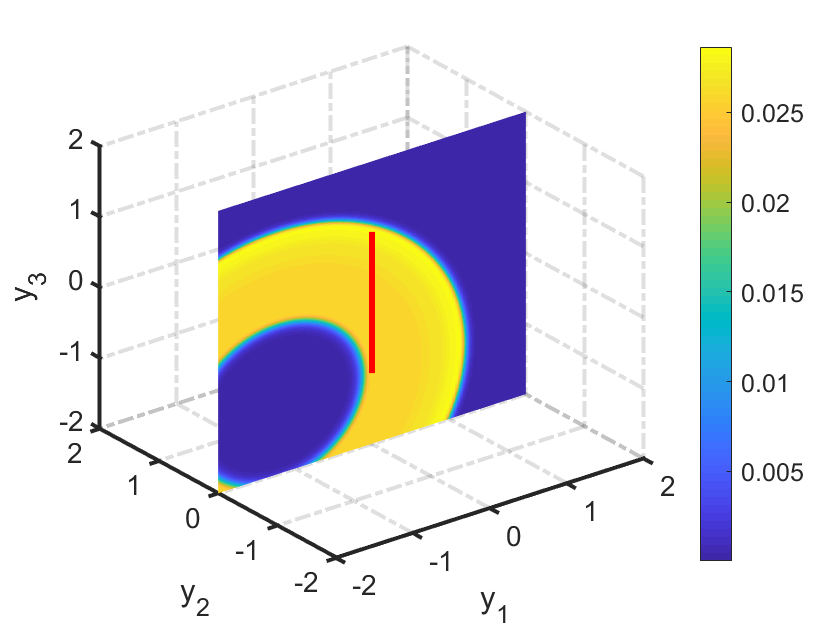}
		
	}
	\subfigure[$\theta=4\pi/4$, $\varphi=7\pi/9$ ]{
		\includegraphics[scale=0.22]{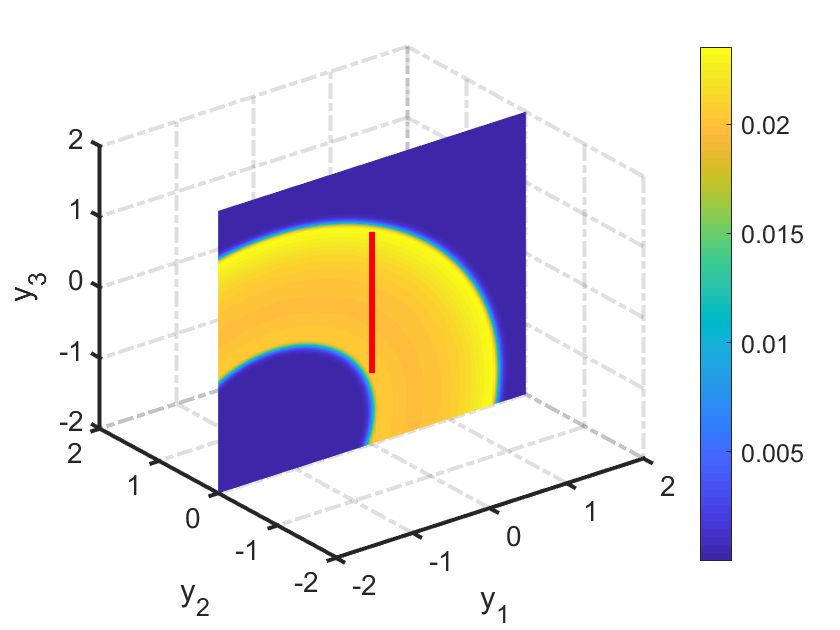}
		
	}
    \subfigure[$\theta=5\pi/4$, $\varphi=8\pi/9$  ]{
		\includegraphics[scale=0.22]{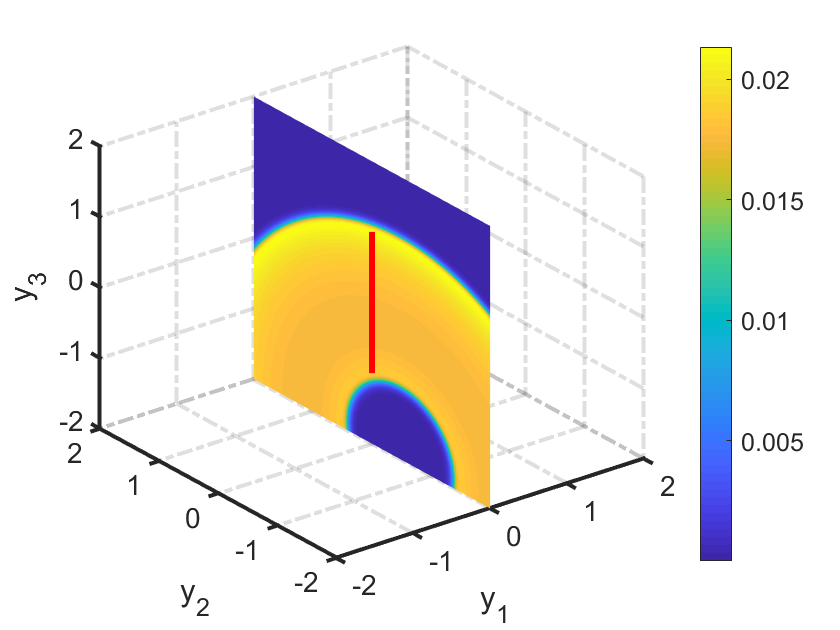}
	
	}
    \subfigure[$\theta=6\pi/4$, $\varphi=6\pi/9$ ]{
		\includegraphics[scale=0.22]{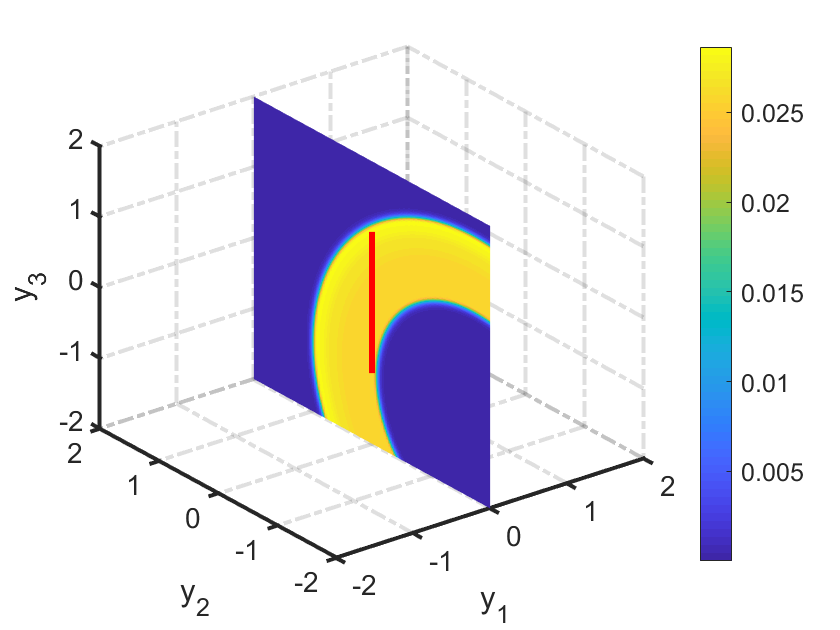}
	}
    \subfigure[$\theta=7\pi/4$, $\varphi=7\pi/9$ ]{
		\includegraphics[scale=0.22]{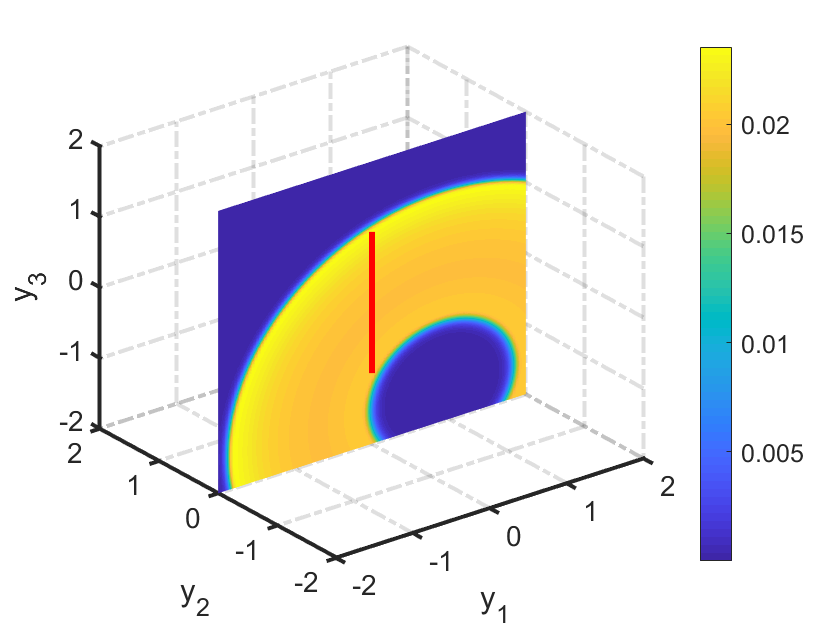}
	}
    \subfigure[$\theta=8\pi/4$, $\varphi=9\pi/9$ ]{
		\includegraphics[scale=0.22]{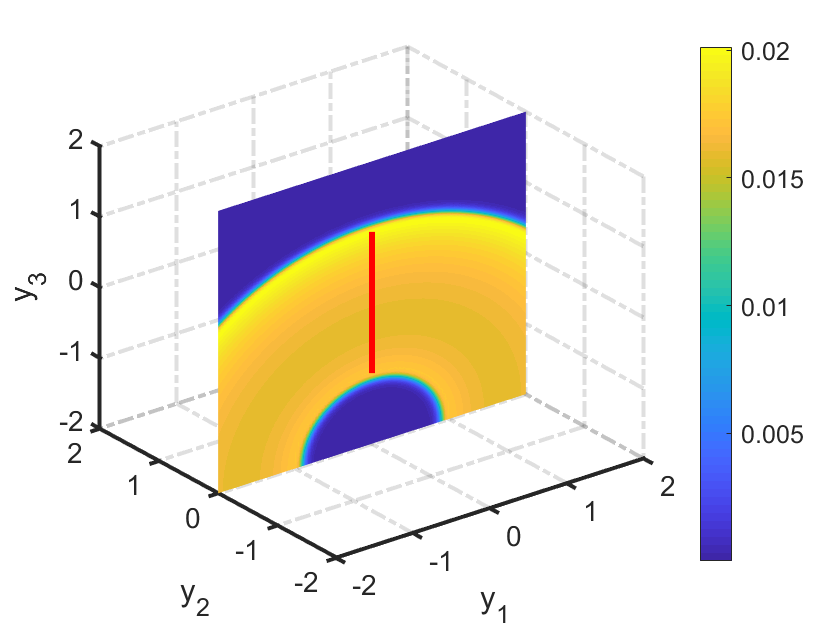}
	}
\caption{Reconstruction from a single observable point $x=(2 \sin \varphi \cos \theta, 2\sin \varphi \sin \theta, 2 \cos \varphi )$  with $\theta\in[0, 2\pi]$   and $\varphi\in[2\pi/3, \pi]$ for a straight line segment $a(t)=(0,0,t-2)$ where $t\in[1,3]$. Here it holds that $A_{\Gamma}^{(x)}= \Lambda_{\Gamma}^{(x)}$.} \label{fig:line1}
\end{figure}

\begin{figure}[H]
	\centering
    \subfigure[$\theta=0$, $\varphi=\pi/2$]{
		\includegraphics[scale=0.22]{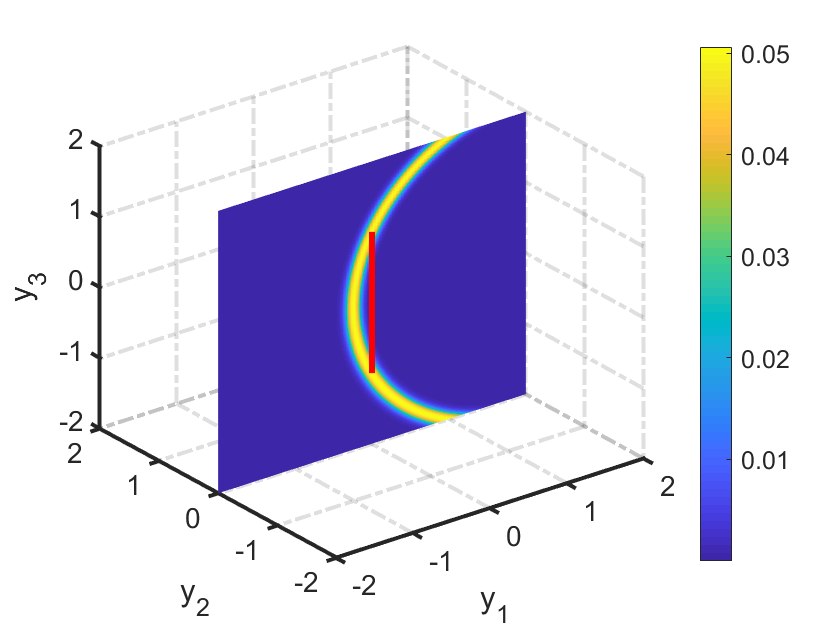}
		
	}
    \subfigure[$\theta=\pi/5$, $\varphi=9\pi/17$]{
		\includegraphics[scale=0.22]{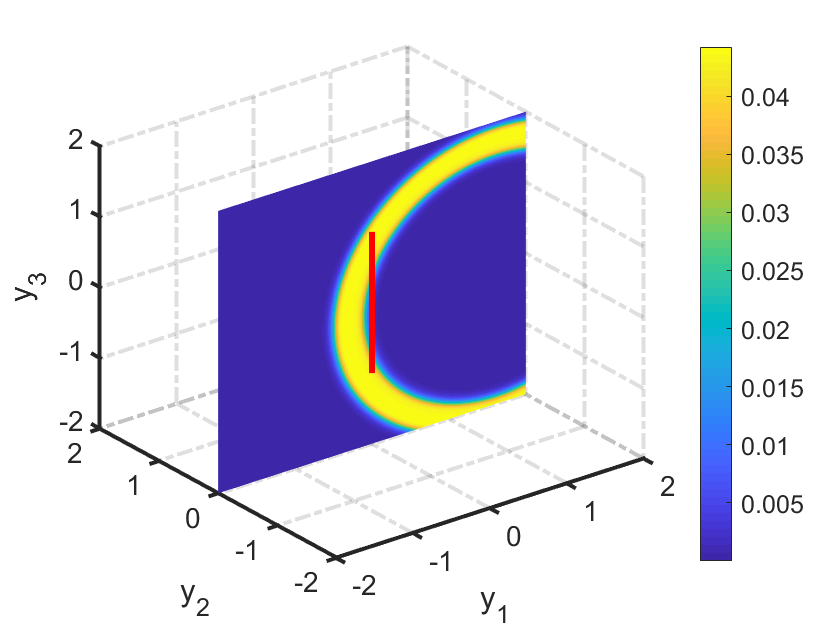}
		
	}
    \subfigure[$\theta=2\pi/5$, $\varphi=8\pi/15$]{
		\includegraphics[scale=0.22]{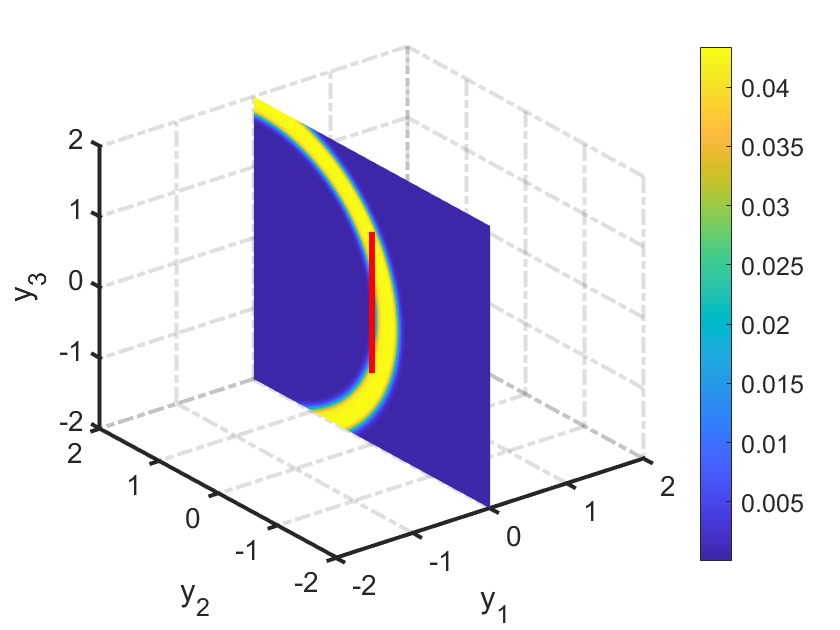}
		
	}
    \subfigure[$\theta=3\pi/5$, $\varphi=\pi/2$]{
		\includegraphics[scale=0.22]{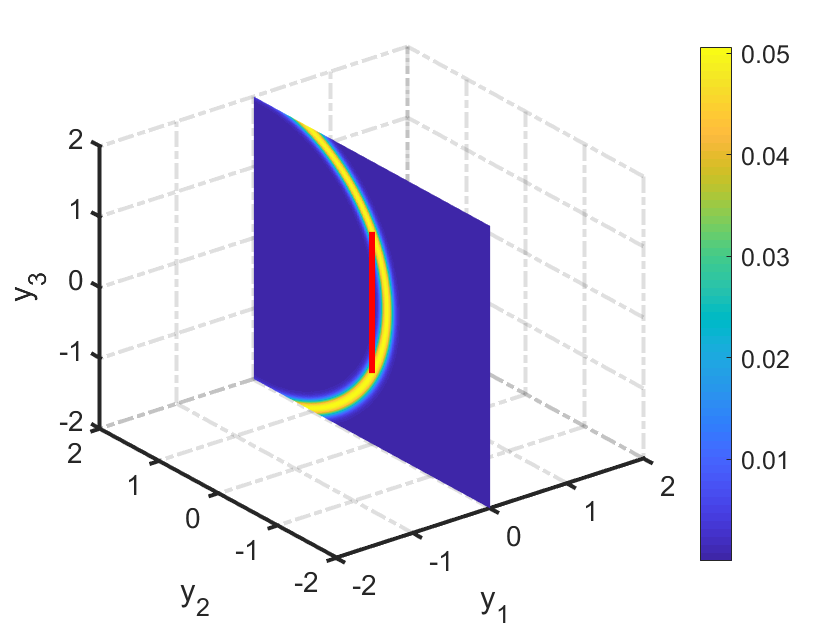}
		
	}
	\subfigure[$\theta=4\pi/5$, $\varphi=9\pi/17$ ]{
		\includegraphics[scale=0.22]{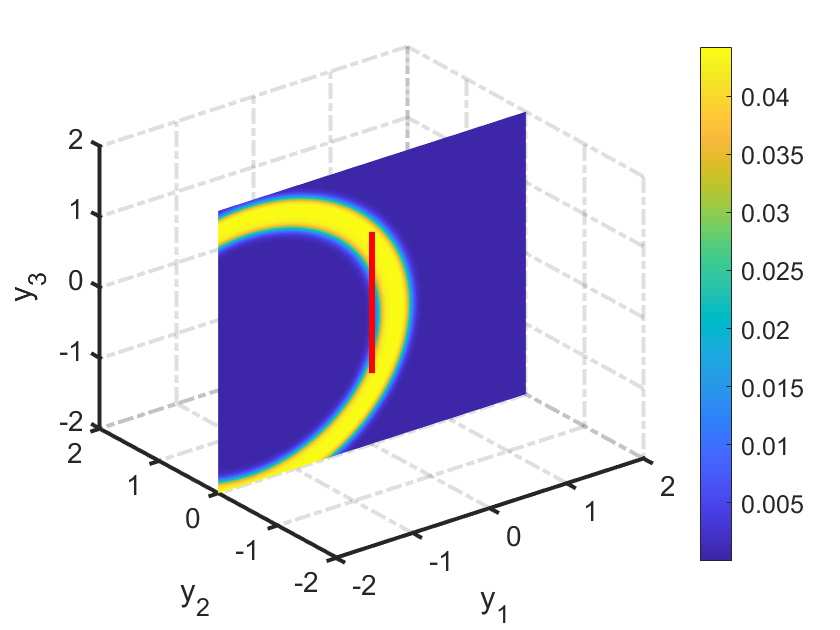}
		
	}
    \subfigure[$\theta=5\pi/5$, $\varphi=8\pi/15$]{
		\includegraphics[scale=0.22]{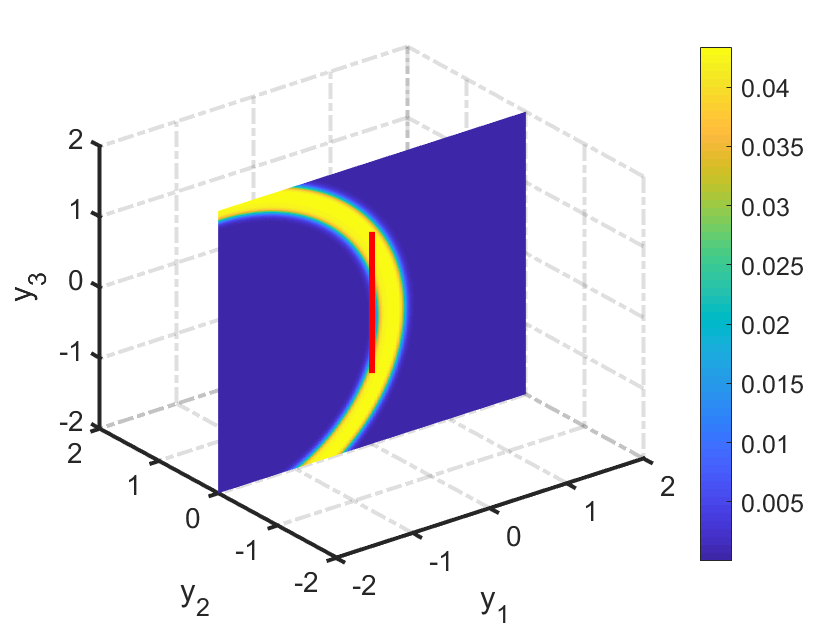}
	
	}
    \subfigure[$\theta=6\pi/5$, $\varphi=\pi/2$ ]{
		\includegraphics[scale=0.22]{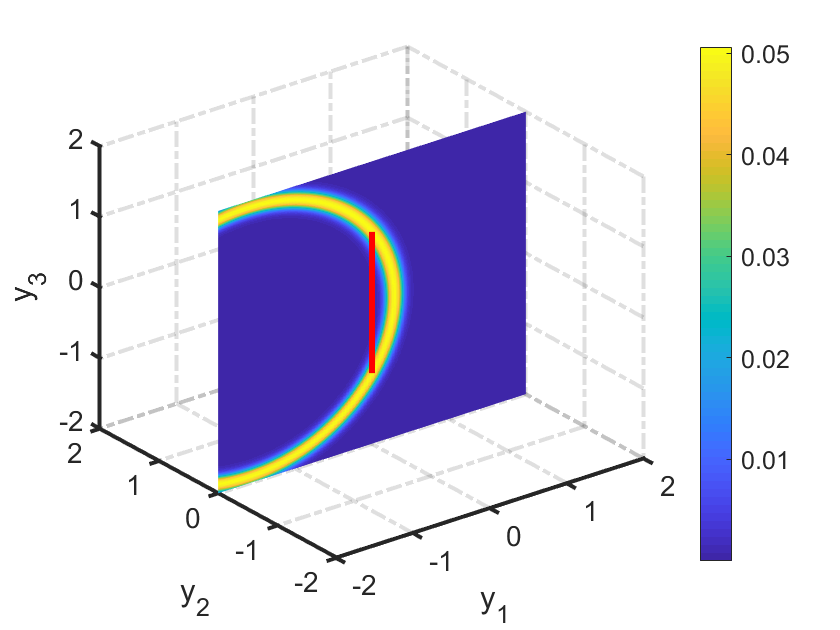}
	}
    \subfigure[$\theta=7\pi/5$, $\varphi=9\pi/17$ ]{
		\includegraphics[scale=0.22]{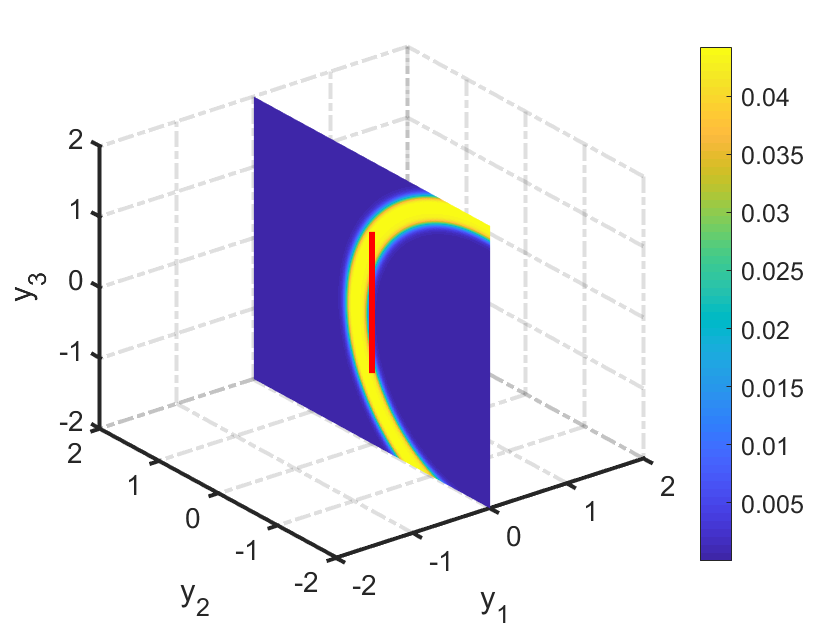}
	}
    \subfigure[$\theta=8\pi/5$, $\varphi=8\pi/15$ ]{
		\includegraphics[scale=0.22]{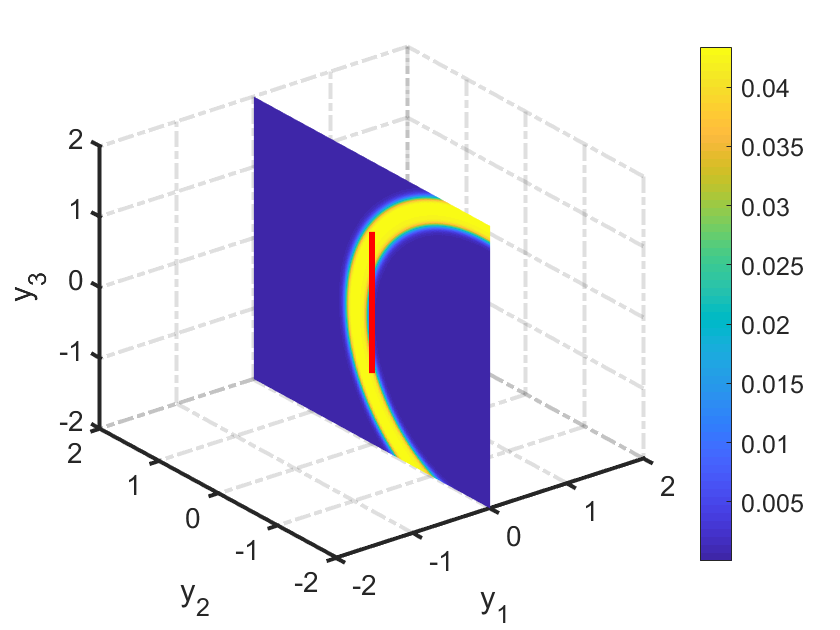}
	}
\caption{Reconstruction from a single observable point $x=(2 \sin \varphi \cos \theta,  2\sin \varphi \sin \theta, 2 \cos \varphi )$  with $\theta\in[0, 2\pi]$   and $\varphi\in[\pi/2,2\pi/3)$ for a straight line segment $a(t)=(0,0,t-2)$ where $t\in[1,3]$. Here it holds that $A_{\Gamma}^{(x)}\subset \Lambda_{\Gamma}^{(x)}$.}\label{fig:line2}
\end{figure}

\begin{figure}[H]
	\centering
    \subfigure[$\theta=0$, $\varphi=0$]{
		\includegraphics[scale=0.22]{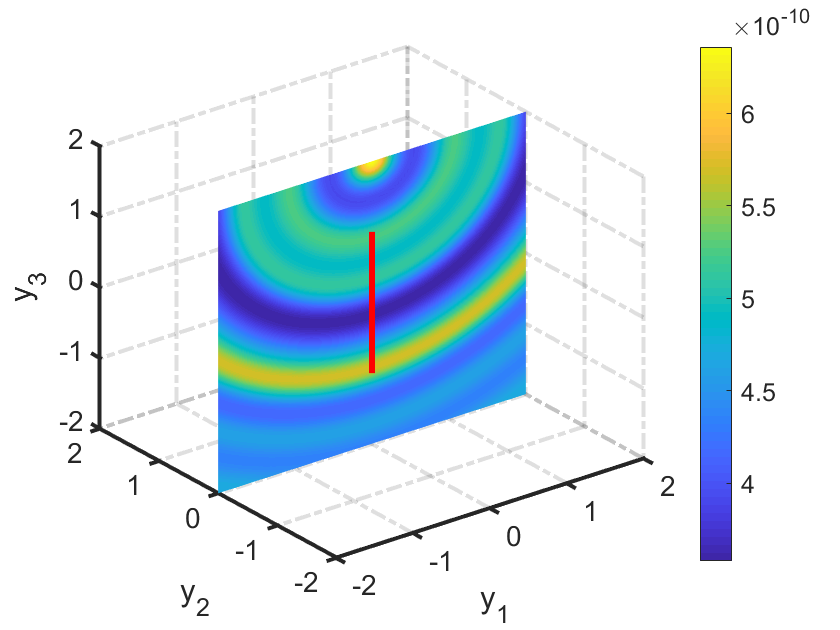}
		
	}
    \subfigure[$\theta=\pi/4$, $\varphi=2\pi/8$]{
		\includegraphics[scale=0.22]{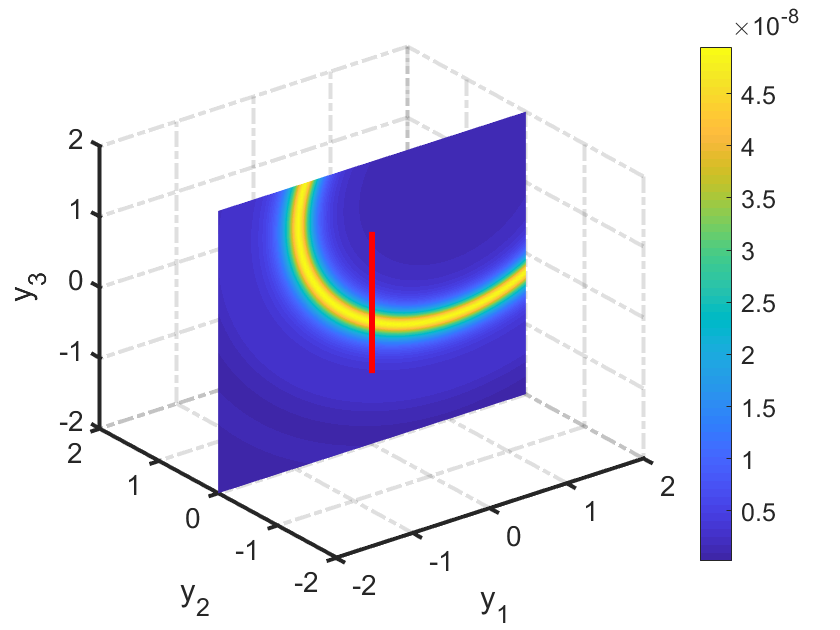}
		
	}
    \subfigure[$\theta=2\pi/4$, $\varphi=3\pi/8$]{
		\includegraphics[scale=0.22]{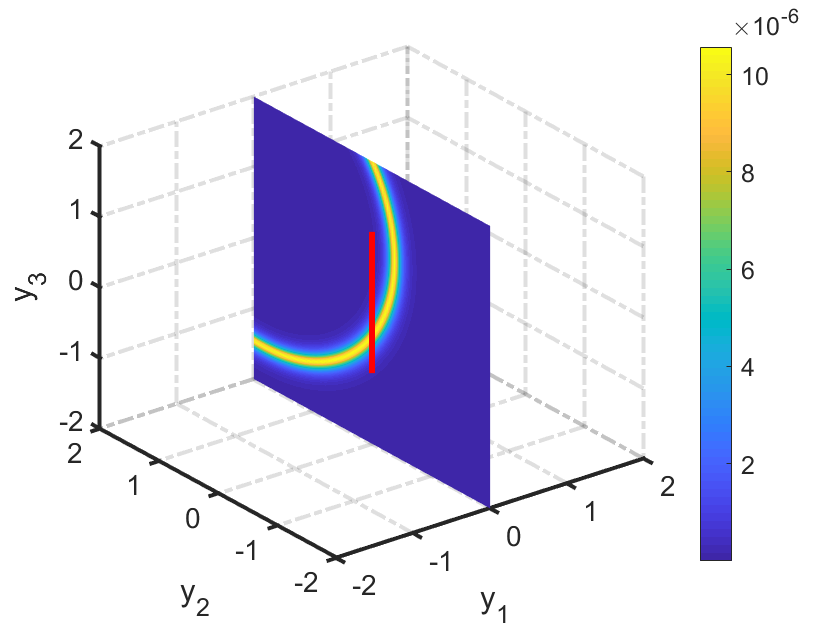}
		
	}
    \subfigure[$\theta=3\pi/4$, $\varphi=\pi/8$]{
		\includegraphics[scale=0.22]{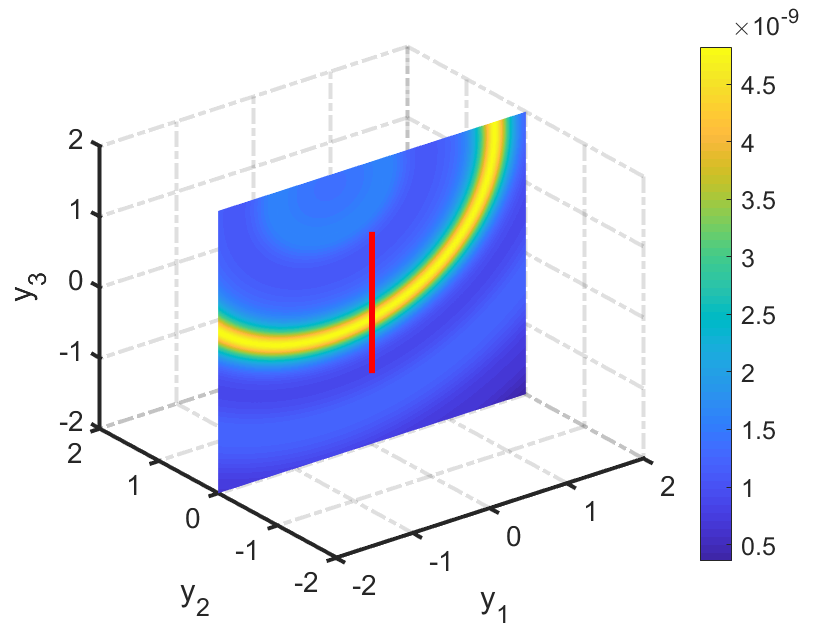}
		
	}
	\subfigure[$\theta=4\pi/4$, $\varphi=2\pi/8$ ]{
		\includegraphics[scale=0.22]{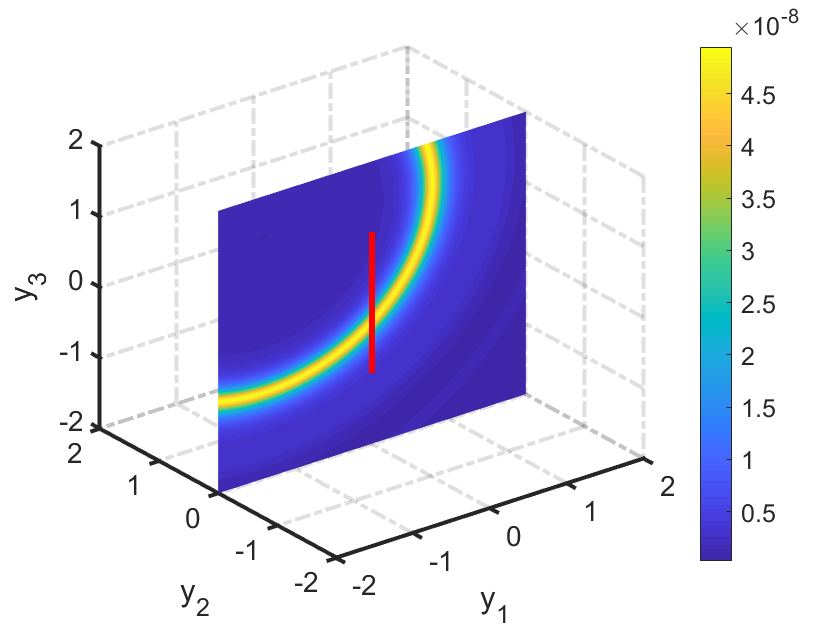}
		
	}
    \subfigure[$\theta=5\pi/4$, $\varphi=3\pi/8$  ]{
		\includegraphics[scale=0.22]{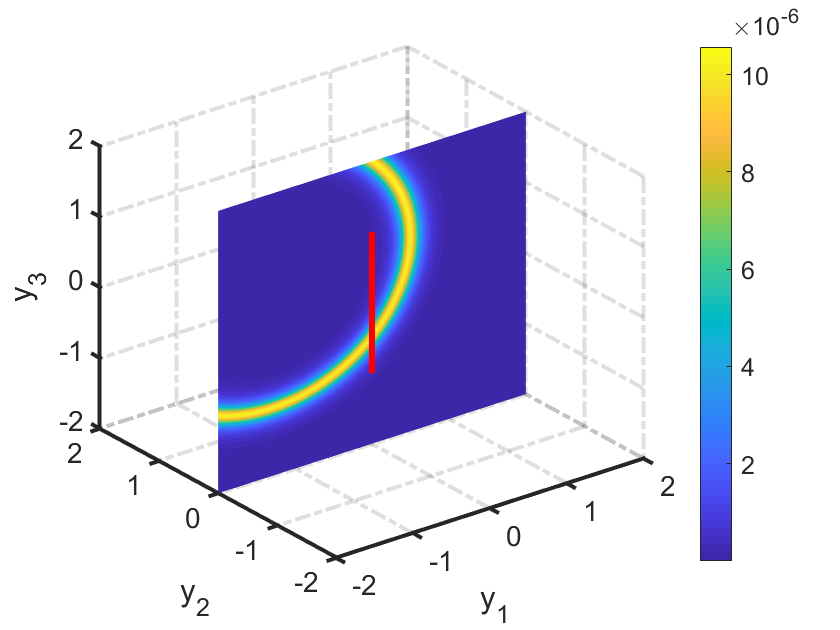}
	
	}
    \subfigure[$\theta=6\pi/4$, $\varphi=\pi/8$ ]{
		\includegraphics[scale=0.22]{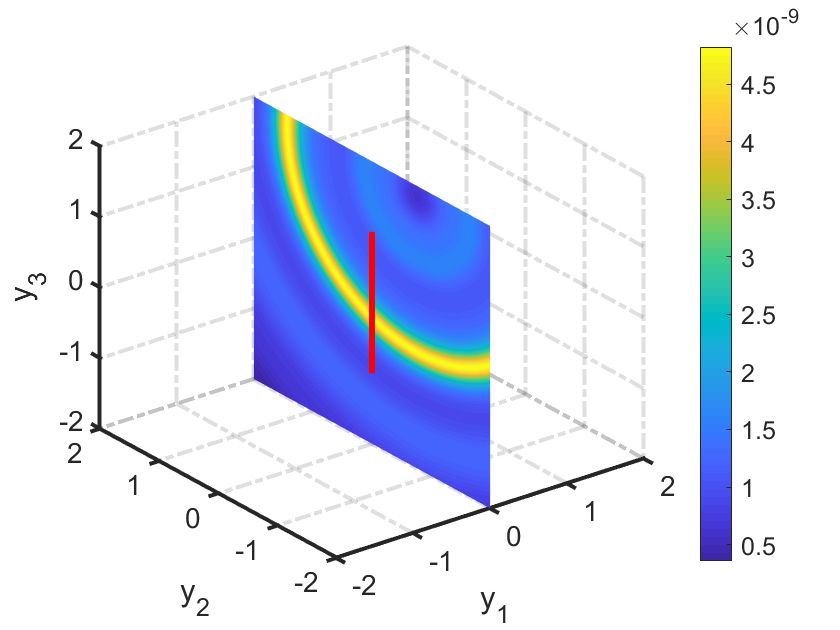}
	}
    \subfigure[$\theta=7\pi/4$, $\varphi=2\pi/8$ ]{
		\includegraphics[scale=0.22]{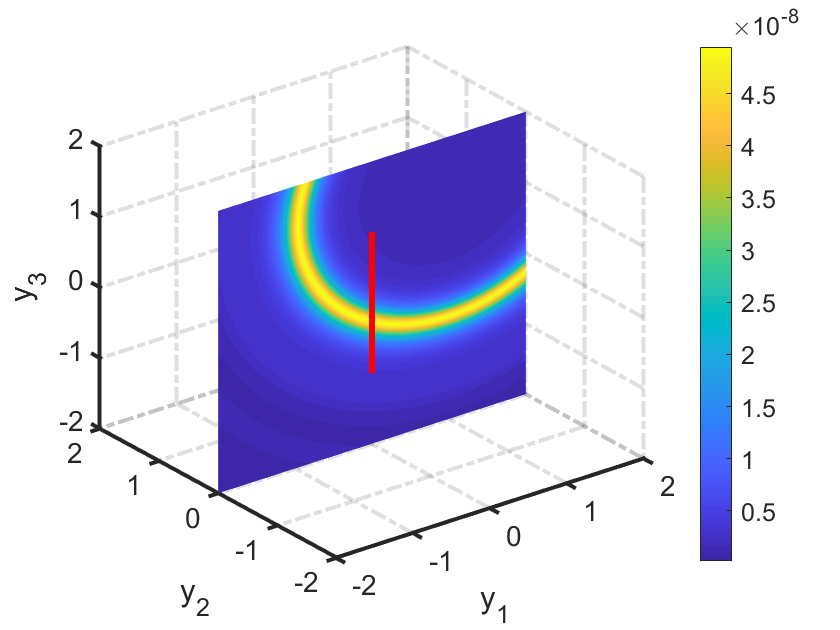}
	}
    \subfigure[$\theta=8\pi/4$, $\varphi=3\pi/8$ ]{
		\includegraphics[scale=0.22]{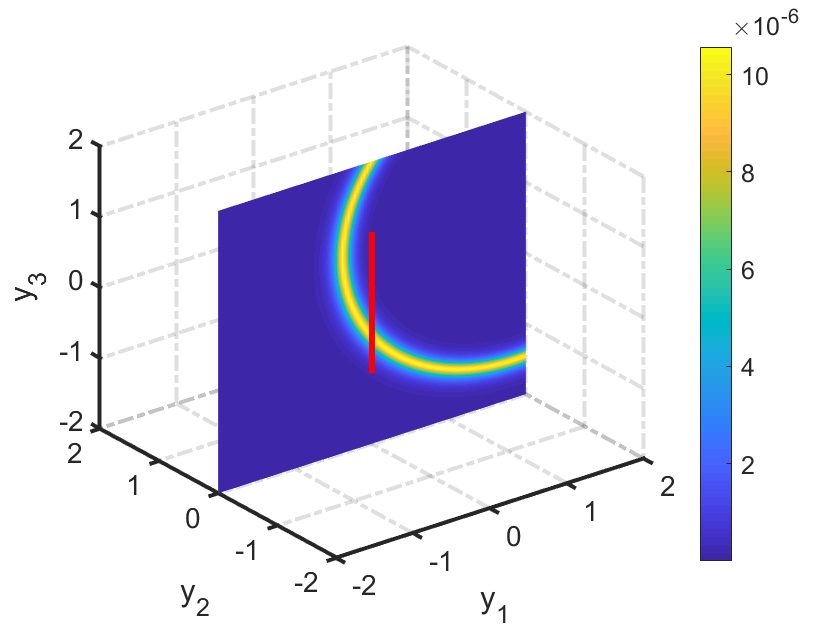}
	}
\caption{Reconstruction from a single non-observable point $x=(2 \sin \varphi \cos \theta,  2\sin \varphi \sin \theta,$ $2 \cos \varphi )$  with $\theta\in[0, 2\pi]$   and $\varphi\in[0,\pi/2)$ for a straight line segment $a(t)=(0,0,t-2)$ with $t\in[1,3]$. }\label{fig:line3}
\end{figure}

We examine a straight line segment from  Example 1, outlined in Section \ref{RangeLx}. Suppose that the trajectory of the moving point source is given by  $a(t)=(0, 0, t-2)$, where $t\in[1,3]$ and $x\in S_2=\{x\in\R^3:|x|=2\}$ represent the observation points. As the first step, we must classify the observable and non-observable points.
According to the orbit function, we have
\ben
&&h(t)=t+|x-a(t)|=t+\sqrt{x_1^2+x_2^2+(x_3-(t-2))^2},\\
&&h^{\prime}(t)=1-\frac{x-a(t)}{|x-a(t)|}\cdot a^{\prime}(t)=1+\frac{t-2 - x_3}{\sqrt{x_1^2+x_2^2+(x_3-(t-2))^2}}.
\enn
As the second term on the right-hand side of the above equation always falls in the range $[-1,1]$, it is evident that $h^{\prime}(t)>0$ for all $t\in[1,3]$, indicating that the function $h(t)$ is monotonically increasing over $[1, 3]$.  Consequently, $\xi_{min}^{(x)}=h(1)$ and $\xi_{max}^{(x)}=h(3)$ as described in \eqref{def:xi}. We know that the points satisfying $h(3)-h(1)\geq 3-1$ are all observable points as illustrated in  Definition \ref{obd}. By simple calculations similar to the proof of Lemma \ref{line-ob}, $x_3\leq0$ can be obtained. Consequently, the observation points $x = (x_1,x_2,x_3)\notin \Gamma$  satisfying $x_3\leq0$ are all observable.  Further more, if $(x-a(t))\cdot a^{\prime}(t)=x_3-(t-2)\leq0$, then  $x_3\in [-2,-1]$. Thus,  the smallest annulus $\Lambda_{\Gamma}^{(x)}$ centered at $x$ and containing the trajectory can be  fully recovered if and only if the observable points $x$ satisfies $-2\leq x_3 \leq -1$. If not, one can only get a slimmer annulus $A_{\Gamma}^{(x)}\subset \Lambda_{\Gamma}^{(x)}$. Moreover, all observation points $x$ where $x_3>0$ are  non-observable. The corresponding numerical results are presented in Figs.\ref{fig:line1}, \ref{fig:line2} and \ref{fig:line3}.

In Fig.\ref{fig:line1}, we examine various observable points $x$, where $\theta\in[0, 2\pi]$ and $\varphi\in[2\pi/3, \pi]$. Specifically, we restrict $\varphi$ to the range $[2\pi/3, \pi]$, which corresponds to $-2\leq x_3\leq-1$. We observe that   $(x-a(t))\cdot a^{\prime}(t)\leq0$ for all $t\in[1,3]$. Consequently, our theoretical predictions ensure that the trajectory of the moving point source can be fully enclosed within the smallest annulus  $\Lambda_{\Gamma}^{(x)}$ centered at $x$. It is worthy noting that our numerical examples demonstrate that $A_{\Gamma}^{(x)}=\Lambda_{\Gamma}^{(x)}$.

In Fig.\ref{fig:line2}, we collect data at the observable points $x=(2 \sin \varphi \cos \theta,  2\sin \varphi \sin \theta, 2 \cos \varphi )$  where $\theta\in[0, 2\pi]$   and $\varphi\in [\pi/2, 2\pi/3)$, such that  $-1< x_3\leq0$. However, we note that $(x-a(t))\cdot a^{\prime}(t)\leq0$ no longer holds for all $t\in[1,3]$.
Even though these observation points $x$ belong to the observable set, the corresponding annulus $A_{\Gamma}^{(x)}$ turn out to be slimmer than the smallest annulus that encloses the trajectory of the moving source and is centered at $x$. This discrepancy arises due to the relationship $A_{\Gamma}^{(x)}\subset \Lambda_{\Gamma}^{(x)}$, as established by Lemma \ref{lem3.9}.

The observation points $x=(2 \sin \varphi \cos \theta,  2\sin \varphi \sin \theta,$ $2 \cos \varphi )$   are non-observable when  $\theta\in[0, 2\pi]$   and $\varphi\in[0,\pi/2)$. Numerical results in Fig.\ref{fig:line3} indicate that the corresponding  indicator values are consistently  much smaller than $10^{-5}$. This is consistent with the outcome of Theorem \ref{Th:factorization}, which implies that it is not possible to reconstruct the annulus centered at $x$ which contains partial or whole information on the trajectory of the moving source. The further the non-observable points are from the observable region, the lower the corresponding indicator values.  Fig.\ref{fig:line3} shows that partial information on the trajectory can still be retrieved by our indicator function even at non-observable points, which is an intriguing observation that requires further investigation.

\vspace{.1in}\textbf{Example 2: An arc in $\R^3$}\vspace{.1in}

\begin{figure}[H]
	\centering
    \subfigure[$\theta=0$, $\varphi=\pi/2$]{
		\includegraphics[scale=0.22]{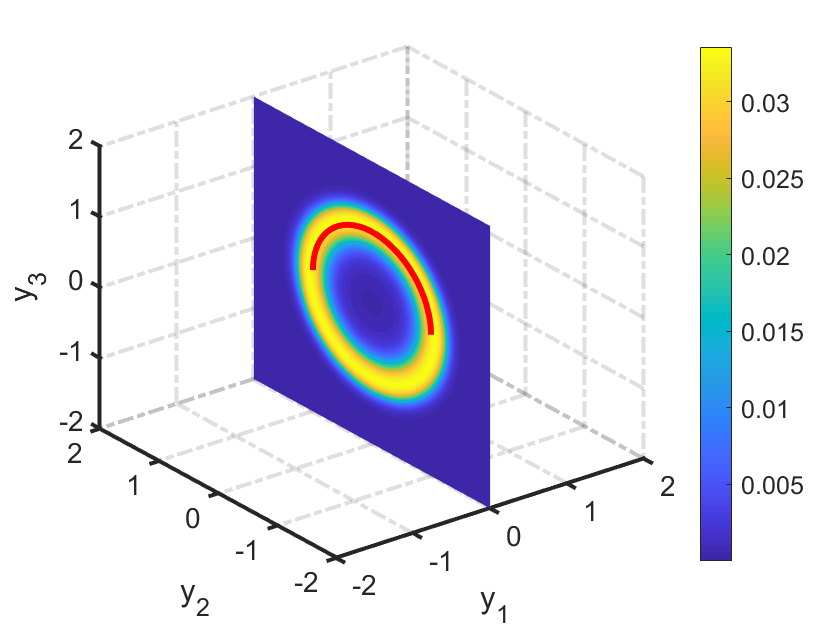}
		
	}
    \subfigure[$\theta=1\pi/8$, $\varphi=\pi/2$]{
		\includegraphics[scale=0.22]{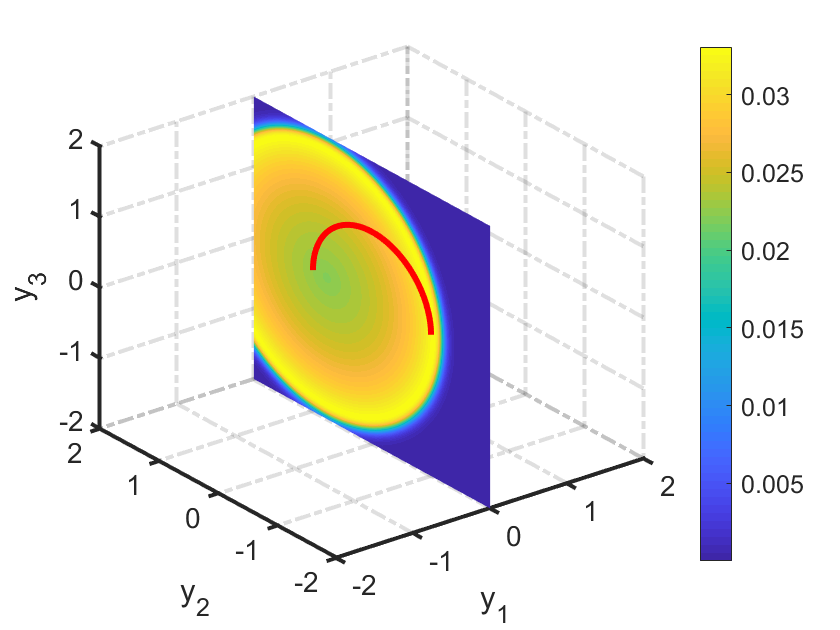}
		
	}
    \subfigure[$\theta=2\pi/8$, $\varphi=\pi/2$]{
		\includegraphics[scale=0.22]{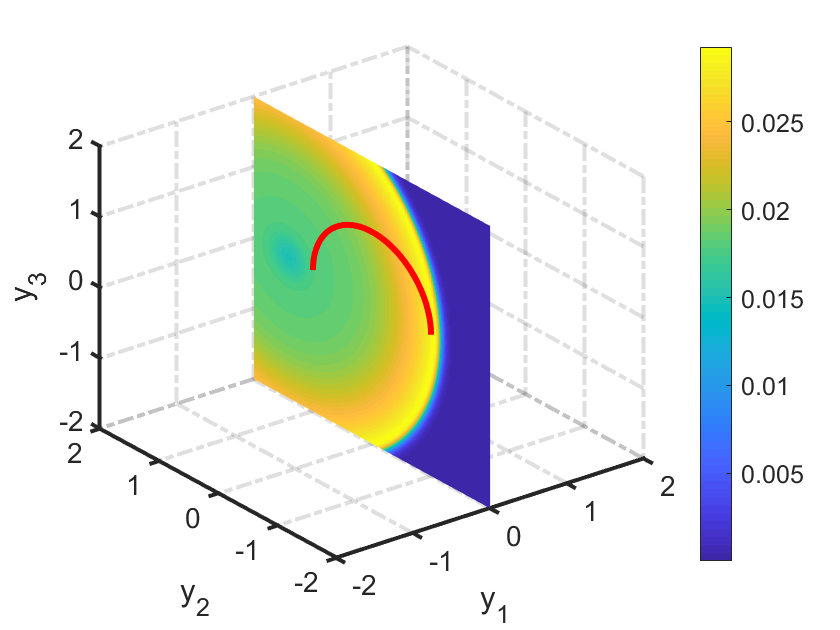}
		
	}
    \subfigure[$\theta=3\pi/8$, $\varphi=\pi/2$]{
		\includegraphics[scale=0.22]{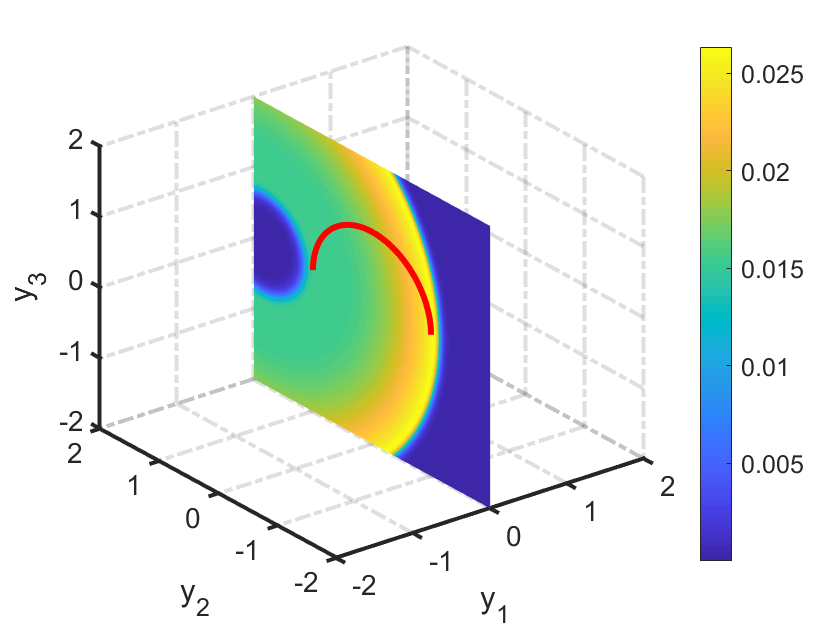}
		
	}
	\subfigure[$\theta=4\pi/8$, $\varphi=\pi/2$ ]{
		\includegraphics[scale=0.22]{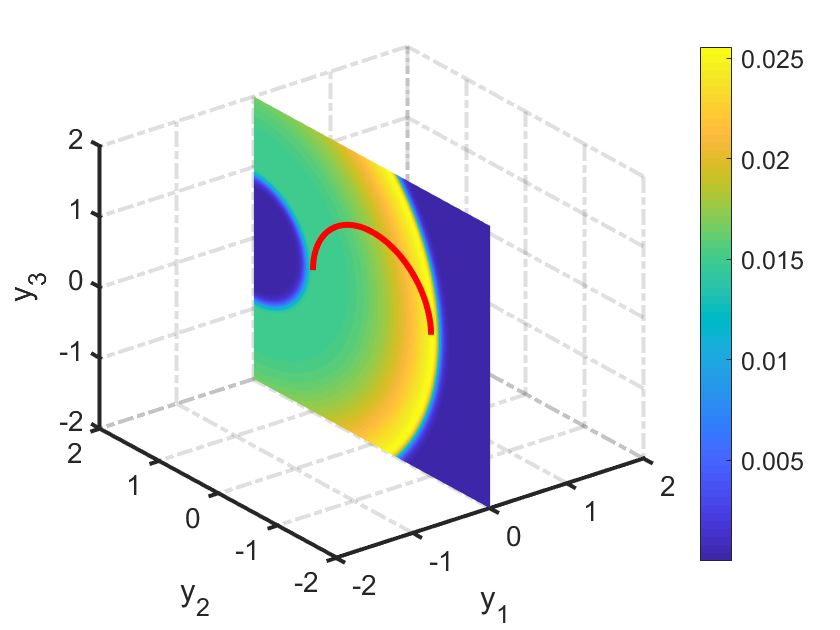}
		
	}
    \subfigure[$\theta=5\pi/8$, $\varphi=\pi/2$  ]{
		\includegraphics[scale=0.22]{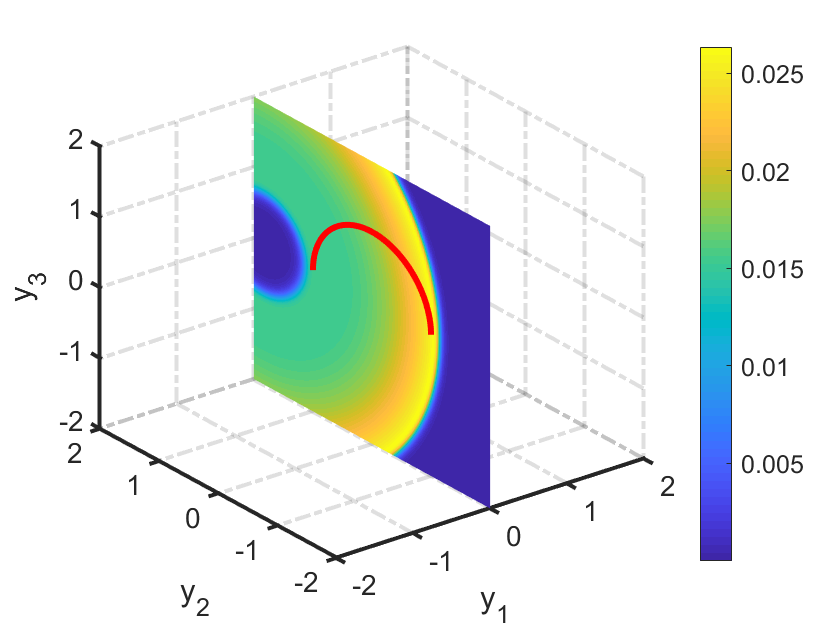}
	
	}
    \subfigure[$\theta=6\pi/8$, $\varphi=\pi/2$ ]{
		\includegraphics[scale=0.22]{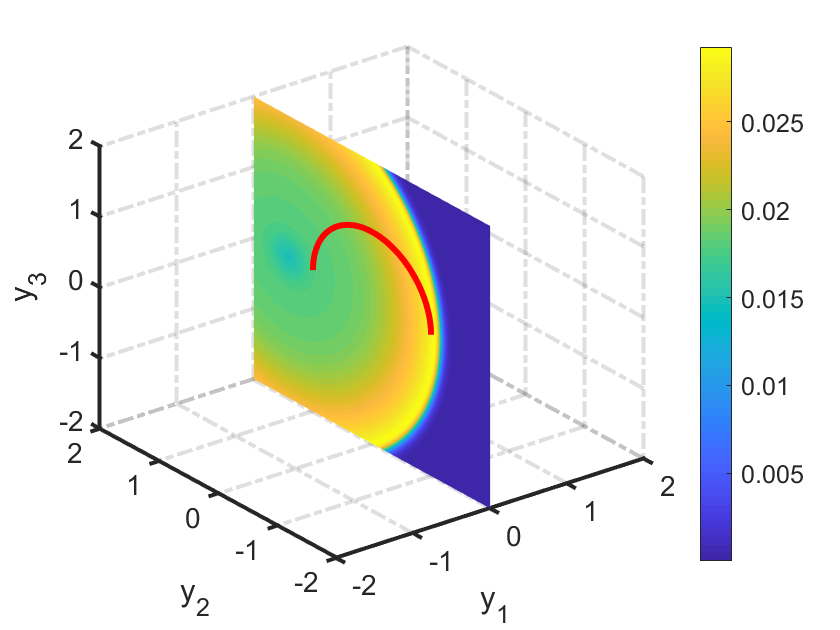}
	}
    \subfigure[$\theta=7\pi/8$, $\varphi=\pi/2$ ]{
		\includegraphics[scale=0.22]{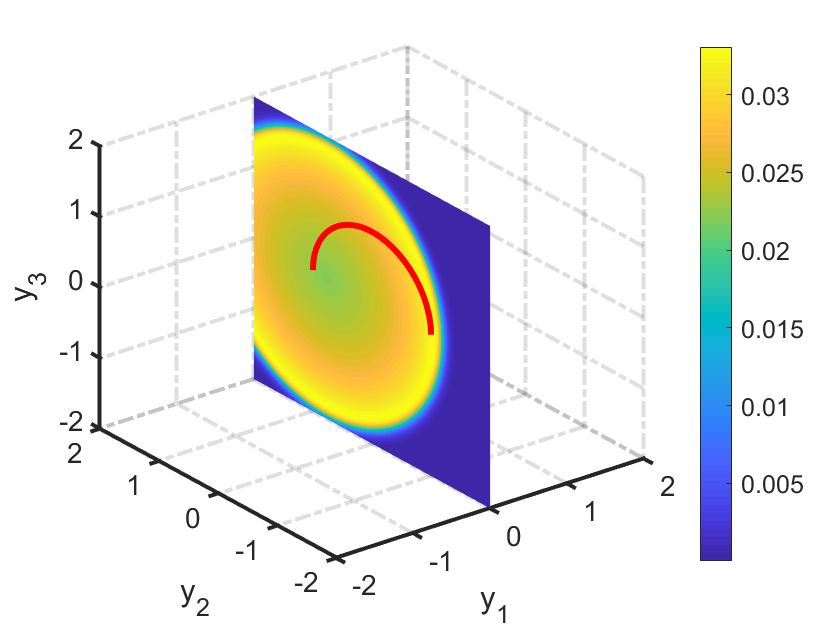}
	}
    \subfigure[$\theta=8\pi/8$, $\varphi=\pi/2$ ]{
		\includegraphics[scale=0.22]{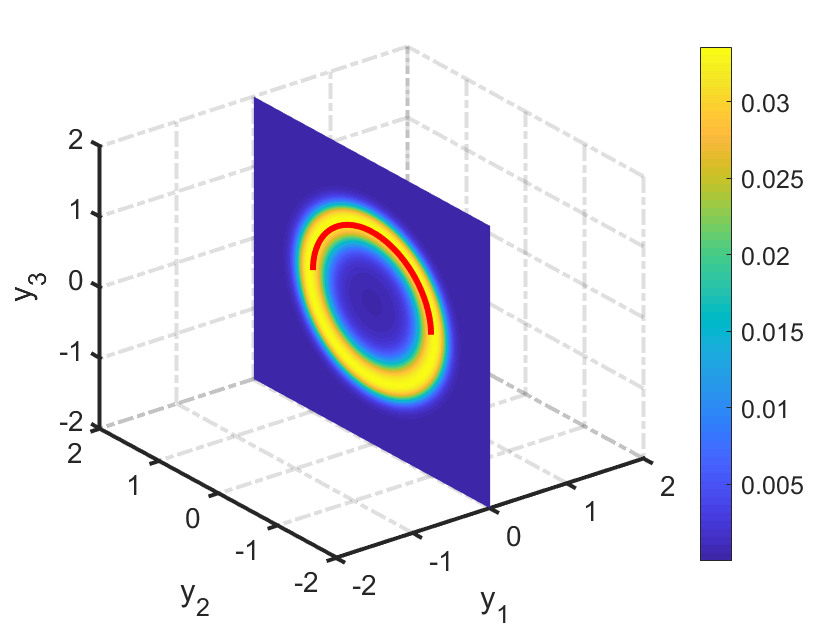}
 	}
\caption{Reconstruction from a single observable point $x=(2 \sin \varphi \cos \theta, 2\sin \varphi \sin \theta, 2 \cos \varphi )$  with $\theta\in[0, \pi]$   and $\varphi=\pi/2$ for an arc  $a(t)=(0,\cos t,\sin t)$ where $t\in[0,\pi]$. Here it holds that $A_{\Gamma}^{(x)}= \Lambda_{\Gamma}^{(x)}$.}\label{fig:circle1}
\end{figure}

\begin{figure}[H]
	\centering
    \subfigure[$\theta=0$, $\varphi=0$]{
		\includegraphics[scale=0.22]{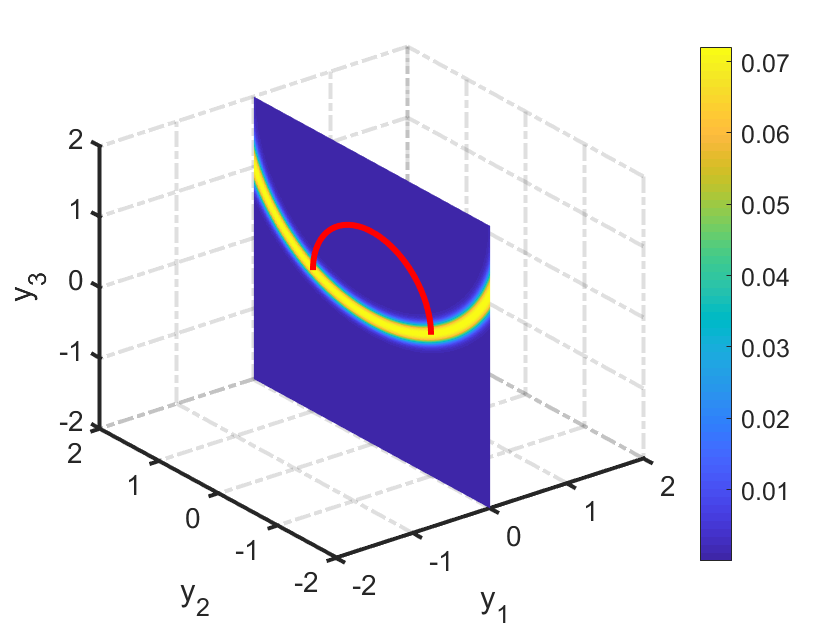}
		
	}
    \subfigure[$\theta=\pi/5$, $\varphi=\pi/4$]{
		\includegraphics[scale=0.22]{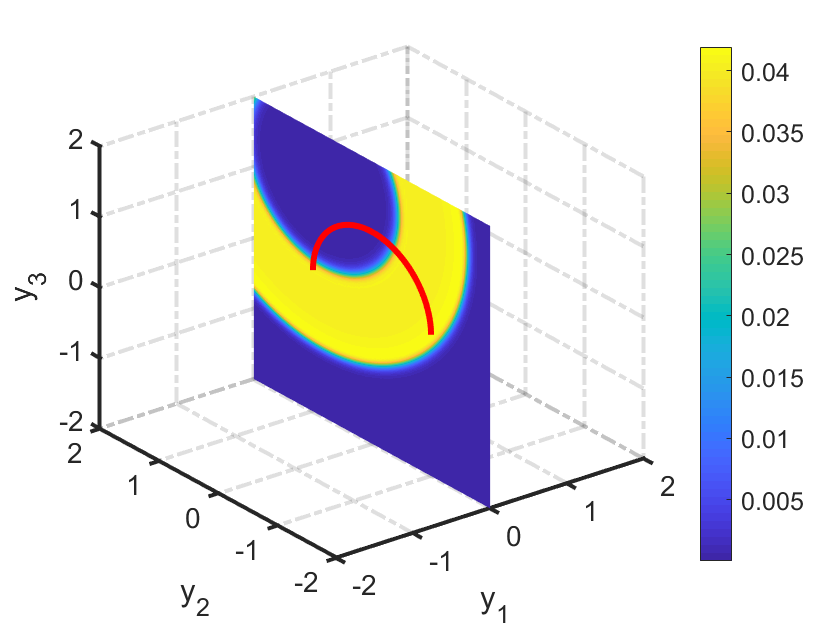}
		
	}
    \subfigure[$\theta=2\pi/5$, $\varphi=\pi/4$]{
		\includegraphics[scale=0.22]{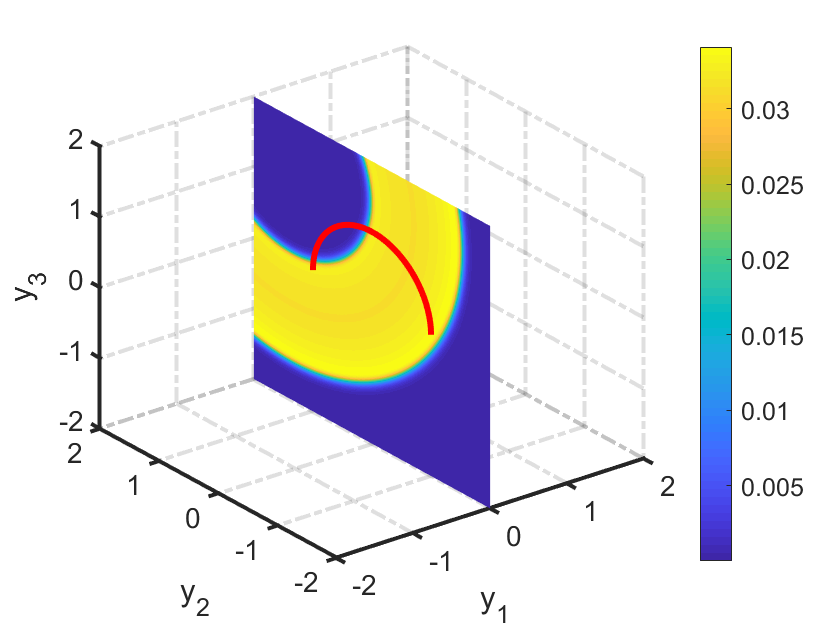}
		
	}
    \subfigure[$\theta=3\pi/5$, $\varphi=\pi/4$]{
		\includegraphics[scale=0.22]{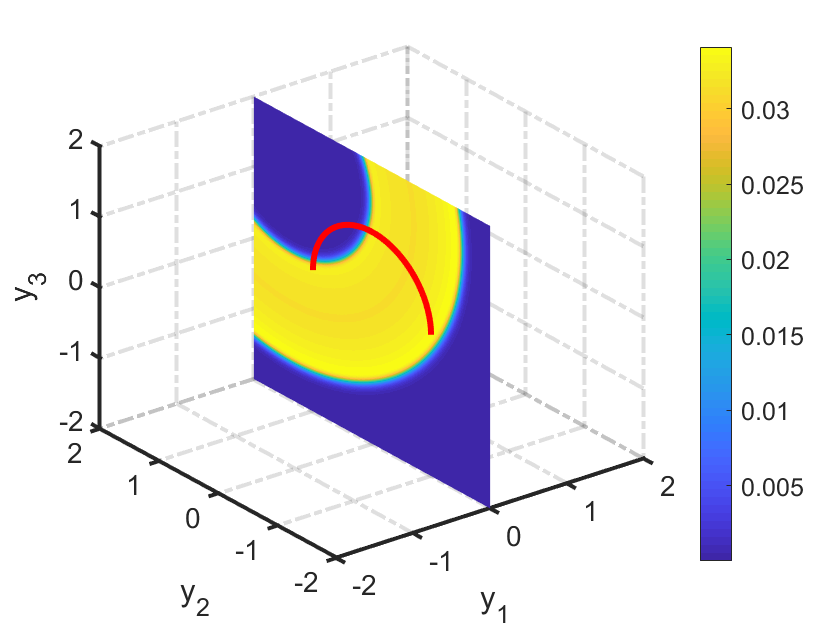}
		
	}
	\subfigure[$\theta=4\pi/5$, $\varphi=\pi/4$ ]{
		\includegraphics[scale=0.22]{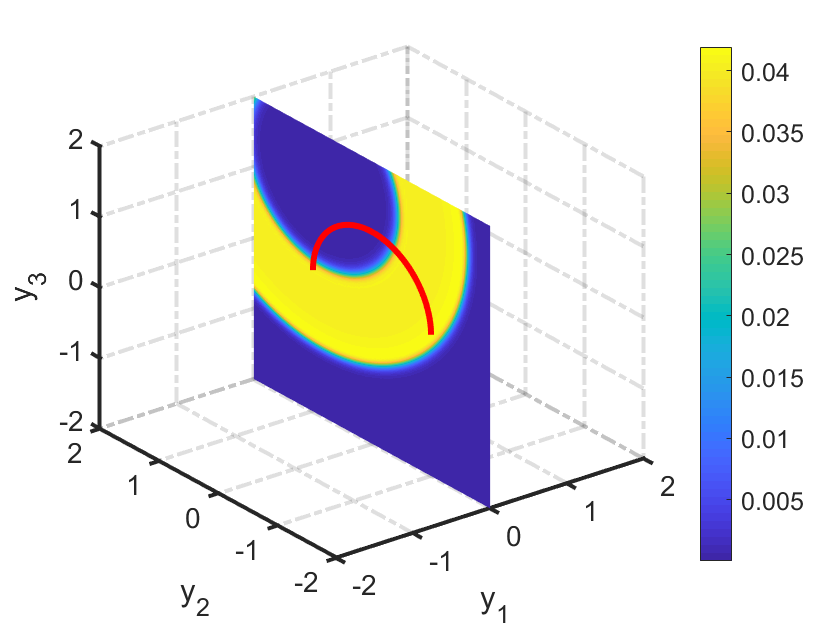}
		
	}
    \subfigure[$\theta=\pi/5$, $\varphi=3\pi/4$  ]{
		\includegraphics[scale=0.22]{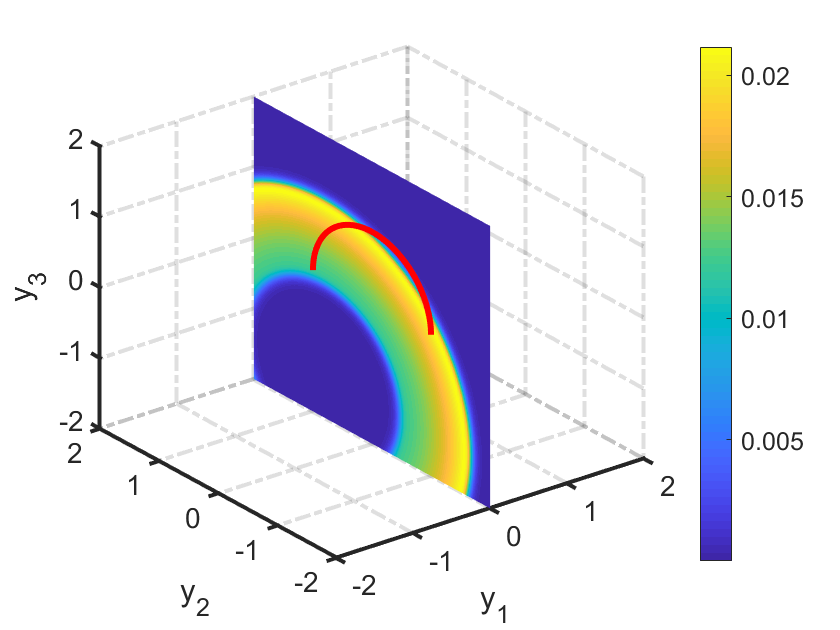}
	
	}
    \subfigure[$\theta=2\pi/5$, $\varphi=3\pi/4$ ]{
		\includegraphics[scale=0.22]{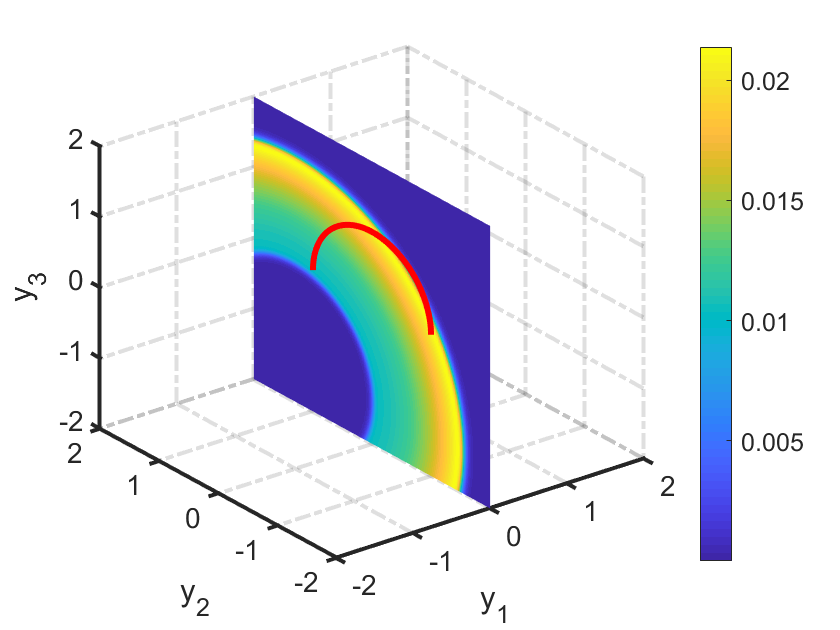}
	}
    \subfigure[$\theta=3\pi/5$, $\varphi=3\pi/4$ ]{
		\includegraphics[scale=0.22]{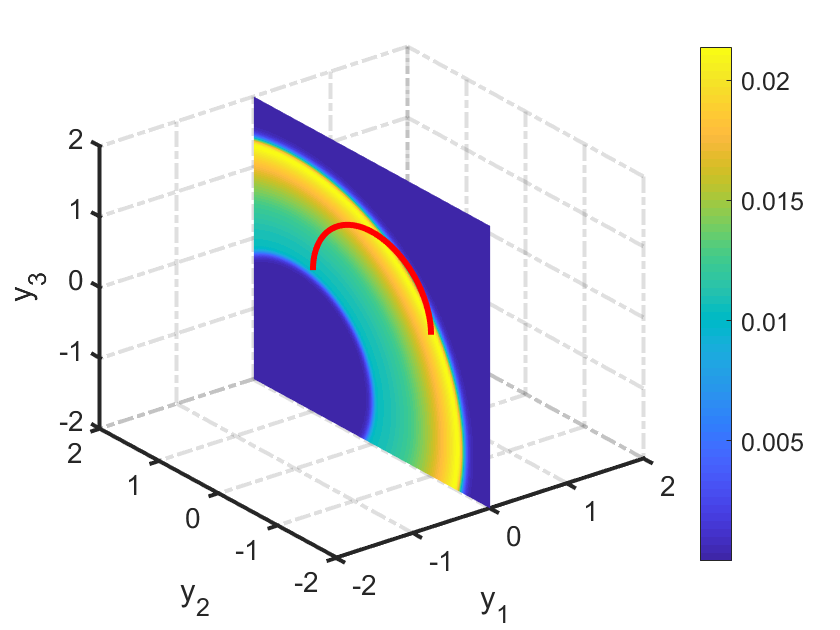}
	}
    \subfigure[$\theta=4\pi/5$, $\varphi=4\pi/4$ ]{
		\includegraphics[scale=0.22]{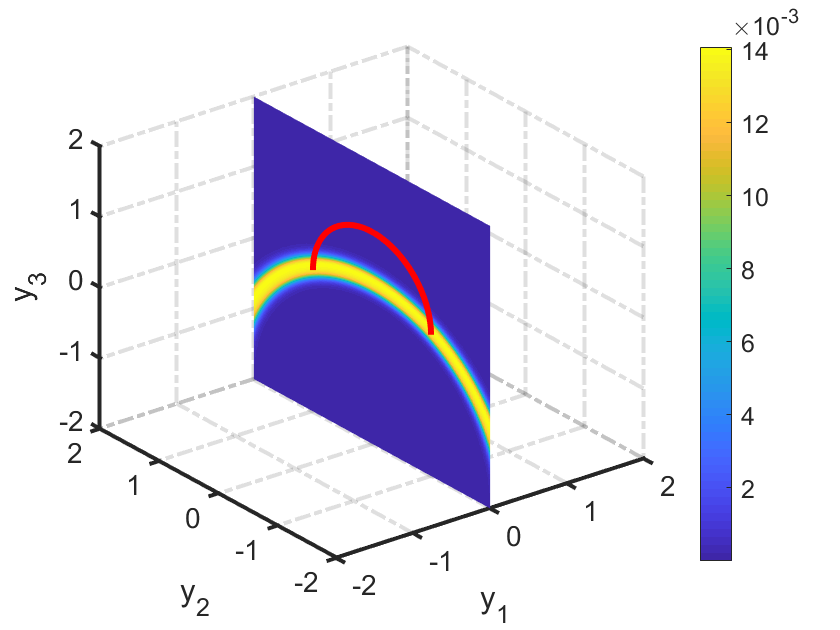}
 	}
\caption{Reconstruction from a single observable point $x=(2 \sin \varphi \cos \theta, 2\sin \varphi \sin \theta, 2 \cos \varphi )$  with $\theta\in[0, \pi]$   and $\varphi=[0,\pi/2)\cup(\pi/2,\pi]$ for an arc  $a(t)=(0,\cos t,\sin t)$ where $t\in[0,\pi]$. Here it holds that $A_{\Gamma}^{(x)}\subset \Lambda_{\Gamma}^{(x)}$. }\label{fig:circle3}
\end{figure}

\begin{figure}[H]
	\centering
    \subfigure[$\theta=5\pi/4$, $\varphi=\pi/4$]{
		\includegraphics[scale=0.22]{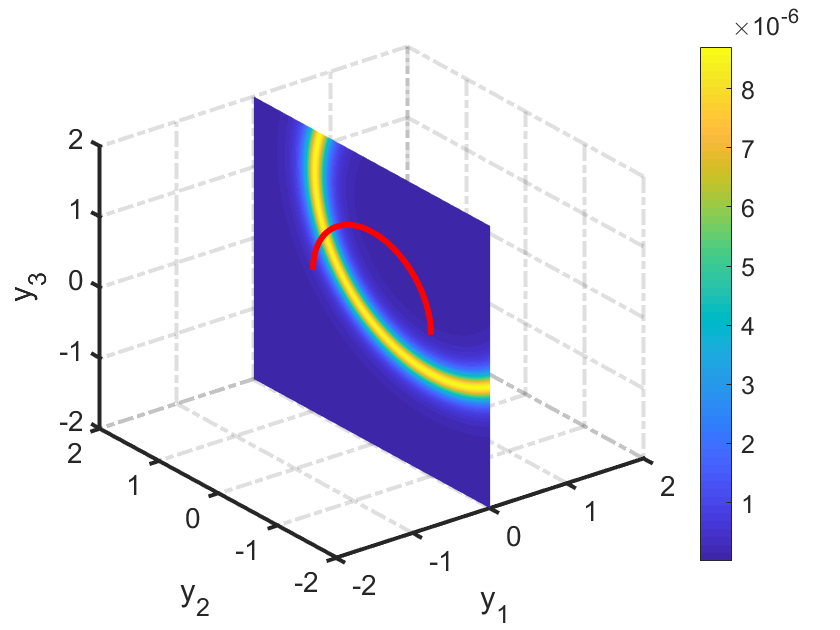}
		
	}
    \subfigure[$\theta=5\pi/4$, $\varphi=2\pi/4$]{
		\includegraphics[scale=0.22]{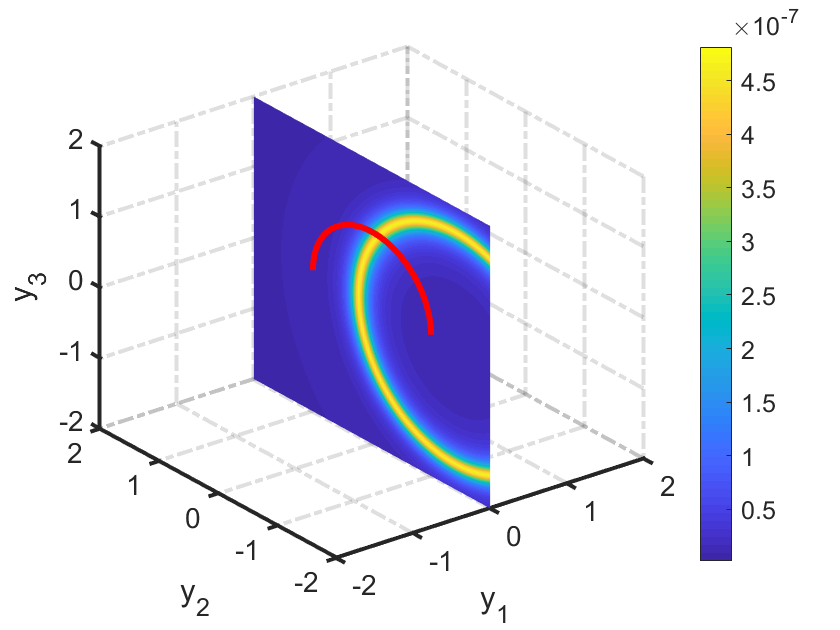}
		
	}
    \subfigure[$\theta=5\pi/4$, $\varphi=3\pi/4$]{
		\includegraphics[scale=0.22]{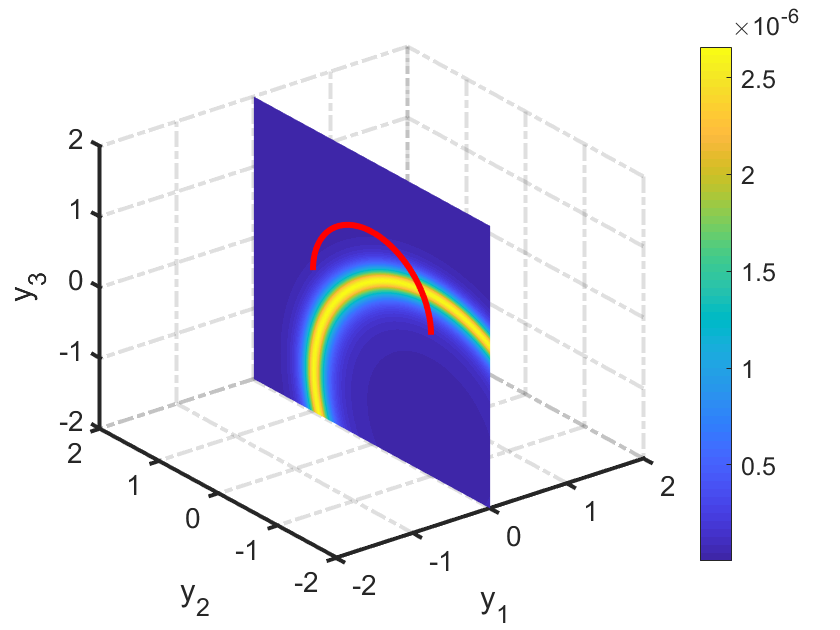}
		
	}
    \subfigure[$\theta=6\pi/4$, $\varphi=\pi/4$]{
		\includegraphics[scale=0.22]{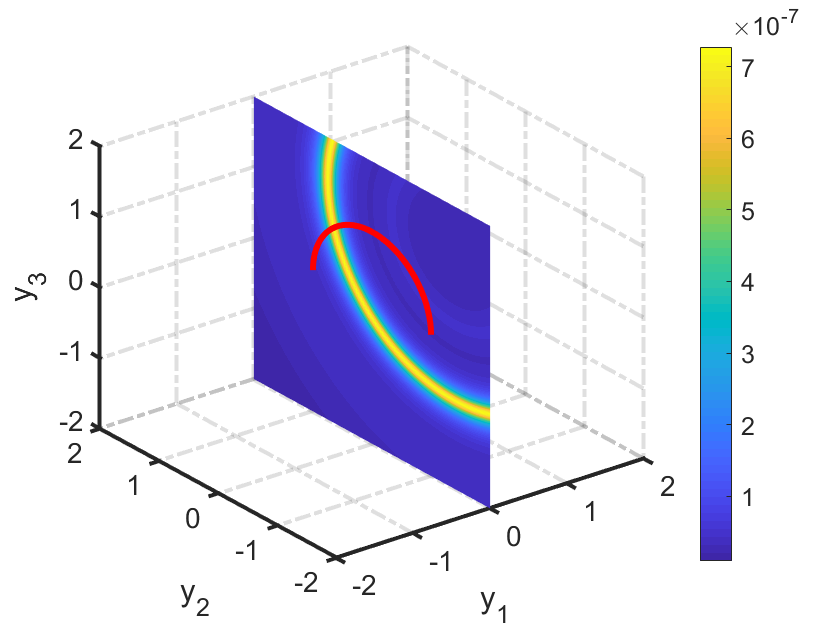}
		
	}
	\subfigure[$\theta=6\pi/4$, $\varphi=2\pi/4$ ]{
		\includegraphics[scale=0.22]{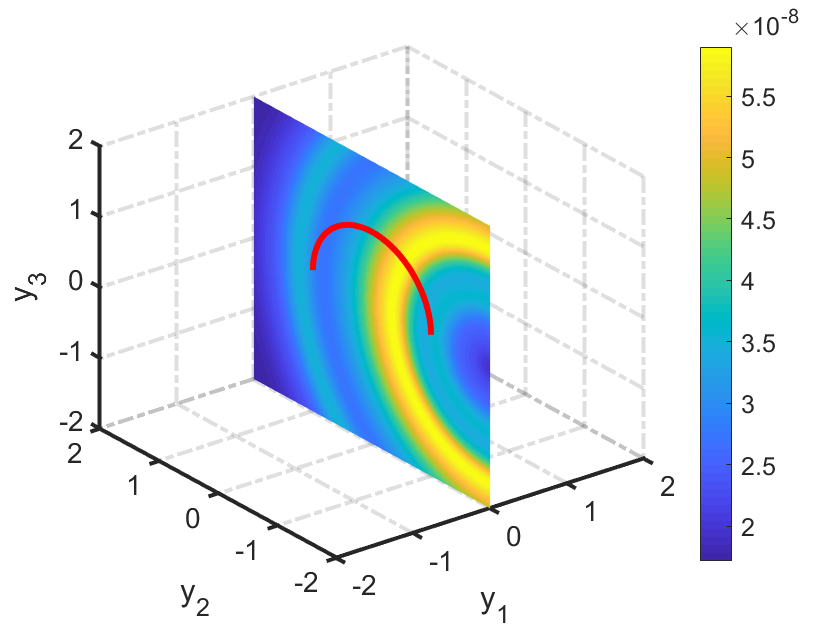}
		
	}
    \subfigure[$\theta=6\pi/4$, $\varphi=3\pi/4$  ]{
		\includegraphics[scale=0.22]{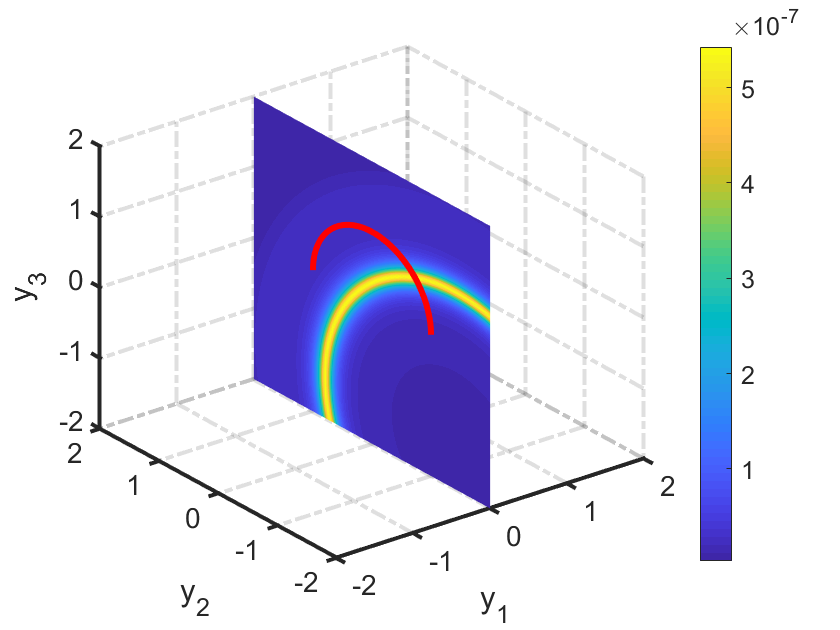}
	
	}
    \subfigure[$\theta=7\pi/4$, $\varphi=\pi/4$ ]{
		\includegraphics[scale=0.22]{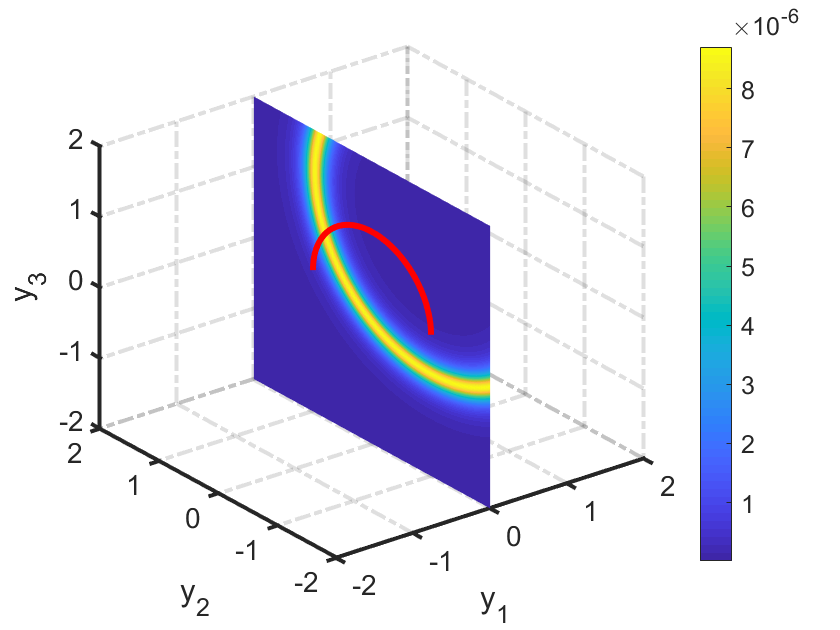}
	}
    \subfigure[$\theta=7\pi/4$, $\varphi=2\pi/4$ ]{
		\includegraphics[scale=0.22]{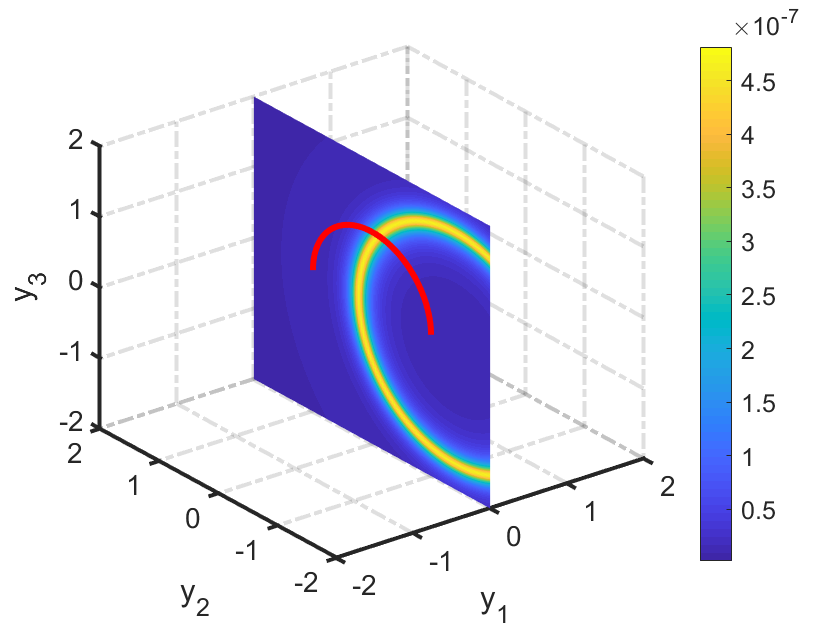}
	}
    \subfigure[$\theta=7\pi/4$, $\varphi=3\pi/4$ ]{
		\includegraphics[scale=0.22]{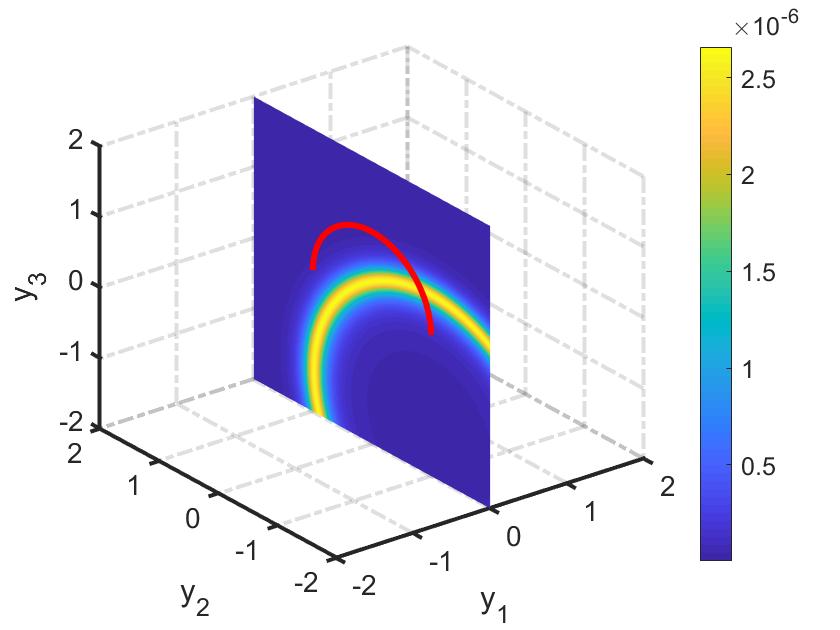}
	}
\caption{Reconstruction from a single non-observable point $x=(2 \sin \varphi \cos \theta, 2\sin \varphi \sin \theta,$  $2 \cos \varphi )$  with $\theta\in(\pi, 2\pi)$   and $\varphi\in[0,\pi]$ for an arc  $a(t)=(0,\cos t,\sin t)$ where $t\in[0,\pi]$.}\label{fig:circle4}
\end{figure}

As  demonstrated in Example 2 of Section 3, we consider the trajectory of the moving point given by $a(t)=(0, \cos t, \sin t)$ where $t\in[0, \pi]$.  From the orbit function, we obtain
\ben
&&h(t)=t+|x-a(t)|=t+\sqrt{x_1^2+(x_2-\cos t)^2+(x_3-\sin t)^2},\\
&&h^{\prime}(t)=1-\frac{x-a(t)}{|x-a(t)|}\cdot a^{\prime}(t)=1-\frac{-x_2 \sin t +x_3 \cos t}{|x-a(t)|}.
\enn
It is evident that $|a^{\prime}(t)| = 1$, thus  $h^{\prime}(t) \geq0$ for all $t\in[0,\pi]$, which indicates the function $h(t)$ monotonically increases over $[0,\pi]$. Subsequently, we have
$\xi_{min}^{(x)}=h(0)$ and $\xi_{max}^{(x)}=h(\pi)$. According to Definition \ref{obd}, observable points are those that satisfy $h(\pi)-h(0)\geq \pi-0$. By simple calculations similar to the proof of Lemma \ref{lem:circle} one can  infer that  $x_2\geq0$. Therefore, the observation points $x = (x_1,x_2,x_3)\notin \Gamma$  that satisfy $x_2\geq0$ are all observable.
Furthermore,  the statement $(x-a(t))\cdot a^{\prime}(t)=-x_2 \sin t+x_3\cos t\leq0$ is equivalent to
\ben
-x_2 \sin t +x_3 \cos t&=&\sqrt{x_2^2+x_3^2} \left( \frac{x_3}{\sqrt{x_2^2+x_3^2}} \cos t-\frac{x_2}{\sqrt{x_2^2+x_3^2}}\sin t\right)\\
&=&\sqrt{x_2^2+x_3^2} \sin(\alpha-t)\leq0,
\enn
where $\sin \alpha =\frac{x_3}{\sqrt{x_2^2+x_3^2}}$ and $ \cos \alpha =\frac{x_2}{\sqrt{x_2^2+x_3^2}}$.  If  it holds that $-x_2 \sin t +x_3 \cos t \leq 0$ for all $t\in[0,\pi],$ then it is evident that $\alpha=2n\pi, n=0,\pm1,\pm2, ...$, meaning $x_3=0$. The smallest annulus $\Lambda_{\Gamma}^{(x)}$, centered at $x$ and containing the trajectory of the moving source, is recoverable only when the observable points $x$ satisfy $x_3=0$.  Otherwise, a slimmer annulus $A_{\Gamma}^{(x)}\subset \Lambda_{\Gamma}^{(x)}$ can be obtained. Additionally, all observation points $x$ with $x_2<0$ are  non-observable. The numerical results are presented in Figs.\ref{fig:circle1}, \ref{fig:circle3} and \ref{fig:circle4}.

\begin{figure}[H]
	\centering
    \subfigure[$\theta=0$, $\varphi=\pi/2$]{
		\includegraphics[scale=0.22]{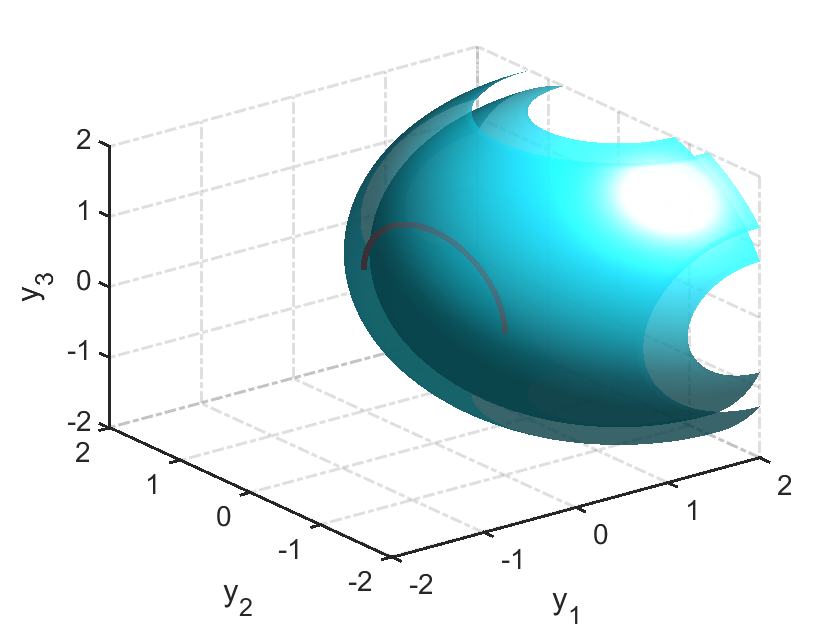}
		
	}
    \subfigure[$\theta=\pi/8$, $\varphi=\pi/2$]{
		\includegraphics[scale=0.22]{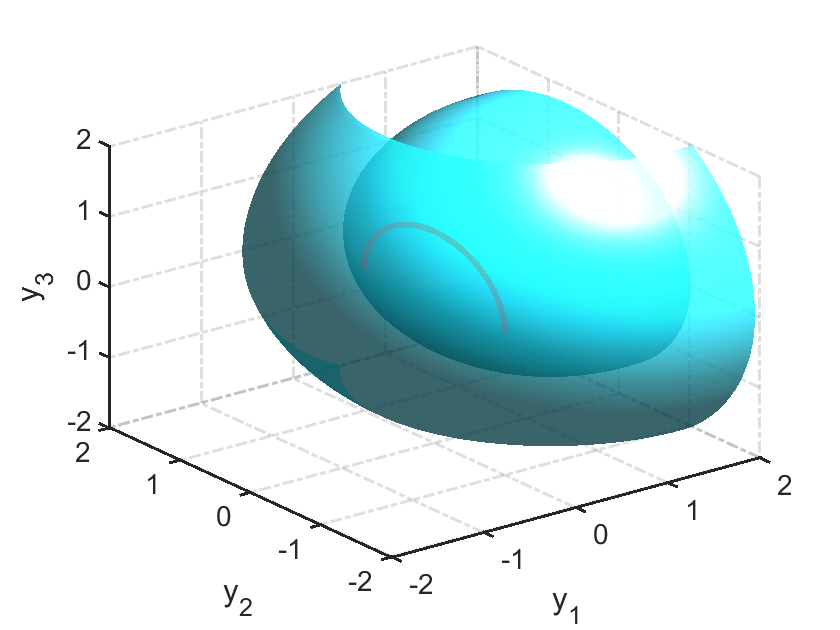}
		
	}
    \subfigure[$\theta=2\pi/8$, $\varphi=\pi/2$]{
		\includegraphics[scale=0.22]{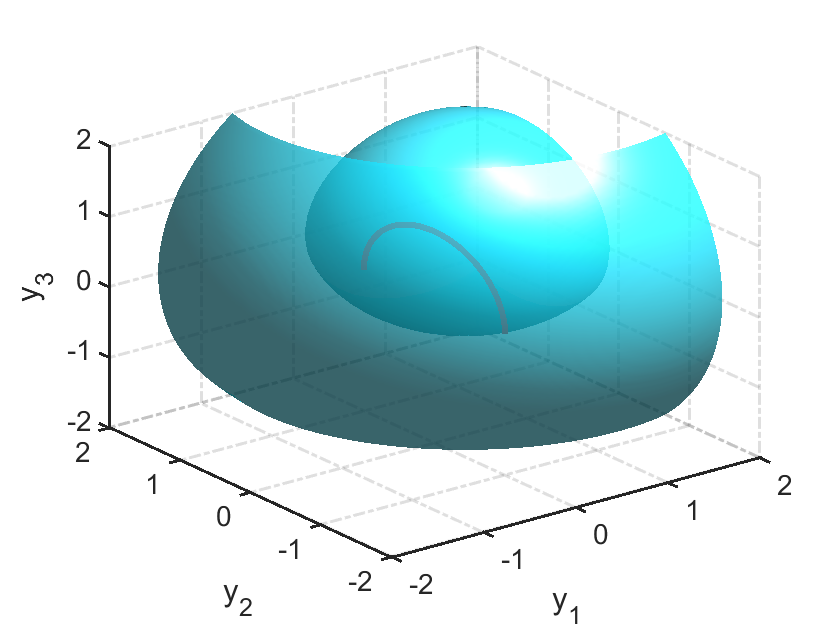}
		
	}
%
%
%
    \subfigure[$\theta=6\pi/8$, $\varphi=\pi/2$ ]{
		\includegraphics[scale=0.22]{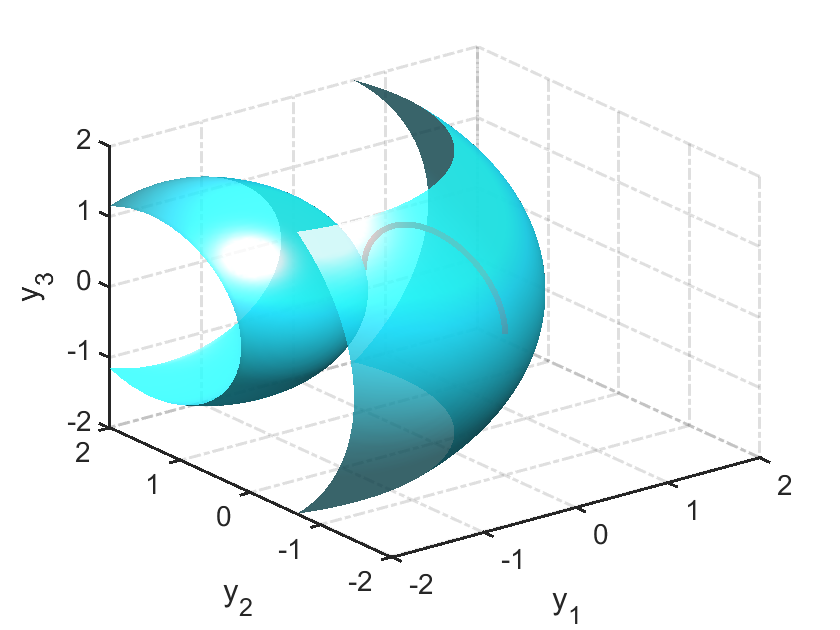}
	}
    \subfigure[$\theta=7\pi/8$, $\varphi=\pi/2$ ]{
		\includegraphics[scale=0.22]{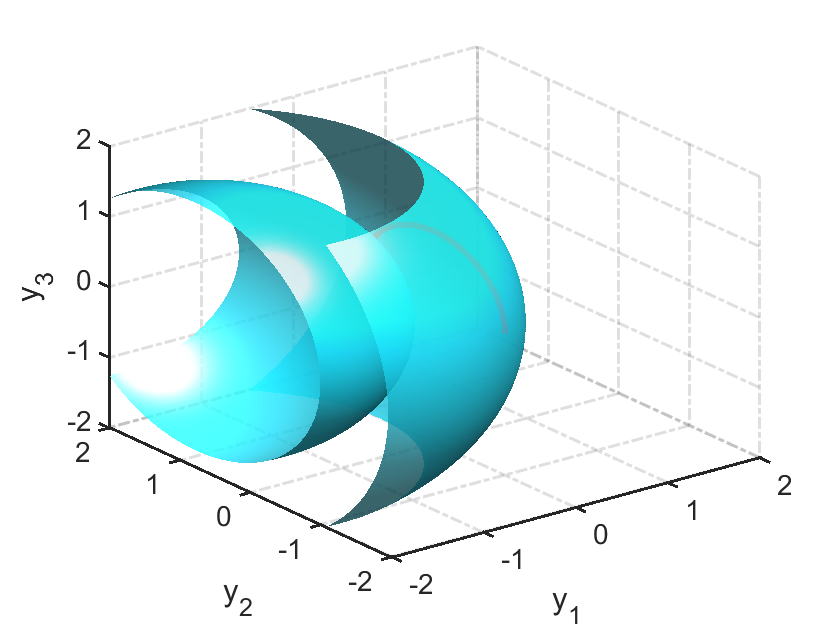}
	}
    \subfigure[$\theta=8\pi/8$, $\varphi=\pi/2$ ]{
		\includegraphics[scale=0.22]{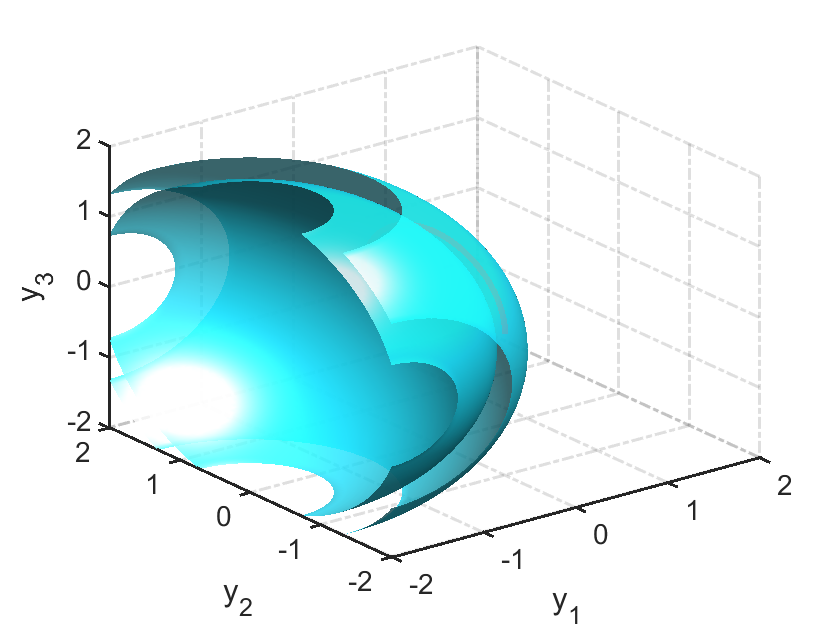}
 	}
\caption{Iso-surfaces of reconstruction from a single observable point $x=(2 \sin \varphi \cos \theta, 2\sin \varphi \sin \theta, 2 \cos \varphi )$  with $\theta\in[0, \pi]$   and $\varphi=\pi/2$ for an arc  $a(t)=(0,\cos t,\sin t)$ where $t\in[0,\pi]$. Here it holds that $A_{\Gamma}^{(x)}= \Lambda_{\Gamma}^{(x)}$.}\label{fig:circle2}
\end{figure}

Fig.\ref{fig:circle1} demonstrates the reconstruction of an arc from  observable points $x=(2 \sin \varphi \cos \theta, 2\sin \varphi \sin \theta, 2 \cos \varphi )$  with $\theta\in[\pi, 2\pi]$ and $\varphi=\pi/2$. We conclude that the trajectory of the moving point source perfectly lies in the smallest annulus centered at $x$ and containing its trajectory. This is due to  the selection of observable points $x$ with $x_3=0$ and $x_2\geq0$, such that $(x-a(t))\cdot a^{\prime}(t)\leq0$. Here, we have $A_{\Gamma}^{(x)}= \Lambda_{\Gamma}$. This effectively demonstrates the effectiveness of our algorithm for imaging an arc in $\R^3$. It is worthy noting that although the arc trajectory of the moving source lies in the annulus depicted in subfigures (a),(b),(c),(g),(h) and (i), they can not be seen clearly since the  slice is set at $y_1=0$. Therefore, corresponding isosurfaces of the reconstruction are plotted in Fig.\ref{fig:circle2}, showing the trajectory of the moving source perfectly located between the isosurfaces as predicted by our theory.

Fig.\ref{fig:circle3} displays the reconstructed annulus $A^{(x)}_{\Gamma}$ which are slimmer than $\Lambda_{\Gamma}^{(x)}$. This is due to our selection of observable points $x=(2 \sin \varphi \cos \theta, 2\sin \varphi \sin \theta, 2 \cos \varphi )$  with $\theta\in[0, \pi]$   and $\varphi=[0,\pi/2)\cup (\pi/2,\pi]$ implying $x_2\geq0$ and $x_3\neq0$, making it unsuitable to apply $(x-a(t))\cdot a^{\prime}(t)\leq0$ for all $t\in [0,\pi]$. This results in $A_{\Gamma}^{(x)}\subset\Lambda_{\Gamma}^{(x)}$, limiting the retrieval of partial trajectory information. Since $h^{\prime}(t)>0$ for $t \in [0,\pi]$, it is possible to capture the starting and end points of the trajectory by $A_{\Gamma}^{(x)}=\{y\in \mathbb{R}^3: |x-a(t_{\max})| \leq |x-y| \leq |x-a(t_{\min})| \}$. It should be noted that  the size of the annulus depends on the location of the observation points. The numerical results  presented in Fig.\ref{fig:circle3}  are in agreement with our theory predictions.

Fig.\ref{fig:circle4} presents indicator functions for  various  non-observable points $x=(2 \sin \varphi \cos \theta, 2\sin \varphi \sin \theta, 2 \cos \varphi )$  with $\theta\in(\pi, 2\pi)$   and $\varphi\in[0,\pi]$, i.e., $x_2<0$. It can be observed that the values of the indicator functions are significantly smaller than $10^{-5}$.

%
%
%
%
%
%
%

\subsection{Sparse observation points}

In this subsection, we extend Examples 1 and 2 to include multi-frequency near-field data measured at sparse points.  To truncate the indicator function (\ref{W}), we introduce the following expression:
\be \label{W1}
W(y):=\left[\sum _{j=1}^M\sum_{n=1}^N\frac{\left| \phi^{(x^{(j)})}_{y}\cdot \overline{\psi_n^{(x^{(j)})} }\right|^2}{ |\lambda_n^{(x^{(j)})}|}\right]^{-1}, \quad y\in \R^3.
\en
In this definition, $M>0$ denotes the number of sparse observation points distributed on $S_2$. Also,  the test function $\phi_y^{(x^{(j)})}$  has the same definition as in (\ref{testn}), and
 $\left\{(\lambda_n^{(x^{(j)})} , \psi_n^{(x^{(j)})}): n=1,\cdots,N\right\} $ denotes an eigensystem of the operator $(\mathcal N^{(x^{(j)})})_\#$.
Notably, $x^{(j)}$ ($j=1,2,\cdots, M$) may contain both observable and non-observable points.  To eliminate the terms similar to
$$\tilde u_j=\sum_{n=1}^N\frac{\left| \phi^{(x^{(j)})}_{y} \cdot \overline{\psi_n^{(x^{(j)})} }\right|^2}{ |\lambda_n^{(x^{(j)})}|},\quad j=1,2,...W$$
 from the sum in (\ref{W1}), we set a threshold $M'>0$. Precisely, if $\min(\tilde u_j (y))>M'$, the point $x^{(j)}$ can be categorized as a non-observable point through the second assertion of Theorem \ref{Th:factorization}.

\begin{figure}[H]
	\centering
    \subfigure[A slice at $y_2=0$]{
		\includegraphics[scale=0.3]{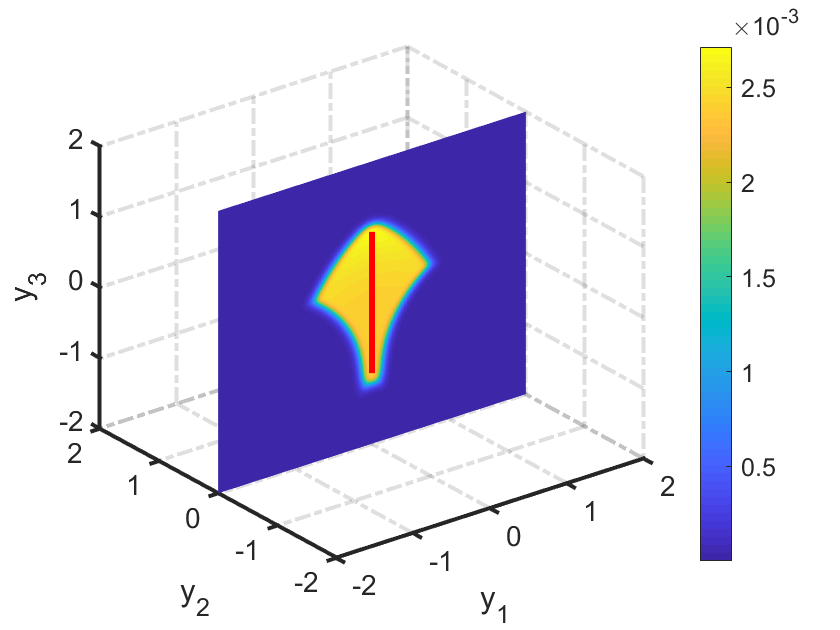}
		
	}
   \subfigure[ Iso-surface level $=1\times10^{-3}$]{
		\includegraphics[scale=0.3]{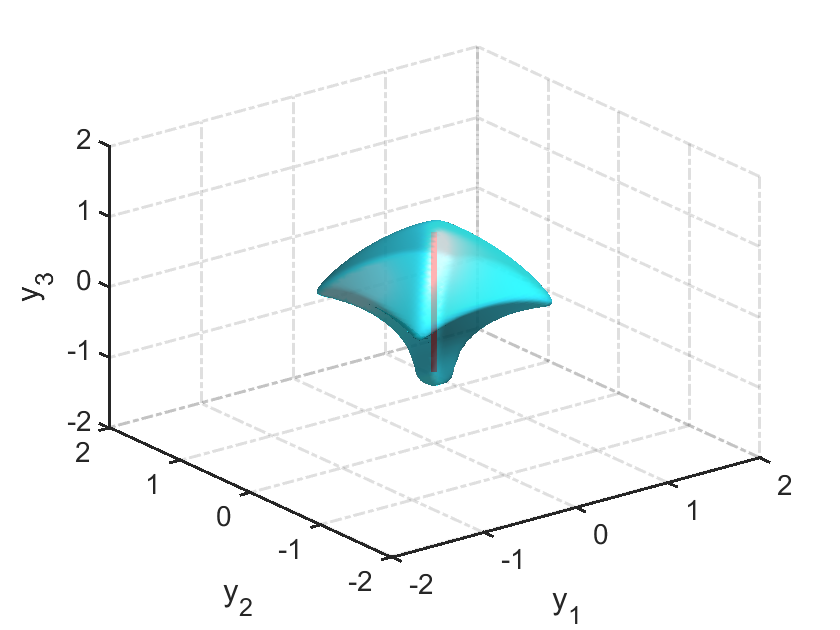}
		
	}
\caption{Reconstruction from sparse observation points $x=(2 \sin \varphi \cos \theta, 2\sin \varphi \sin \theta,$  $2 \cos \varphi )$ with $\theta\in[0, 2\pi)$   and $\varphi\in[2\pi/3,\pi]$ for a straight line segment  $a(t)=(0,0, t-2)$ where $t\in[1,3]$. Here $M=13$ denotes the number of the observation points and we take $\theta=(j-1)\pi/2$, $j=1,\cdots,4$ and $\varphi=(j+5)\pi/9$, $j=1,\cdots,4$ such that $A_{\Gamma}^{(x_j)}= \Lambda_{\Gamma}^{(x_j)}$.  }\label{fig:line-good-1}
\end{figure}

\begin{figure}[H]
	\centering
    \subfigure[A slice at $y_1=0$]{
		\includegraphics[scale=0.3]{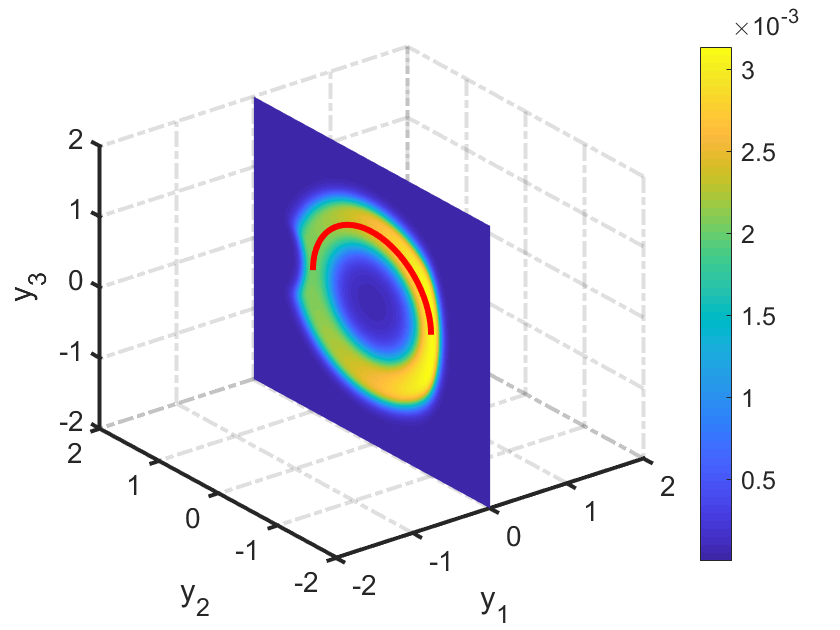}
		
	}
   \subfigure[ Iso-surface level $=1\times10^{-3}$]{
		\includegraphics[scale=0.3]{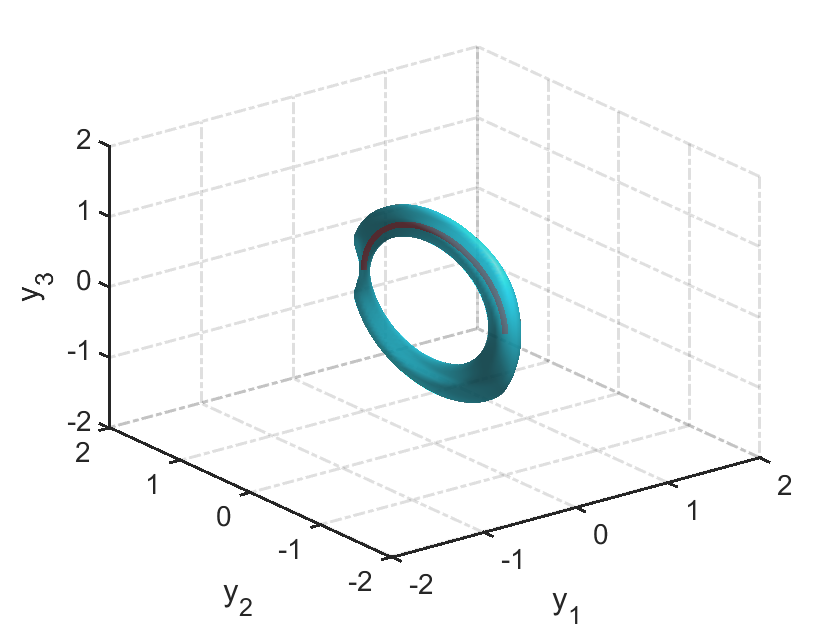}
		
	}
\caption{Reconstruction from sparse observation points $x=(2 \sin \varphi \cos \theta, 2\sin \varphi \sin \theta,$  $2 \cos \varphi )$  with $\theta\in[0, \pi]$   and $\varphi=\pi/2$ for an arc $a(t)=(0,\cos t, \sin t)$ where $t\in[0,\pi]$. Here $M=9$ denotes the number of the observation points and we take $\theta=(j-1)\pi/8$, $j=1,\cdots,M$ and $\varphi=\pi/2$ such that $A_{\Gamma}^{(x_j)}= \Lambda_{\Gamma}^{(x_j)}$.  }\label{fig:circle-good-2}
\end{figure}

Firstly, assuming that all the selected observation points are observable, and that the angle between the vector connecting these observable points and the trajectory points, and the velocity vector of the moving point source, lies within the range of $[\pi/2,3\pi/2]$, implying $(x-a(t))\cdot a^{\prime}(t)\leq0.$  For every observation point,  it is possible to extract the smallest annulus centered at the observation point and containing the trajectory of the moving point source. In Figs.\ref{fig:line-good-1} and \ref{fig:circle-good-2}, we use 13 observation points to reconstruct a straight line segment $a(t)=(0,0, t-2)$ with $t\in[1,3]$ in Example 1 and 9 observation points to reconstruct an arc  $a(t)=(0,\cos t, \sin t)$ with $t\in[0,\pi]$ in Example 2, respectively. The reconstructed slices  and iso-surface are shown in Figs.\ref{fig:line-good-1} and \ref{fig:circle-good-2}, where it is evident that the trajectories are enclosed by the intersections of the smallest annulus $\Lambda_{\Gamma}^{(x_j)}$ centered at $x_j$ and containing their own trajectories. However, since we have only chosen a partial set of observation points, we are not able to reconstruct the trajectories perfectly.

\begin{figure}[H]
	\centering
    \subfigure[$M=4$]{
		\includegraphics[scale=0.22]{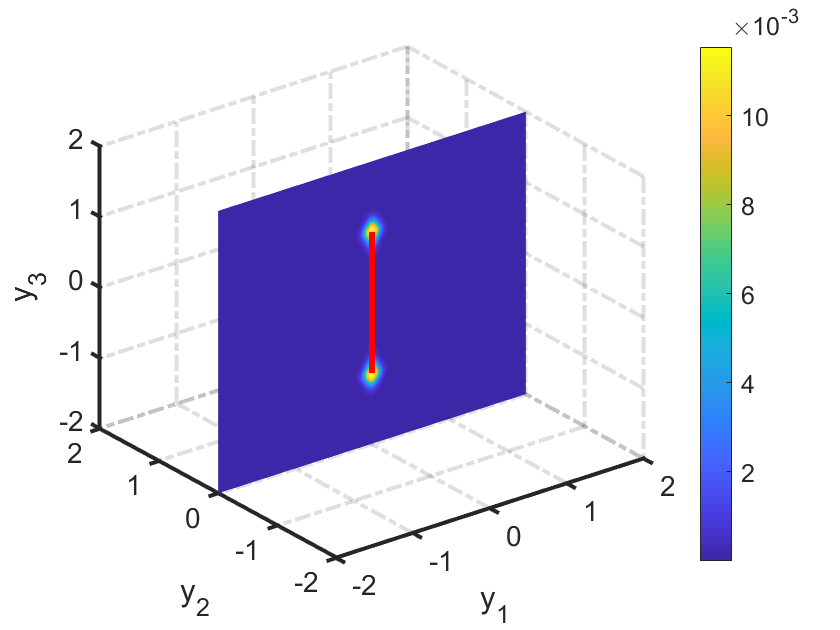}
		
	}
    \subfigure[$M=6$]{
		\includegraphics[scale=0.22]{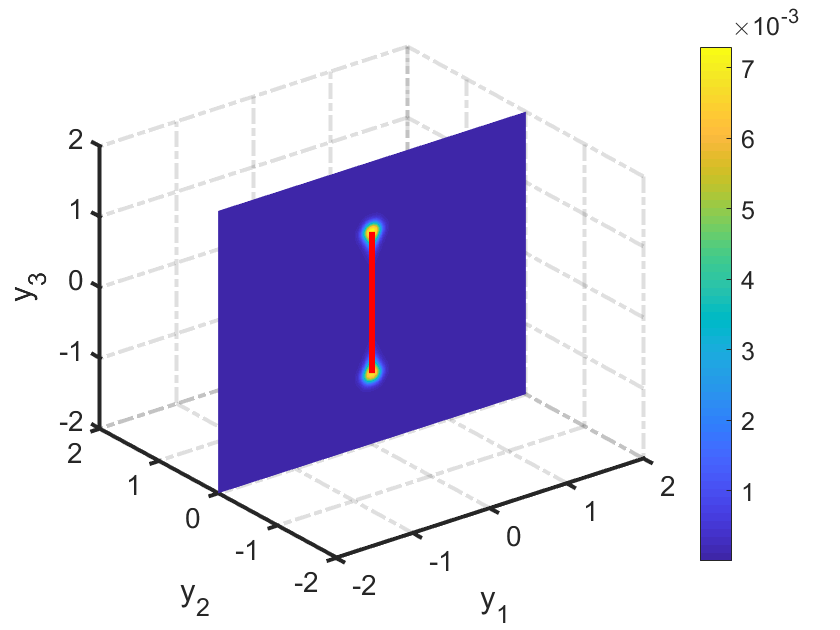}
		
	}
   \subfigure[$M=12$ ]{
		\includegraphics[scale=0.22]{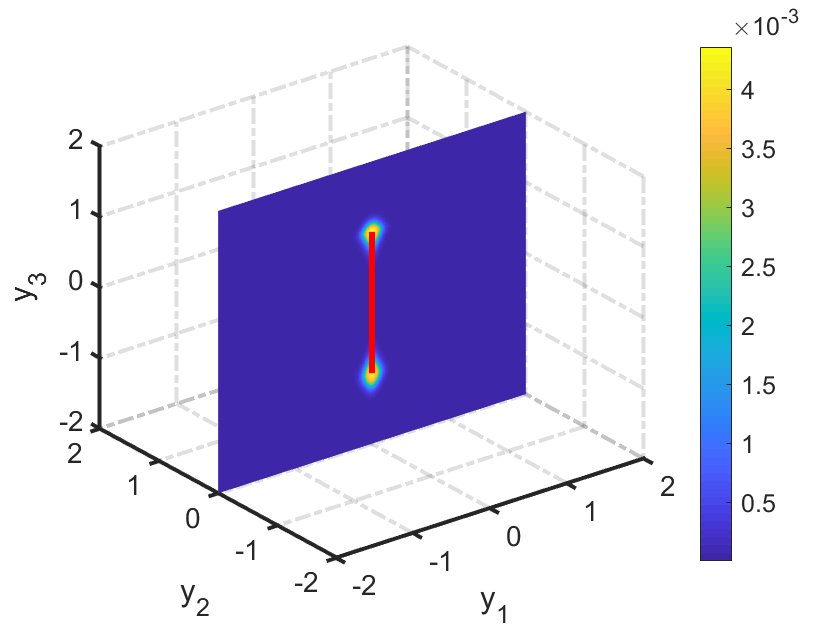}
		
	}
\caption{Reconstruction from sparse observation points $x=(2 \sin \varphi \cos \theta, 2\sin \varphi \sin \theta,$  $2 \cos \varphi )$  with $\theta\in[0, 2\pi)$   and $\varphi\in[0,\pi]$ for a straight line segment $a(t)=(0,0,t-2)$ where $t\in[1,3]$. Here $M$ denotes the number of the observation points. (a) $\theta=(j-1)*2\pi/4$, $j=1,\cdots,4$ and $\varphi=\pi/2$;  (b)  $\varphi=0$; $\theta=(j-1)*2\pi/4$, $j=1,\cdots,4$ and $\varphi=\pi/2$;  $\varphi=\pi$;  (c) $\theta=(2j-1)*\pi/4$, $j=1,\cdots,4$, $\varphi=\pi/3$; $\theta=(j-1)*2\pi/4$, $j=1,\cdots,4$, $\varphi=\pi/2$;  $\theta=(4j-3)*\pi/8$, $j=1,\cdots,5$, $\varphi=2\pi/3$.}\label{fig:line-m-1}
\end{figure}

\begin{figure}[H]
	\centering
    \subfigure[$M=4$]{
		\includegraphics[scale=0.22]{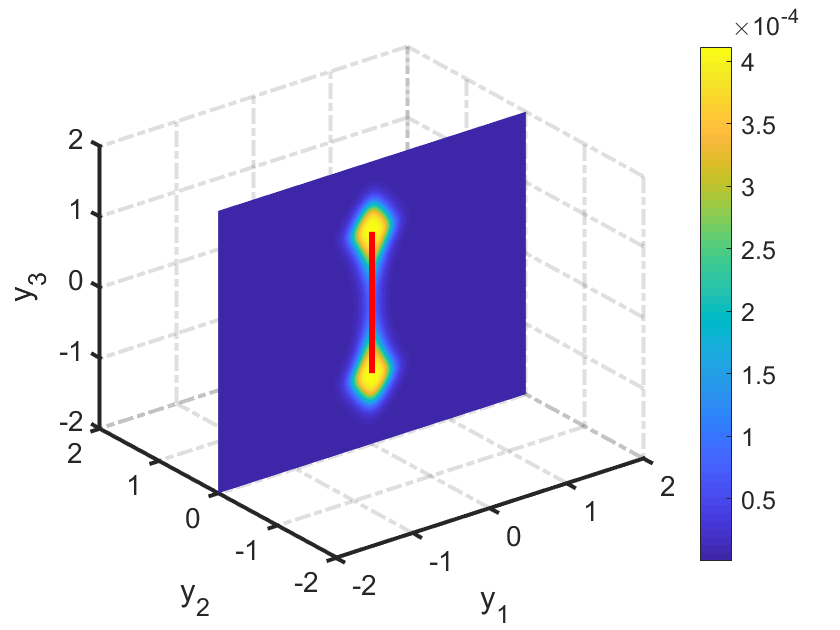}
		
	}
   \subfigure[$M=6$]{
		\includegraphics[scale=0.22]{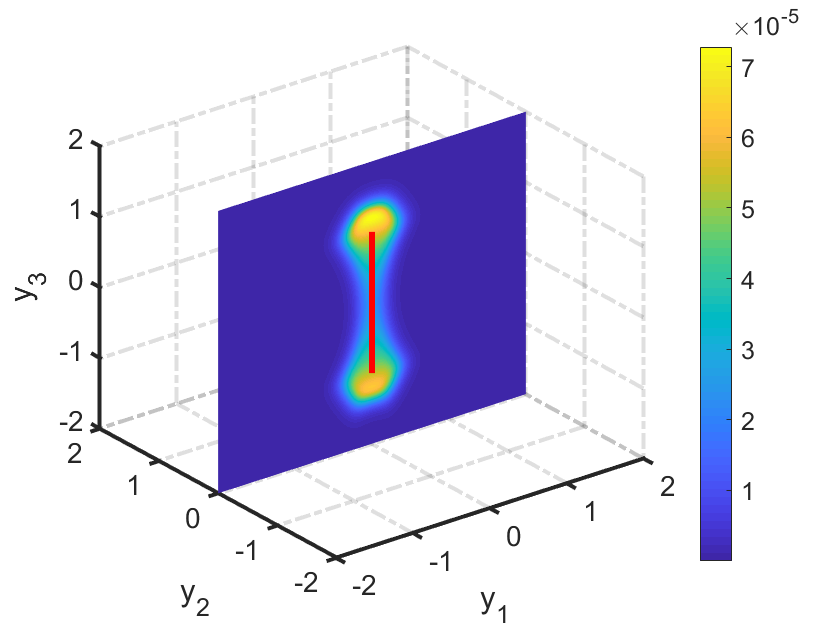}
		
	}
    \subfigure[$M=12$]{
		\includegraphics[scale=0.22]{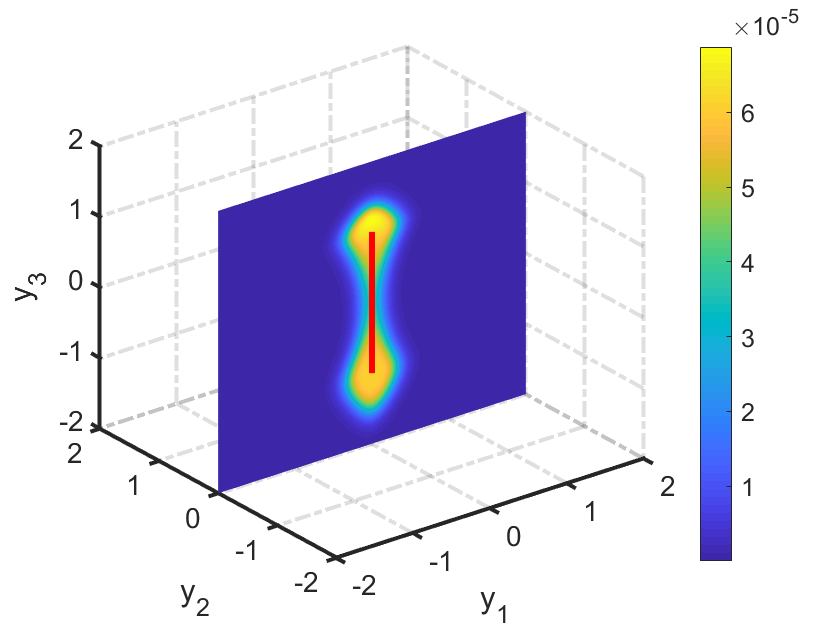}
		
	}
\caption{Reconstruction from sparse observation points $x=(2 \sin \varphi \cos \theta, 2\sin \varphi \sin \theta,$  $2 \cos \varphi )$  with $\theta\in[0, 2\pi)$   and $\varphi\in[0,\pi]$ for a straight line segment $a(t)=(0,0,t-2)$ where $t\in[1,3]$. Here $M$ denotes the number of the observation points. (a) $\theta=(j-1)*2\pi/4$, $j=1,\cdots,4$ and $\varphi=\pi/2$;  (b)  $\varphi=0$; $\theta=(j-1)*2\pi/4$, $j=1,\cdots,4$ and $\varphi=\pi/2$;  $\varphi=\pi$;  (c) $\theta=(2j-1)*\pi/4$, $j=1,\cdots,4$, $\varphi=\pi/3$; $\theta=(j-1)*2\pi/4$, $j=1,\cdots,4$, $\varphi=\pi/2$;  $\theta=(4j-3)*\pi/8$, $j=1,\cdots,4$, $\varphi=2\pi/3$. We take $k_{min}=1$ and $k_{max}=6$.}\label{fig:line-m-2}
\end{figure}

\begin{figure}[H]
	\centering
    \subfigure[$M=4$ ]{
		\includegraphics[scale=0.3]{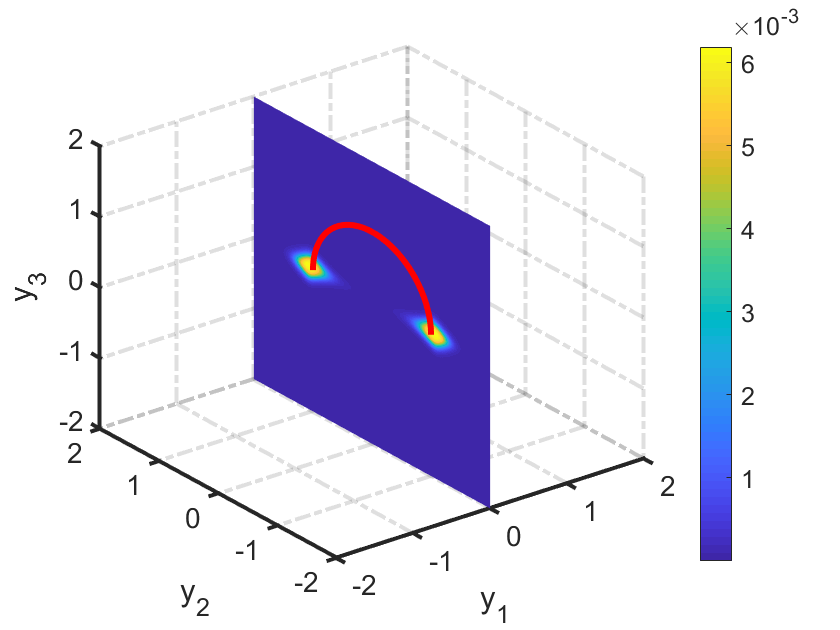}
		
	}
    \subfigure[$M=6$ ]{
		\includegraphics[scale=0.3]{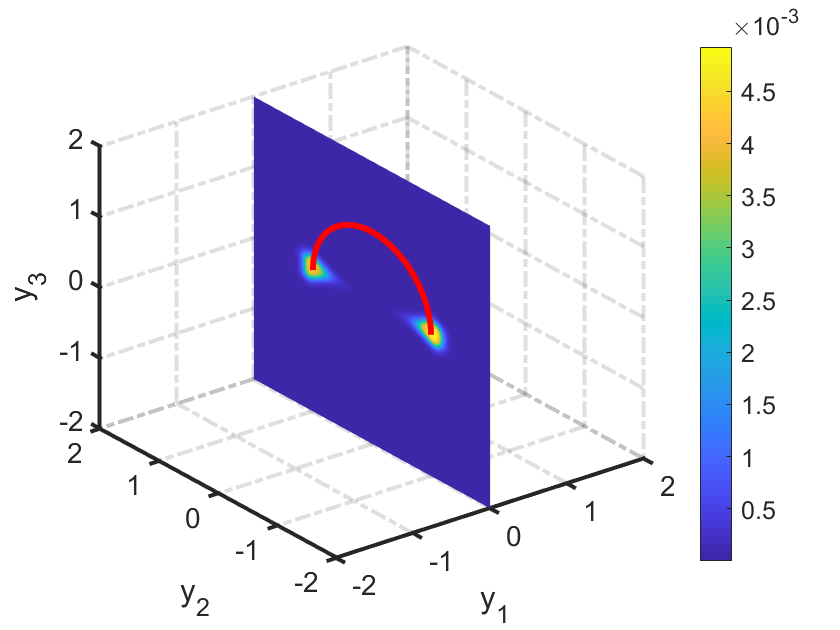}
		
	}
\caption{Reconstruction from sparse observation points $x=(2 \sin \varphi \cos \theta, 2\sin \varphi \sin \theta,$  $2 \cos \varphi )$  with $\theta\in[0, 2\pi)$   and $\varphi\in[0,\pi]$ for a straight line segment $a(t)=(0,\cos t, \sin t)$ where $t\in[0,\pi]$. Here $M$ denotes the number of the observation points. (a) $\theta=(j-1)*2\pi/4$, $j=1,\cdots,4$ and $\varphi=\pi/2$;  (b)  $\theta=(j-1)*2\pi/4$, $j=1,\cdots,4$ and $\varphi=\pi/2$; $\varphi=0$; $\varphi=\pi$.} \label{fig:circle-m-1}
\end{figure}

\begin{figure}[H]
	\centering

    \subfigure[$M=4$]{
		\includegraphics[scale=0.3]{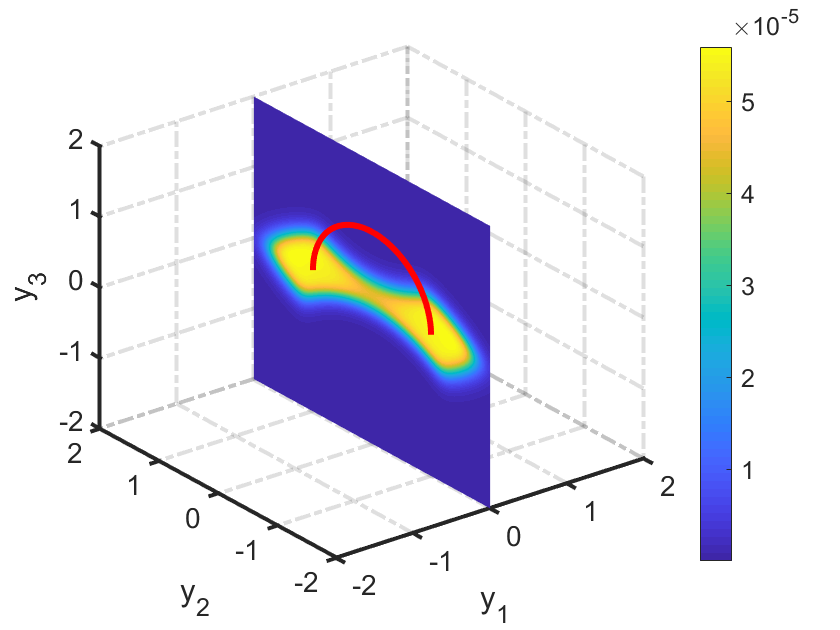}
		
	}
   \subfigure[$M=6$]{
		\includegraphics[scale=0.3]{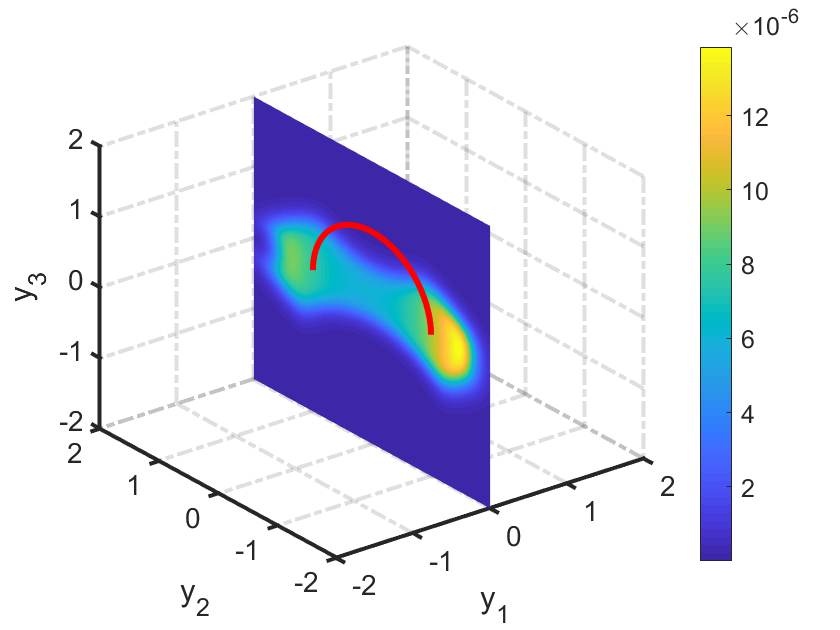}
		
	}
\caption{Reconstruction from sparse observation points $x=(2 \sin \varphi \cos \theta, 2\sin \varphi \sin \theta,$  $2 \cos \varphi )$  with $\theta\in[0, 2\pi)$   and $\varphi\in[0,\pi]$ for a straight line segment $a(t)=(0,\cos t, \sin t)$ where $t\in[0,\pi]$. Here $M$ denotes the number of the observation points. (a) $\theta=(j-1)*2\pi/4$, $j=1,\cdots,4$ and $\varphi=\pi/2$;  (b)  $\theta=(j-1)*2\pi/4$, $j=1,\cdots,4$ and $\varphi=\pi/2$; $\varphi=0$; $\varphi=\pi$. We take $k_{min}=1$ and $k_{max}=6$.}\label{fig:circle-m-2}
\end{figure}

Next,  assuming that all the selected observation points may contain both  observable and non-observable points. In the presented numerical examples below, we set the threshold value as $M'=10^5$.
Figs.\ref{fig:line-m-1} and \ref{fig:line-m-2} present the reconstructed trajectory for orbit functions $a(t)=(0,0,t-2)$ with $t\in[1,3]$ using different sparse observation points and different frequency bandwidths. Although sparse observation points data $M=4, 6, 12$ are used, the trajectory cannot be fully determined from Figs.\ref{fig:line-m-1} and \ref{fig:line-m-2}. The reason is that there consistently exist observation points $x_j$ that satisfy $A_{\Gamma}^{(x_j)}\subset \Lambda_{\Gamma}^{(x_j)}$, just as  what is depicted in Fig.\ref{fig:line2}. Further more, for these specific observation points, the reconstructed annular $A_{\Gamma}^{(x_j)}$ might be excessively slim. Consequently, the intersections of the annular $A_\Gamma^{(\hat{x}_j)}$ can only reconstruct the starting and ending points of the moving point source's trajectory. Figs.\ref{fig:circle-m-1} and \ref{fig:circle-m-2} further illustrate visualizations of the reconstructed trajectory pertaining to orbit functions $a(t)=(0,\cos t, \sin t)$ with $t\in[0,\pi]$, employing sparse observation points. These visualizations enable the determination of the starting and ending points of the arc-shaped trajectory.  Due to the lack of low-frequency data $\{u(x,\omega), k\in(0,1)\}$,  the inversion results in Figs.\ref{fig:line-m-2} and \ref{fig:circle-m-2} with $\omega\in(1,6)$ are not as accurate as those in Figs.\ref{fig:line-m-1} and \ref{fig:circle-m-1} with $\omega\in(0,6)$.

\begin{figure}[H]
	\centering
		\includegraphics[scale=0.5]{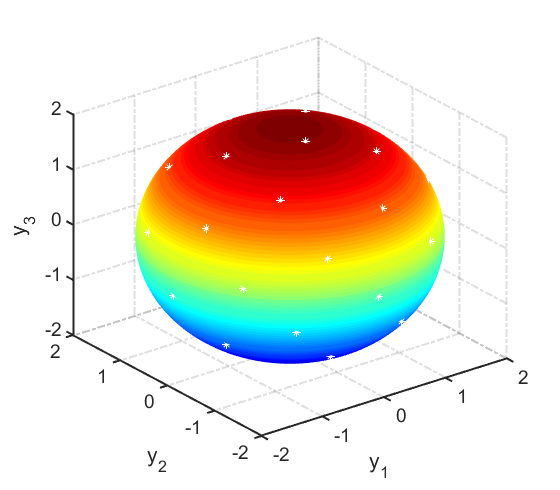}
\caption{40 uniformly distributed points on a sphere with a radius of 2.}\label{fig:m}
\end{figure}

\begin{figure}[H]
	\centering
    \subfigure[$M=20$]{
		\includegraphics[scale=0.22]{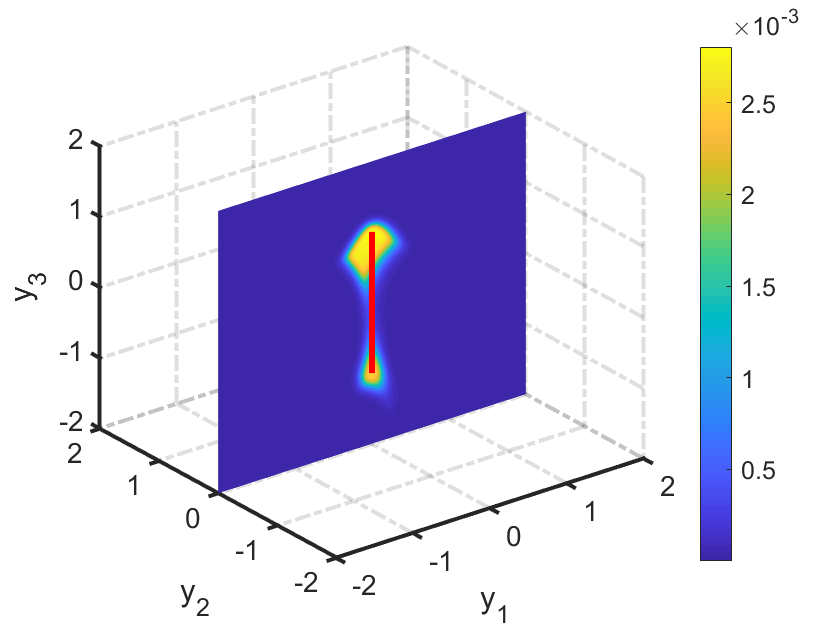}
		
	}
    \subfigure[$M=30$]{
		\includegraphics[scale=0.22]{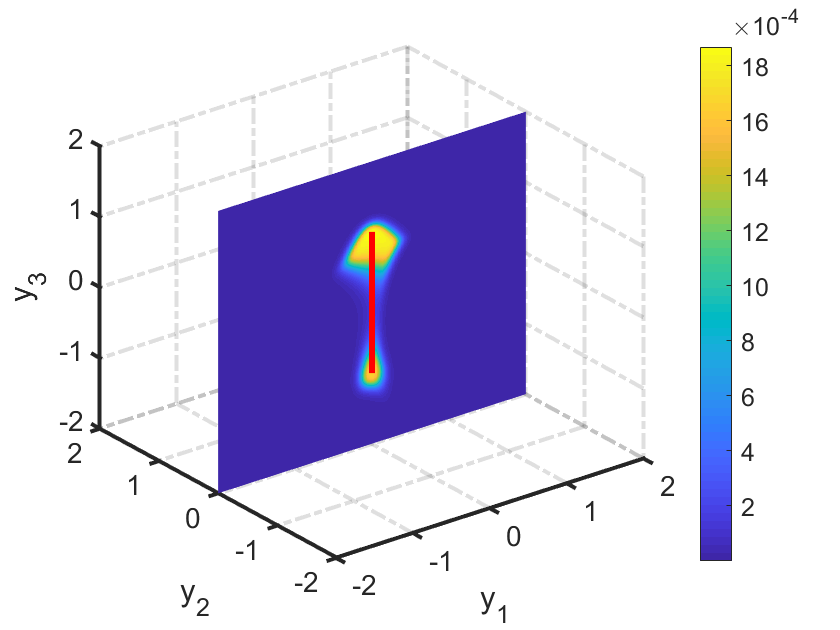}
		
	}
    \subfigure[$M=40$]{
		\includegraphics[scale=0.22]{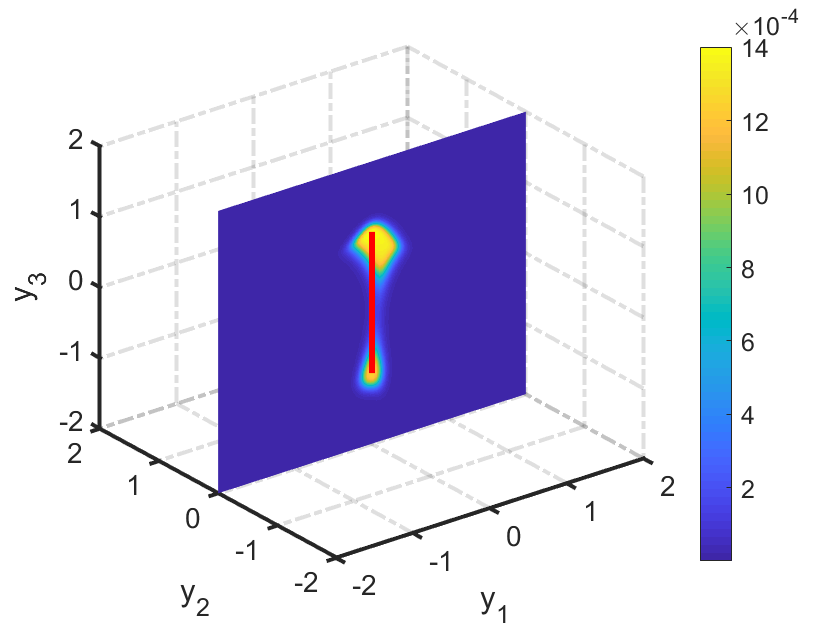}
		
	}
\caption{Reconstruction from sparse uniformly distributed observation points on a sphere with  a radius of 2 for a straight line segment $a(t)=(0,0,t-2)$ where $t\in[1,3]$. Here we take different numbers of observation points $M$. }\label{fig:m-1}
\end{figure}

\begin{figure}[H]
	\centering
    \subfigure[$M=20$]{
		\includegraphics[scale=0.22]{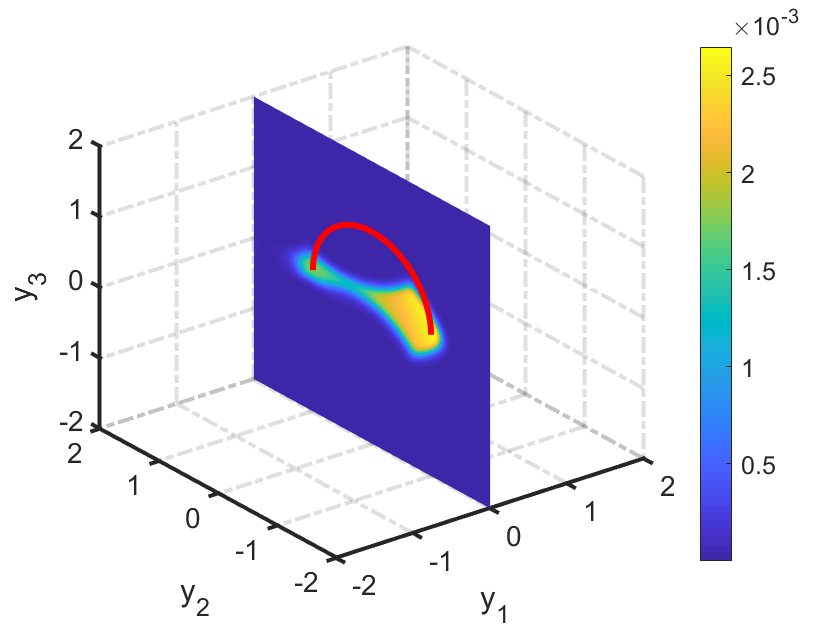}
		
	}
    \subfigure[$M=30$]{
		\includegraphics[scale=0.22]{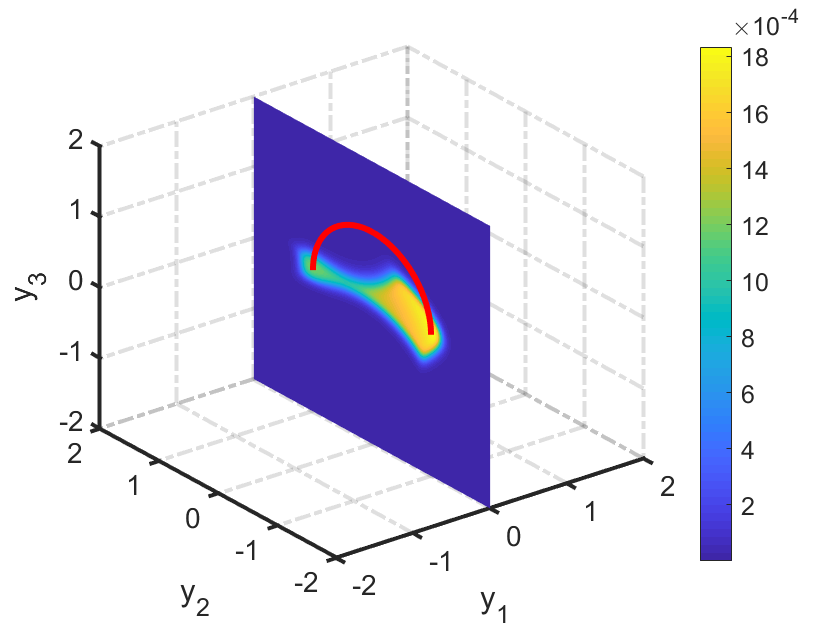}
		
	}
   \subfigure[$M=40$]{
		\includegraphics[scale=0.22]{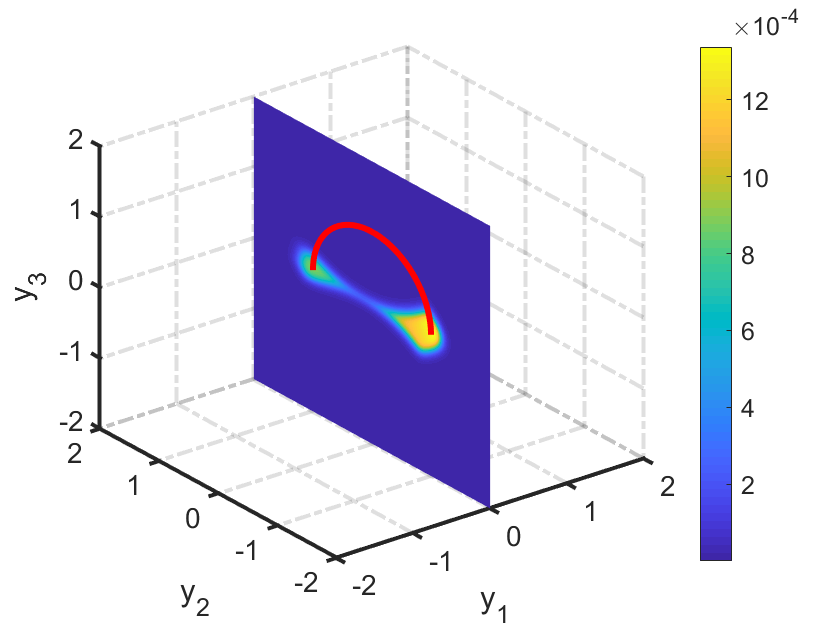}
		
	}
\caption{Reconstruction from sparse uniformly distributed observation points on a sphere with  a radius of 2 for an arc  $a(t)=(0,\cos t,\sin t)$ where $t\in[0,\pi]$. Here we take different numbers of observation points $M$. }\label{fig:m-2}
\end{figure}

Finally, our goal is to reconstruct the trajectories of a straight line segment and an arc using uniformly distributed observation points situated on a sphere with a radius of $2$. We explore sets of $20$, $30$, and $40$ points for this purpose. The positioning of the 40 uniformly distributed points is presented in Fig.\ref{fig:m}.
 With the increasing number of observation points, Figs.\ref{fig:m-1} and \ref{fig:m-2} illustrate that  the starting and ending points of the trajectories can be obtained. This limitation stems from the characteristics of the observation points within the observable set, where $A_{\Gamma}^{(x)}\subset \Lambda_{\Gamma}^{(x)}$.

\subsection{Noisy test}

 We evaluate sensitivity with respect to the noisy data by selecting Example 1, which involves a line segment recovery. The near-field data are corrupted with Gaussian noise, as shown below:
 \ben
 u_\delta(x,\omega)\coloneqq \real [u(x,\omega)]\;\big(1+\delta\, \gamma_1\big) + i {\rm Im} [u(x,\omega)]\;\big(1+\delta\, \gamma_2\big)
 \enn
 where $\delta>0$ represents the noise level  and $\gamma_j\in[-1,1]$ $(j=1,2)$ denote Gaussian random variables.

\begin{figure}[H]
	\centering
    \subfigure[$\delta=0\%$]{
		\includegraphics[scale=0.3]{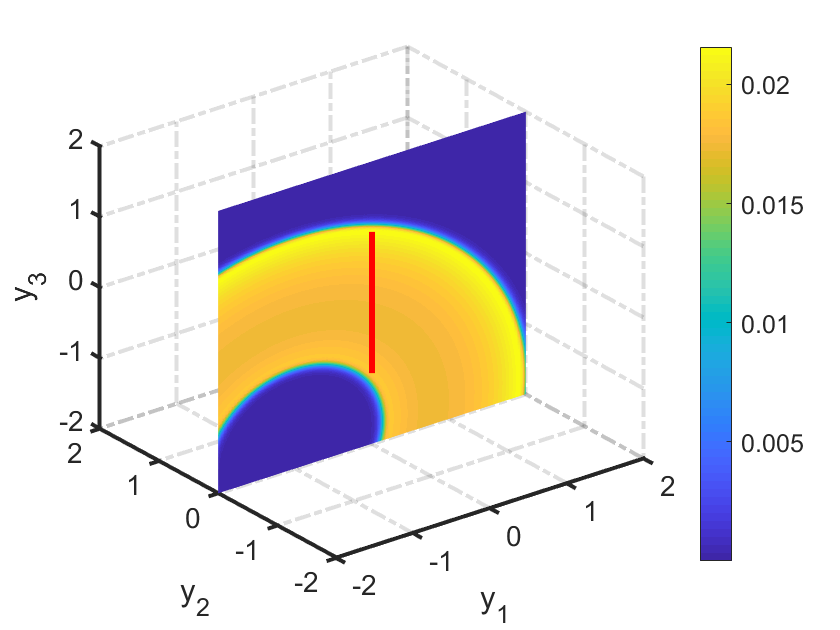}
		
	}
    \subfigure[$\delta=5\%$]{
		\includegraphics[scale=0.3]{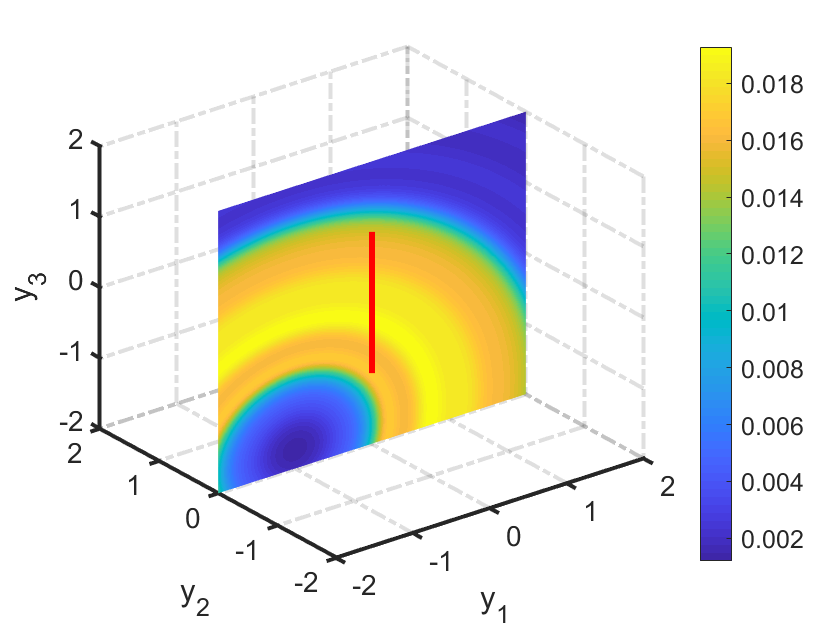}
		
	}
   \subfigure[$\delta=10\%$]{
		\includegraphics[scale=0.3]{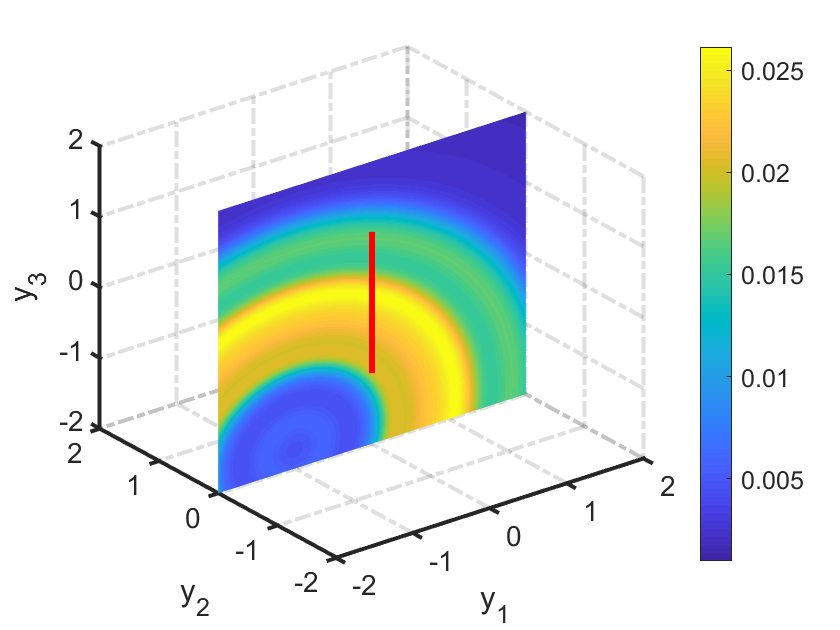}
		
	}
\subfigure[$\delta=20\%$]{
		\includegraphics[scale=0.3]{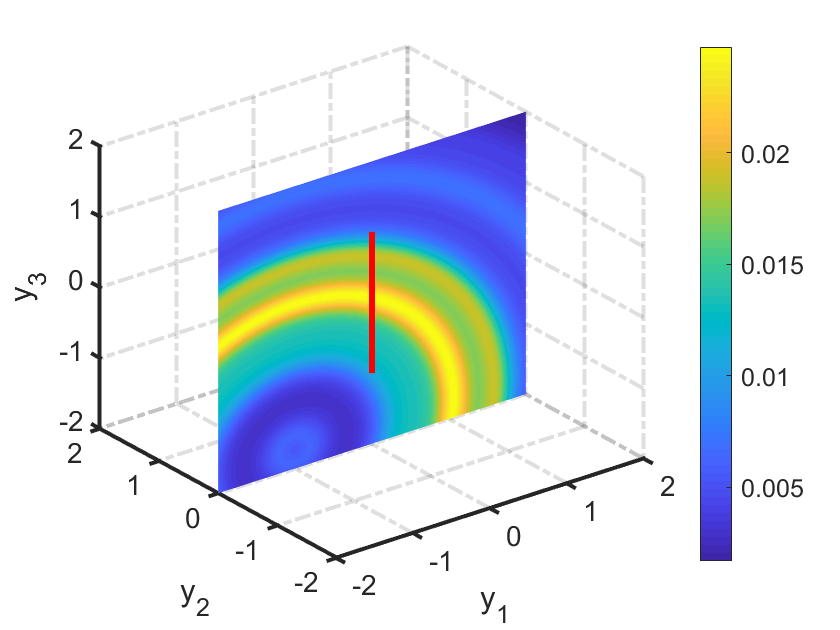}
		
	}
\caption{ Reconstruction of a straight line segment $a(t) = (0,0,t-2)$, $t\in [1,3]$ from noisy data with different levels  $\delta$ measured at a single observable point  $x=(2 \sin \varphi \cos \theta, 2\sin \varphi \sin \theta,$  $2 \cos \varphi )$ with $\theta=\pi$ and $\varphi=5\pi/6$.} \label{fig:noise}
\end{figure}

\begin{figure}[H]
	\centering
    \subfigure[$\delta=0\%$]{
		\includegraphics[scale=0.3]{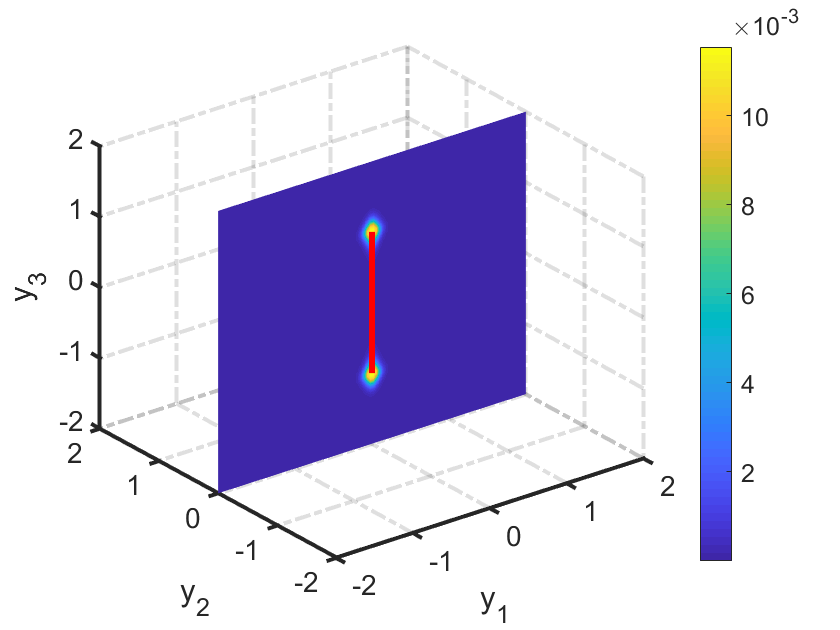}
		
	}
    \subfigure[$\delta=5\%$]{
		\includegraphics[scale=0.3]{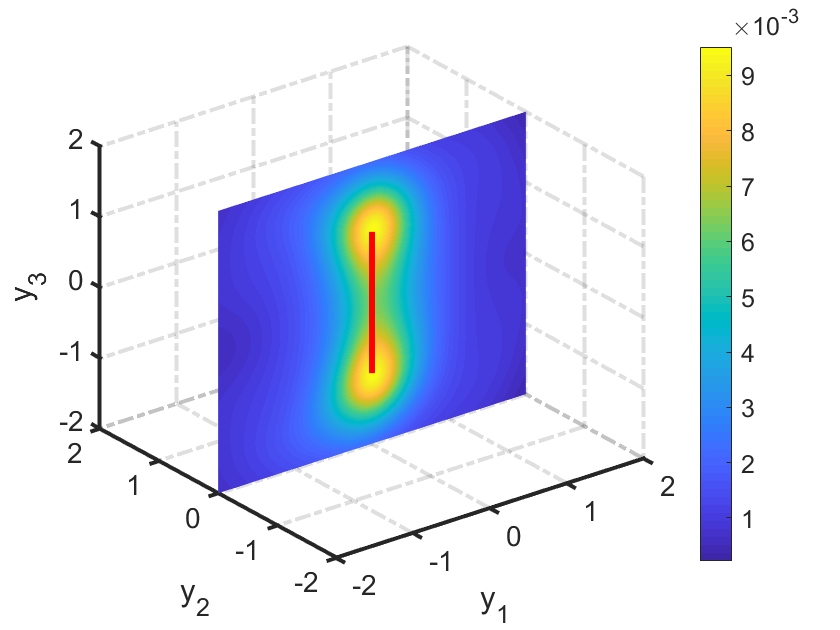}
		
	}
   \subfigure[$\delta=10\%$]{
		\includegraphics[scale=0.3]{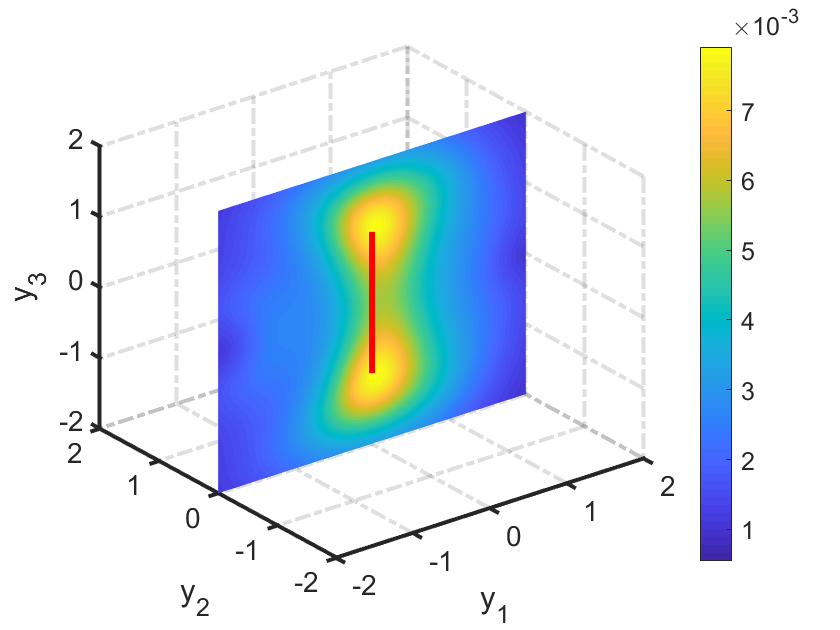}
		
	}
\subfigure[$\delta=20\%$]{
		\includegraphics[scale=0.3]{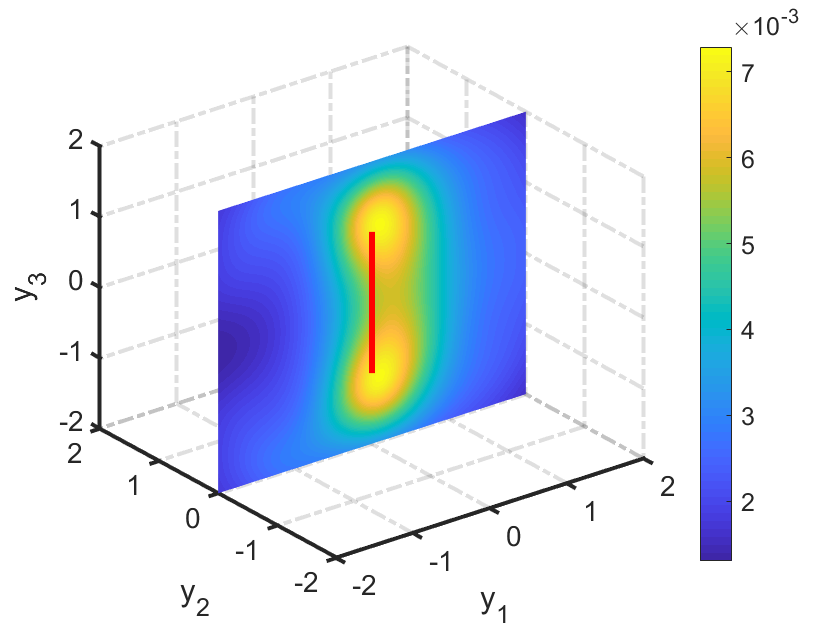}
		
	}
\caption{Reconstruction of a straight line segment $a(t) = (0,0,t-2)$, $t\in [1,3]$ from noisy data with different levels  $\delta$ measured at sparse observable points  $x=(2 \sin \varphi \cos \theta, 2\sin \varphi \sin \theta,$  $2 \cos \varphi )$ with $\theta=0, \pi/2, \pi, 3\pi/2$ and $\varphi=\pi/2$.}\label{fig:noisem}
\end{figure}

To accomplish this test, we assigned $\delta=0\%, 5\%, 10\%, 20\%$ and plot the indicator functions in Figs.\ref{fig:noise} and \ref{fig:noisem} using one and sparse observation points, respectively.
The images are clearly getting distorted at higher noise levels, but the starting and the ending points of the trajectory of the moving source using the data measured at sparse points can still be captured.

\section*{Acknowledgements}

G. Hu is partially supported by the National Natural Science Foundation of China (No. 12071236) and the Fundamental Research Funds for Central Universities in China (No. 63233071).

\end{document}